\documentclass[11pt]{article}
\usepackage{amssymb}
\usepackage{latexsym}
\usepackage{amsmath}
\usepackage{epsfig}
\usepackage{multicol}
\newcommand{\IP}{\mathbb P}

\newcommand{\IR}{\mathbb R}
\newcommand{\IZ}{\mathbb Z}
\newcommand{\IN}{\mathbb N}
\newcommand{\IT}{\mathbb T}

\newcommand \N {\mathbb{N}}
\newcommand \Z {\mathbb{Z}}
\newcommand \R {\mathbb{R}}
\newcommand \T {\mathbb{T}}
\newcommand \prob {\mathbb{P}}
\newcommand \E {\mathbb{E}}
\newcommand \ov {\overline}
\newcommand \e {\varepsilon}
\newcommand \integ[1]{\lfloor #1 \rfloor}
\newcommand \ind[1] {\mathbf{1}_{\{ #1 \}}}
\newcommand \A {\mathcal{A}}
\newcommand \GA {\Gamma_{\mathcal{A}}}
\newcommand \vp {\varpi}
\newcommand \G {\mathcal{G}}

\newtheorem{thm}{Theorem}[section]
\newtheorem{propn}[thm]{Proposition}
\newtheorem{lemma}[thm]{Lemma}

\newtheorem{defn}[thm]{Definition}
\newtheorem{remark}[thm]{Remark}
\newtheorem{notn}[thm]{Notation}

\newtheorem{assumption}[thm]{Assumption}

\begin{document}
\title{A new model for evolution in a spatial continuum}

\bigskip
\author{N.H. Barton\thanks{NHB supported in part by EPSRC Grant EP/E066070/1}\\Institute of Science and Technology
\\Am Campus I\\A-3400 Klosterneuberg, Austria\\
\texttt{Nick.Barton@ist-austria.ac.at}
\\~\\ A.M. Etheridge\thanks{AME supported in part by EPSRC Grant EP/E065945/1}
\\Department of Statistics\\University of Oxford\\1 South Parks Road\\Oxford OX1 3TG, UK \\
\texttt{etheridg@stats.ox.ac.uk}
\\~\\ A. V\'eber\\D\'epartement de math\'ematiques\\
Universit\'e Paris 11\\91405 Orsay Cedex, France\\
\texttt{amandine.veber@math.u-psud.fr}}
\date{}
 \maketitle

\begin{abstract}
We investigate a new model for populations evolving in a spatial continuum.  This model can be thought of as a spatial version of the $\Lambda$-Fleming-Viot process. It explicitly incorporates both small scale reproduction events and large scale extinction-recolonisation events. The lineages ancestral to a sample from a population evolving according to this model can be described in terms of a spatial version of the $\Lambda$-coalescent. Using a technique of Evans~(1997), we prove existence and uniqueness in law for the model. We then investigate the asymptotic behaviour of the genealogy of a finite number of individuals sampled uniformly at random (or more generally `far enough apart') from a two-dimensional torus of sidelength $L$ as $L\rightarrow\infty$. Under appropriate conditions (and on a suitable timescale) we can obtain as limiting genealogical processes a Kingman coalescent, a more general $\Lambda$-coalescent or a system of coalescing Brownian motions (with a non-local coalescence mechanism).



\bigskip
\noindent Submitted to EJP on April 3, 2009, final version accepted
February 5, 2010.
\end{abstract}

\baselineskip=6mm

\section{Introduction}
\label{intro}

In 1982, Kingman introduced a process called the {\em coalescent}. \nocite{kingman:1982} This process provides a simple and elegant description of the genealogical (family) relationships amongst a set of neutral genes in a randomly mating (biologists would say {\em panmictic}) population of constant size. Since that time, spurred on by the flood of DNA sequence data, considerable effort has been spent extending Kingman's coalescent to incorporate things like varying population size, natural selection and spatial (and genetic) structure of populations. Analytic results for these coalescent models can be very hard to obtain, but it is relatively easy, at least in principle, to simulate them and so they have become fundamental tools in sequence analysis. However, models of spatial structure have largely concentrated on subdivided populations and a satisfactory model for the ancestry of a population evolving in a two-dimensional spatial continuum has remained elusive. Our aim in this paper is to present the first rigorous investigation of a new model that addresses some of the difficulties of existing models for spatially extended populations while retaining some analytic tractability. The rest of this introduction is devoted to placing this research in context. The reader eager to skip straight to the model and a precise statement of our main results should proceed directly to Section \ref{model}.

Our concern here is with the extension of the coalescent to spatially structured populations. In this setting it is customary to assume that the population is subdivided into {\em demes} of (large) constant size, each situated at a vertex of a graph $G$, and model the genealogical trees using the {\em structured} coalescent. As we trace backwards in time, within each deme the ancestral lineages follow Kingman's coalescent, that is each pair of lineages merges (or {\em coalesces}) into a single lineage at a constant rate, but in addition lineages can migrate between demes according to a random walk on the graph $G$. The genealogical trees obtained in this way coincide with those for a population whose forwards in time dynamics are given by Kimura's stepping stone model (Kimura~1953) \nocite{kimura:1953} or, as a special case, if $G$ is a complete graph, by Wright's island model (Wright~1931). \nocite{wright:1931}

The stepping stone model is most easily described when the population consists of individuals of just two types, $a$ and $A$ say.  It can be extended to incorporate selection, but let us suppose for simplicity that these types are selectively neutral. Labelling the vertices of the graph $G$ by the elements of the (finite or countable) set $I$ and writing $p_i$ for the proportion of individuals in deme $i$ of type $a$, say, we have
\begin{equation}
\label{stepstone model} dp_i(t)=\sum_{j\in I}m_{ji}\left(p_j(t)-p_i(t)\right)dt+\sqrt{\gamma p_i(t)\left(1-p_i(t)\right)}dW_i(t)
\end{equation}
where $\{W_i(t); t\geq 0\}_{i\in I}$ is a collection of independent Wiener processes, $\gamma$ is a positive constant and $\{m_{ij}\}_{i,j\in I}$ specifies the rates of a continuous time random walk on $G$.  The graph $G$, chosen to caricature the spatial structure of the population, is typically taken to be $\IZ^2$ (or its intersection with a two-dimensional torus) and then one sets $m_{ij}=\kappa\mathbf{1}_{\{\|i-j\|=1\}}$, corresponding to simple random walk.

Although the stepping stone model is widely accepted as a model for structured populations, in reality, many populations are not subdivided, but instead are distributed across a spatial continuum. Wright~(1943) and Mal\'ecot~(1948) \nocite{wright:1943} \nocite{malecot:1948} derived expressions for the probability of identity of two individuals sampled from a population dispersed in a two-dimensional continuum by assuming on the one hand that genes reproduce and disperse independently of one another, and on the other hand that they are scattered in a stationary Poisson distribution. However, these assumptions are incompatible (Felsenstein~1975, Sawyer \& Fleischmann~1979). \nocite{felsenstein:1975} \nocite{sawyer/fleischmann:1979} The assumption of independent reproduction will result in `clumping' of the population and some local regulation will be required to control the local population density.

A closely related approach is to assume that the genealogical trees can be constructed from Brownian motions which coalesce at an instantaneous rate given by a function of their separation. The position of the common ancestor is typically taken to be a Gaussian centred on the midpoint between the two lineages immediately before the coalescence event (although other distributions are of course possible). However, the coalescent obtained in this way does not exhibit {\em sampling consistency}. That is, if we construct the genealogical tree corresponding to a sample of size $n$ and then examine the induced genealogical tree for a randomly chosen subsample of size $k<n$, this will not have the same distribution as the tree we obtain by constructing a system of coalescing lineages directly from the subsample. The reason is that whenever one of the lineages in the subsample is involved in a coalescence event in the full tree it will jump. Furthermore, just as in Mal\'ecot's setting, there is no corresponding {\em forwards} in time model for the evolution of the population.

Barton et al.~(2002) extend the formulae of Wright and Mal\'ecot to \nocite{barton/depaulis/etheridge:2002} population models which incorporate local structure. The probability of identity is obtained from a recursion over timeslices of length $\Delta t$. Two related assumptions are made. First, the ancestral lineages of genes that are sufficiently well separated are assumed to follow independent Brownian motions (with an effective dispersal rate which will in general differ from the forwards in time dispersal rate) and their chance of coancestry in the previous timeslice is negligible. Second, it must be possible to choose $\Delta t$ sufficiently large that the changes in the population over successive timeslices are uncorrelated. (For general $\Delta t$ this will not be the case. The movements of ancestral lineages in one time step may be correlated with their movements in previous steps if, for example, individuals tend to disperse away from temporarily crowded clusters.) Over all but very small scales, the resulting probability of identity can be written as a function of three parameters: the {\em effective dispersal rate}, the {\em neighbourhood size} and the {\em local scale}. However the usefulness of this result is limited due to a lack of explicit models for which the assumptions can be validated and the effective parameters established. Moreover, as explained in Barton et al.~(2002), although one can in principle extend the formula to approximate the distribution of genealogies amongst larger samples of well-separated genes, additional assumptions need to be made if such genealogies are to be dominated by pairwise coalescence. If several genes are sampled from one location and neighbourhood size is small then multiple coalescence (by which we mean simultaneous coalescence of {\em three} or more lineages) could become significant.

Multiple merger coalescents have received considerable attention from mathematicians over the last decade. Pitman~(1999) and Sagitov~(1999) introduced what we now call \nocite{pitman:1999} \nocite{sagitov:1999} {\em $\Lambda$-coalescents}, in which more than two ancestral lineages can coalesce in a single event, but {\em simultaneous} coalescence events are not allowed. Like Kingman's coalescent, these processes take their values among partitions of $\IN$ and their laws can be prescribed by specifying the restriction to partitions of $\{1,2,\ldots ,n\}$ for each $n\in\IN$. For our purposes, the $\Lambda$-coalescent describes the ancestry of a population whose individuals are labelled by $\IN$. Each block in the partition at time $t$ corresponds to a single ancestor at time $t$ before the present, with the elements of the block being the descendants of that ancestor. Tracing backwards in time, the evolution of the $\Lambda$-coalescent is as follows: if there are currently $p$ ancestral lineages, then each transition involving $j$ of the blocks merging into one happens at rate
\begin{equation}
\label{betas in the coalescent} \beta_{p,j}^{\Lambda}=\int_{[0,1]}u^{j-2}(1-u)^{p-j}\Lambda(du),
\end{equation}
and these are the only possible transitions. Here, $\Lambda$ is a finite measure on $[0,1]$. Kingman's coalescent corresponds to the special case $\Lambda=\delta_0$, the point mass at the origin.
\begin{remark}
More generally, one can consider processes with simultaneous multiple coalescence events. Such coalescents were obtained as the genealogies of suitably rescaled population models by M\"ohle \& Sagitov~(2001). Independently, Schweinsberg~(2000) obtained the same class of coalescents and characterised the possible rates of mergers in terms of a single measure $\Xi$ on an infinite simplex. Coalescents which allow {\em simultaneous} multiple mergers are now generally referred to as {\em $\Xi$-coalescents}. \nocite{mohle/sagitov:2001} \nocite{schweinsberg:2000}
\end{remark}
Kingman's coalescent can be thought of as describing the genealogy of a random sample from a Fleming-Viot process. In the same way, a $\Lambda$-coalescent describes the genealogy of a random sample from a generalised Fleming-Viot process. This process takes its values among probability measures on $[0,1]$. We shall describe it in terms of its generator, $\mathcal R$ acting on functions of the form
$$
F(\rho)=\int f(x_1,\ldots ,x_p)\rho(dx_p)\ldots\rho(dx_1),
$$
where $p\in\IN$ and $f:[0,1]^p\rightarrow\IR$ is measurable and bounded. First we need some notation. If $x=(x_1,\ldots ,x_p)\in [0,1]^p$ and $J\subseteq \{1,\ldots ,p\}$ we write
$$
x_i^J=x_{\min J}\mbox{ if }i\in J, \mbox{ and } x_i^J=x_i \mbox{ if }i\notin J, \ \ i=1,\ldots ,p.
$$
Then for $\Lambda$ a finite measure on $[0,1]$, a $\Lambda$-Fleming-Viot process has generator
$$
{\mathcal R}F(\rho)=\sum_{J\subseteq\{1,\ldots ,p\}, |J|\geq 2} \beta_{p,|J|}^{\Lambda}\int\left(f(x_1^J,\ldots ,x_p^J)-f(x_1,\ldots ,x_p)\right) \rho(dx_p)\ldots \rho(dx_1),
$$
where $\beta_{p,j}^{\Lambda}$ is defined in Equation~(\ref{betas in the coalescent}). When $\Lambda(\{0\})=0$, this can also be written
$$
{\mathcal R}F(\rho)=\int_{(0,1]}\int_{[0,1]} \Big(F\big((1-u)\rho+u\delta_k\big)- F(\rho)\Big) \rho(dk)u^{-2}\Lambda(du).
$$
(When $\Lambda(\{0\})>0$, one must add a second term corresponding to a classical Fleming-Viot process and somehow dual to the Kingman part of the $\Lambda$-coalescent.) In this case, an intuitive way to think about the process is to consider a Poisson point process on $\IR_+\times (0,1]$ with intensity measure $dt\otimes u^{-2}\Lambda (du)$, which picks jump times and sizes for $\rho(t)$.  At a jump time $t$ with corresponding jump size $u$, a type $k$ is chosen according to $\rho(t-)$, an atom of mass $u$ is inserted at $k$ and $\rho(t-)$ is scaled down by $(1-u)$ so that the total mass remains equal to one, i.e.,
\begin{equation}
\label{lambda fv} \rho(t)=(1-u)\rho(t-)+u\delta_k.
\end{equation}
The duality between $\Lambda$-coalescents and $\Lambda$-Fleming-Viot processes was first proved by Bertoin \& Le Gall~(2003). Their approach uses a correspondence between the $\Lambda$-coalescents and stochastic flows of bridges. The duality can also be understood via the Donnelly \& Kurtz~(1999) `modified lookdown construction' and indeed is implicit there. An explicit explanation can be found in Birkner et al.~(2005). \nocite{bertoin/legall:2003} \nocite{birkner/blath/capaldo/etal:2005} \nocite{donnelly/kurtz:1999}

In recent work (described briefly in Etheridge~2008), \nocite{etheridge:2008} Barton \& Etheridge have proposed a new class of consistent forwards and backwards in time models for the evolution of allele frequencies in a population distributed in a two-dimensional (or indeed $d$-dimensional) spatial continuum which, in the simplest setting, can be thought of as spatial versions of the $\Lambda$-Fleming-Viot and $\Lambda$-coalescent models (although we emphasize that these are not the same as the spatial $\Lambda$-coalescents considered by Limic \& Sturm~2006). They \nocite{limic/sturm:2006} share many of the advantages of the classical models for spatially structured populations while overcoming at least some of the disadvantages. The idea is simple. Just as in the $\Lambda$-Fleming-Viot process, reproduction events are determined by a Poisson point process but now, in addition to specifying a time and a value $u$, this process prescribes a region of space which will be affected by the event. In what follows, the region will be a ball with random centre and radius. Within that region the effect is entirely analogous to Equation~(\ref{lambda fv}).

This approach differs from existing spatial models in three key ways. First, density dependent reproduction is achieved by basing reproduction events on neighbourhoods (whose locations are determined by the Poisson point process), rather than on individuals. Second, the offspring of a single individual can form a significant proportion of the population in a neighbourhood about the parent, capturing the essentially finite nature of the local population size. Third, large scale extinction-recolonisation events are explicitly incorporated. This reflects the large scale fluctuations experienced by real populations in which the movement and reproductive success of many individuals are correlated. For example, climate change has caused extreme extinction and recolonisation events that dominate the demographic history of humans and other species (e.g. Eller et al.~2004). \nocite{eller/hawks/relethford:2004}

The spatial $\Lambda$-Fleming-Viot process, like its classical counterpart, can be obtained as a limit of individual based models. Those prelimiting models are discussed in Berestycki et al.~(2009). \nocite{berestycki/etheridge/hutzenthaler:2009} In the (backwards in time) spatial $\Lambda$-coalescent, ancestral lineages move around according to dependent L\'evy processes (in fact they will be compound Poisson processes), jumping whenever they are affected by a reproduction event. Two or more lineages can coalesce if they are all affected by the same reproduction event.

Our first aim here is to provide a precise mathematical description of the spatial $\Lambda$-Fleming-Viot process and the corresponding spatial $\Lambda$-coalescent model and address questions of existence and uniqueness. This is achieved through adapting the work of Evans~(1997). \nocite{evans:1997} The idea is to first construct the dual (backwards in time) process of coalescing L\'evy processes corresponding to a finite sample from the population at time zero, and then to use a functional duality to define the forwards in time model. The principal difference between our setting and that of Evans is that, in his work, ancestral lineages evolve {\em independently} until they meet.

The system of coalescing L\'evy processes that describes the genealogy of a sample from the population, mirrors the system of coalescing random walks that plays the same r\^ole for the stepping stone model. For systems of coalescing walks a number of studies have investigated conditions under which, when viewed on an appropriate timescale, and for sufficiently well-separated samples, the effect of the geographical structure of the population can be summarised as a single `effective' parameter and the system of coalescing lineages converges to Kingman's coalescent. One of the first works along these lines is due to Cox~(1989), who considers random walks on a torus $\IT(L)\cap\IZ^d$ of sidelength $L$ with the walks coalescing instantly on meeting. This corresponds to taking $G=\IT(L)\cap\IZ^d$ and $\gamma=\infty$ in Equation~(\ref{stepstone model}). He shows that if one starts walks from any finite number $n\in\IN$ of points chosen independently and uniformly at random from $\IT(L)\cap\IZ^d$, then in suitable time units, as $L\rightarrow\infty$, the number of surviving lineages is determined by Kingman's coalescent. For two spatial dimensions, this analysis was extended by Cox \& Durrett~(2002) and Z\"ahle et al.~(2005) to random walks on $\IT(L)\cap\IZ^2$ with delayed coalescence (corresponding to $\gamma<\infty$). \nocite{cox/durrett:2002} \nocite{zahle/cox/durrett:2005} It is natural to ask whether similar results are true here. Our second aim then is to establish conditions under which the genealogy of a sample taken at random from a large torus will converge to a non-spatial coalescent. We shall concentrate on the most difficult, but also most biologically relevant, case of two spatial dimensions. If reproduction events only affect bounded neighbourhoods, then, not surprisingly, we recover a Kingman coalescent limit. However, we also consider the more general situation in which in addition to `small' events that affect only bounded neighbourhoods we allow `large' extinction-recolonisation events (see Section \ref{section result} for the precise setting). Unless these events affect a non-negligible proportion of the torus, on a suitable timescale, asymptotically we once again recover a Kingman coalescent. The timescale is determined by the relative rates of `large' and `small' events. However, if we have extinction-recolonisation events that affect regions with sidelength of order ${\mathcal O}(L)$, then, again depending on the relative rates of `large' and `small' events, we can obtain a more general (non-spatial) {\em $\Lambda$-coalescent} limit or a system of coalescing Brownian motions (where the coalescence is non-local).

The rest of the paper is laid out as follows. In Section \ref{model} we define the model. In Section \ref{section result}, we give a precise statement of the conditions under which we obtain convergence of the genealogy of a random sample from a (two-dimensional) torus of side $L$ as $L\rightarrow\infty$. The corresponding convergence results are Theorem~\ref{result alpha<1} and Theorem~\ref{result alpha=1}. In Section \ref{section existence} we establish existence of the process and prove uniqueness in law. In Section \ref{levy processes} we gather the necessary results on L\'evy processes in preparation for our proofs of Theorem~\ref{result alpha<1} and Theorem~\ref{result alpha=1} in Sections \ref{alpha<1} and \ref{alpha=1}. Finally, Appendices \ref{appendix 1} and \ref{appendix 2} contain the proofs of the technical lemmas stated in Sections \ref{levy processes} and \ref{alpha<1}.

\section{The model}
\label{model}

First we describe a prelimiting model. Individuals in our population are assumed to have a {\em type} taken from $[0,1]$ and a spatial position in a metric space $E$ that we shall usually take to be $\IR^2$ (or the torus $\IT(L)$ in $\IR^2$). Even though it will be clear that existence and uniqueness of the process holds in much greater generality, the model is primarily motivated by considerations for populations evolving in two-dimensional continua. The dynamics are driven by a Poisson point process $\Pi$ on $\IR_+\times \IR^2 \times (0,\infty)$ with intensity $dt\otimes dx\otimes \mu(dr)$. If $(t,x,r)\in\Pi$, the first component represents the time of a reproduction event. The event will affect only individuals in $B(x,r)$, the closed ball of centre $x$ and radius $r$. We require two more ingredients. The first, $m$, is a fixed positive constant which we shall refer to as the {\em intensity} of the model. Second, associated to each fixed radius $r>0$ there is a probability measure $\nu_r$ on $[0,1]$. In the sequel, we assume that the mapping $r\mapsto \nu_r$ is measurable with respect to $\mu$.

For definiteness, suppose that the population is initially distributed according to a spatially homogeneous Poisson process. The dynamics of our prelimiting model are described as follows. Suppose that $(t,x,r)\in\Pi$. Consider the population in $B(x,r)$ at time $t-$. If the ball is empty, then nothing happens. Otherwise, independently for each event:
\begin{enumerate}
\item{Select a `parent' uniformly at random from those individuals in $B(x,r)$ at time $t-$ and sample $u\in [0,1]$ at random according to $\nu_r$.}
\item{Each individual in $B(x,r)$, independently, dies with probability $u$, otherwise it is unaffected by the reproduction event.}
\item{Throw down offspring in the ball, with the same type as the selected parent (who may now be dead), according to an independent Poisson point process with intensity $\left.u\, m\,\mathrm{Leb}\right|_{B(x,r)}$ where $\mathrm{Leb}$ denotes Lebesgue measure.}
\end{enumerate}
We shall refer to these events as {\em reproduction events}, even though they are also used to model large-scale extinction-recolonisation events. Notice that recolonisation is modelled as being instantaneous even after a large scale extinction.
\begin{remark}
For simplicity we have described only a special version of the model in which, even when the reproduction event affects a large region, recolonisation is through a single founder. This guarantees that if we look at the genealogy of a sample from this population, although we may see more than two lineages coalescing in a single event, we do not see {\em simultaneous} mergers. More generally it would be natural to take a random number of colonists and then, on passing to the limit, the corresponding model would yield a spatial $\Xi$-coalescent.
\end{remark}

Any reproductive event has positive probability of leaving the corresponding region empty, but because the neighbourhoods determined by different reproduction events overlap, an empty region can subsequently become recolonised. Provided the measure $\mu(dr)$ decays sufficiently quickly as $r\rightarrow\infty$, Berestycki et al.~(2009) show that there is a critical value of $m$ above which the population, when started from a translation invariant initial condition, survives with probability one. The difficulty is that it is not easy to find an explicit expression for the distribution of the genealogical trees relating individuals in a sample from the population. Knowing that an ancestral lineage is in a given region of space gives us information about the rate at which that region was hit by reproduction events as we trace backwards in time. On the other hand, simulations reveal that this effect is rarely significant. Mathematically, we overcome this difficulty by considering a model in which the intensity $m$ is infinite, but we preserve some of the signature of a finite local population size by retaining the reproduction mechanism so that a non-trivial proportion of individuals in a neighbourhood are descended from a common ancestor. In particular, this will result in multiple coalescences of ancestral lineages.

Now let us describe the model that arises from letting $m\rightarrow \infty$. (That the prelimiting model really does converge to this limit will be proved elsewhere.) At each point $x\in\IR^2$, the model specifies a probability measure on type space which we shall write $\rho(t,x,\cdot)$, or sometimes for brevity $\rho_x$. The interpretation is that if we sample an individual from $x$, then its type will be determined by sampling from $\rho_x$. The reproduction mechanism mirrors that for our discrete time model:
\begin{defn}[Spatial $\Lambda$-Fleming-Viot process] \label{def slfv}The {\em spatial $\Lambda$-Fleming-Viot process}, $\{\rho(t,x,\cdot), x\in\IR^2, t\geq 0\}$ specifies a probability measure on the type space $[0,1]$ for every $t\geq 0$ and every $x\in\IR^2$. With the notation above, the dynamics of the process are as follows. At every point $(t,x,r)$ of the Poisson point process $\Pi$, we choose $u\in [0,1]$ independently according to the measure $\nu_r(du)$. We also select a point $z$ at random from $B(x,r)$ and a type $k$ at random according to $\rho(t-,z,\cdot)$. For all $y\in B(x,r)$,
$$
\rho(t,y,\cdot)=(1-u)\rho(t-,y,\cdot)+u\delta_k.
$$
Sites outside $B(x,r)$ are not affected, that is $\rho(t,y,\cdot)=\rho(t-,y,\cdot)$ for every $y\notin B(x,r)$.
\end{defn}
\begin{remark}
There are many variants of this model, some of which are outlined in Etheridge~(2008). \nocite{etheridge:2008} The model presented here should be regarded as fitting into a general framework in which the key feature is that reproduction events are driven by a Poisson point process determining their times and spatial locations, rather than on individuals. Barton et al.~(2009) investigate a version of the model in which, instead of replacing a portion $u$ of the population in a disc at the time of a reproduction event, the proportion of individuals affected decays (in a Gaussian distribution) with the distance from the `centre' $x$ of the event. \nocite{barton/kelleher/etheridge:2009} Whereas in the disc based approach in the prelimiting (individual based) model we had to suppress reproduction events that affected empty regions, this is not necessary in the Gaussian model. Moreover, (in contrast to the disc model) in that setting the prelimiting model has the Poisson point process in $\IR^2$ with constant intensity $m$ as a stationary distribution. Although the proofs would be rather involved, analogues of our results here should carry over to the Gaussian setting.
\end{remark}
Of course we must impose restrictions on the intensity measure if our process is to exist. To see what these should be, consider first the evolution of the probability measure $\rho(t,x,\cdot)$ defining the distribution of types at the point $x$. This measure experiences a jump of size $y\in A\subseteq (0,1]$ at rate
$$
\int_{(0,\infty)} \int_A \pi r^2 \nu_r(du)\mu(dr).
$$
By analogy with the $\Lambda$-Fleming-Viot process, we expect to require that
\begin{equation}
\label{condition 1} \Lambda(du)=\int_{(0,\infty)}u^2r^2 \nu_r(du)\mu(dr)
\end{equation}
defines a finite measure on $[0,1]$.  In fact, in the spatial setting we require a bit more. To see why, suppose that $\psi$ is a bounded measurable function on $[0,1]$ and consider the form that the infinitesimal generator of the process must take on test functions of the form $\langle \rho(x,dk),\psi(k)\rangle$ (with angle brackets denoting integration). Denoting the generator, if it exists, by $G$ we shall have
\begin{eqnarray*}
G(\langle \rho,\psi\rangle)&=& \int_{\IR^2}\int_{(0,\infty)}\int_{[0,1]}\int_{[0,1]} \frac{L_r(x,y)}{\pi r^2}\big(\langle (1-u)\rho(x,\cdot)+u\delta_k,\psi\rangle- \langle\rho(x,\cdot),\psi\rangle \big)\\
&&\phantom{AAAAAAAAAAAAAAAAAAAAAAAAAA}\rho(y,dk)\nu_r(du)\mu(dr)dy\\
&=& \int_{\IR^2}\int_{(0,\infty)}\int_{[0,1]} \frac{L_r(x,y)}{\pi r^2}\,u\big(\langle \rho(y,\cdot),\psi\rangle- \langle\rho(x,\cdot),\psi\rangle\big) \nu_r(du)\mu(dr)dy,
\end{eqnarray*}
where $L_r(x,y)$ denotes the volume of the set $B(x,r)\cap B(y,r)$. Notice in particular that $L_r(x,y)\leq\pi r^2\mathbf{1}_{\{|x-y|\leq 2r\}}$. In the non-spatial case, this term vanishes (set $y=x$), but here if we want the generator to be well-defined on these test functions we make the stronger
\begin{assumption}
\begin{equation}
\label{condition for convergence} \tilde{\Lambda} (du)=\int_{(0,\infty)}ur^2 \nu_r(du)\mu(dr)
\end{equation}
defines a finite measure on $[0,1]$.
\end{assumption}
Condition~(\ref{condition for convergence}) controls the jumps of $\rho$ at a single point. Since we are going to follow Evans~(1997) in constructing our process via the dual process of coalescing lineages ancestral to a sample from the population, we should check that such a process is well-defined. First we define the coalescent process more carefully.

In order to make sense of the genealogy of a sample at any time, we extend the Poisson point process $\Pi$ of reproduction events to the whole time line $(-\infty,+\infty)$.  We need some notation for (labelled) partitions.
\begin{notn}[Notation for partitions]
\label{notation for partitions}
\begin{enumerate}
\item{For each integer $n\geq 1$, let ${\mathcal P}_n$ denote the set of partitions of $\{1,\ldots,n\}$, and define a labelled partition of $\{1,\ldots,n\}$, with labels from a set $E$, to be a set of the form $\{(\pi_1,x_{\pi_1}),\ldots,(\pi_k,x_{\pi_k})\}$, where $\{\pi_1,\ldots,\pi_k\}\in {\mathcal P}_n$ and $(x_{\pi_1},\ldots,x_{\pi_k})\in E^k$. Let ${\mathcal P}_n^{\ell}$ be the set of all labelled partitions of $\{1,\ldots,n\}$. }
\item{For each $n\in \N$, let $\wp_n$ denote the partition of $\{1,\ldots,n\}$ into singletons. Moreover, if $E$ is the space of labels and $\mathbf{x}\equiv(x_1,\ldots,x_n)\in E^n$, let $\wp_n(\mathbf{x})$ denote the element $\{(\{1\},x_1),\ldots,(\{n\},x_n)\}$ of ${\mathcal P}_n^{\ell}$. }
\item{If $\pi \in {\mathcal P}_n^{\ell}$ for some $n\in \N$, then $\mathrm{bl}(\pi)$ will refer to the unlabelled partition of $\{1,\ldots,n\}$ induced by $\pi$ and if $a\in \mathrm{bl}(\pi)$, $x_a$ will be our notation for the label of $a$.}
\end{enumerate}
\end{notn}
Our genealogical process will be a labelled partition. As in classical representations of genealogical processes, a block of the partition at genealogical time $t\geq 0$ contains the indices of the initial lineages which share a common ancestor $t$ units of time in the past, and its label gives the current location of this ancestor in $E=\R^2$.

From the description of the forwards-in-time dynamics, the evolution of a sample of ancestral lineages represented by a labelled partition should be the following. We start with a finite collection of lineages at time $0$. At each point $(-t,x,r)\in \Pi$ (with $t\geq 0$ here, since genealogical time points towards the past), given that $u\in [0,1]$ is the result of the sampling according to $\nu_r$ each lineage present in the ball $B(x,r)$, independently, is affected (resp., is not affected) with probability $u$ (resp., $1-u$). A site $y$ is chosen uniformly in $B(x,r)$, and the blocks of all affected lineages merge into a single block labelled by $y$. The other blocks and their labels are not modified. We write $\{\A(t),\ t\geq 0\}$ for the Markov process of coalescing lineages described in this way. Its state space is $\bigcup_{n\geq 1}{\mathcal P}_n^{\ell}$. Note that $\A$ is constructed on the same probability space as that of the Poisson point process of reproduction events. Writing $\IP$ for the probability measure on that space, we abuse notation slightly by writing $\prob_A$ to indicate that $\A(0)=A$, $\prob_A$-a.s. Now let us verify that our Condition~(\ref{condition for convergence}) is sufficient to ensure that the process $\{\A(t),t\geq 0\}$ is well-defined. Since two lineages currently at separation $y\in\IR^2$ will coalesce if they are {\em both} involved in a replacement event, which happens at instantaneous rate
\begin{equation} \label{condition 3}
\int_{(|y|/2,\infty)}L_r(y,0)\left(\int_{[0,1]}u^2\nu_r(du)\right)\mu(dr),
\end{equation}
Condition~(\ref{condition for convergence}) is more than enough to bound the rate of coalescence of ancestral lineages. To guarantee that we can fit together the measures $\rho$ at different points in a consistent way, we also need to be able to control the spatial motion of ancestral lineages. Consider the (backwards in time) dynamics of a single ancestral lineage. It evolves in a series of jumps with intensity
\begin{equation}\label{jump intensity}
dt\otimes\int_{(|x|/2,\infty)}\int_{[0,1]}\frac{L_r(x,0)}{\pi r^2}\, u\,\nu_r(du)\mu(dr)dx
\end{equation}
on $\IR_+\times\IR^2$. If we want this to give a well-defined L\'evy process, then we require
\begin{equation}
\label{condition 2} \int_{\IR^2}(1\wedge |x|^2)\left(\int_{(|x|/2,\infty)}\int_{[0,1]} \frac{L_r(x,0)}{\pi r^2}\,u\,\nu_r(du)\mu(dr)\right)dx<\infty.
\end{equation}
But Condition~(\ref{condition for convergence}) certainly guarantees this. In fact it ensures that the rate of jumps of each ancestral lineage is {\em finite}. In other words, ancestral lineages follow compound Poisson processes.
\begin{remark}
At first sight it is disappointing that we have to take Condition~(\ref{condition for convergence}) and hence obtain a system of coalescing compound Poisson processes rather than more general symmetric L\'evy processes that (\ref{condition 1}) and (\ref{condition 2}) would allow. However, biologically there is not much loss. The `gap' between Condition~(\ref{condition for convergence}) and the weaker Condition~(\ref{condition 1}) is that the latter would allow one to include very large numbers of extremely small jumps (in which only a tiny proportion of the population is affected) as the radius of the area affected by a reproduction event tends to zero. But in our population model, for small $r$ we expect that a {\em large} proportion of the population in the neighbourhood be replaced.
\end{remark}
\begin{remark}
Notice that the locations of ancestral lineages are {\em not} independent of one another. Knowing that one lineage has jumped tells us that a reproduction event has taken place that could have affected other lineages ancestral to our sample. Wilkins \& Wakeley~(2002) consider a somewhat analogous model in which a linear population evolves in discrete generations (see Wilkins~2004 for a two-dimensional analogue). Each individual in the parental generation scatters an infinite pool of gametes in a Gaussian distribution about themselves, and the next generation is formed by sampling from the pool of gametes at each point. Individuals are assumed to have a finite linear width to avoid the pathologies that arise when common ancestry in a continuum model requires two ancestral lineages to have a physical separation of zero. They observe that ``conditional on not coalescing in the previous generation, two lineages are slightly more likely to be further apart than closer together''.  In their setting a change of coordinates settles the problem: the distance apart and the average position of two lineages do evolve independently. For us the dependencies between lineages are more complex because the presence of a jump contains the information that a reproduction event has taken place, whereas the conditioning obviously tells us nothing about the timing of events in the discrete generation model. \nocite{wilkins/wakeley:2002} \nocite{wilkins:2004}
\end{remark}

\section{The genealogy of points sampled uniformly from a large torus}
\label{section result}

We now turn our attention to populations evolving on a two-dimensional torus of sidelength $L$. Our goal is to describe the genealogy of a finite number of individuals sampled uniformly at random from the torus and subject to events of very different scales, as $L\rightarrow\infty$

To this end, we now consider a family of models indexed by $\N$. For each $L\in \N$, we consider a population evolving on the torus $\T(L)\subset \R^2$ of sidelength $L$. We identify $\T(L)$ with the subset $[-L/2,L/2]^2$ of $\R^2$ and use the Euclidean norm $|\cdot |$ induced on $\T(L)$ by this identification. Although $B_{\T(L)}(x,r)$ will be our notation for the ball in $\T(L)$ centred in $x$ and with radius $r$, we shall omit the subscript when there is no risk of confusion.

The population will be subject to two different classes of events that we call \emph{small} and \emph{large}. The region affected by each small event will be uniformly bounded (independently of the size of the torus). Large events will affect regions whose diameter is on the order of $\psi_L$ which will be taken to grow with $L$, but they will be less frequent. We shall assume that the rate at which a given ancestral lineage is affected by a large event is proportional to $1/\rho_L$ with $\rho_L$ also chosen to grow with $L$.

Now let us make the model more precise. Let $(\psi_L)_{L\geq 1}$ be an increasing sequence such that there exists $\alpha \in (0,1]$ satisfying
\begin{equation}\label{def alpha}
\lim_{L\rightarrow \infty}\frac{\log \psi_L}{\log L}=\alpha,
\end{equation}
and assume that $|\alpha\log L-\log \psi_L|=o((\log L)^{-1/2})$ as $L\rightarrow \infty$.
\begin{remark}
The latter assumption is not necessary since all our results would still hold with each occurrence of $(1-\alpha)\log L$ replaced by $\log (L\psi_L^{-1})$ (see the end of the proof of Proposition \ref{prop gathering}), but it is weak and considerably simplifies the presentation.
\end{remark}
Let $(\rho_L)_{L\geq 1}$ be an increasing sequence with values in $(0,+\infty]$, tending to infinity as $L\rightarrow \infty$. Finally, let $\mu^s(dr)$ and $\mu^B(dr)$ be two $\sigma$-finite Borel measures on $(0,\infty)$, independent of $L$, such that there exist some positive constants $R^s$ and $R^B$ satisfying
$$
\inf\big\{R:\mu^s\big((R,\infty)\big)=0\big\}=R^s <\infty \quad \mbox{ and }\quad \inf\big\{R:\mu^B\big((R,\infty)\big)=0\big\}=R^B <\infty.
$$
(For convenience, we ask that $R^B\leq 1/\sqrt{2}$ if $\alpha=1$.) To every $r\geq 0$, we associate two probability measures $\nu_r^s(du)$ and $\nu_r^B(du)$ on $[0,1]$, and we assume that for $\star \in \{B,s\}$ and for each $\e\in (0,R^\star)$,
\begin{equation}\label{coal at boundary}
\mu^\star\big(\big\{r\in [R^\star-\e,R^\star]:\nu_r^\star\neq \delta_0\big\}\big)>0.
\end{equation}
If Condition (\ref{coal at boundary}) does not hold, we decrease the corresponding radius $R^\star$ since otherwise the largest events never affect a lineage.

Let us suppose that for each $L\geq 1$, the reproduction events of the forwards in time model can be of two types :
\begin{itemize}
\item \textbf{Small events}, given by a Poisson point process $\Pi^s_L$ on $\R\times \T(L) \times (0,\infty)$ with intensity measure $dt\otimes dx \otimes \mu^s(dr)$. If $(t,x,r)$ is a point of $\Pi^s_L$, then the centre of the reproduction event is $x$, its radius is $r$ and the fraction of individuals replaced during the event is chosen according to $\nu_r^s$.
\item \textbf{Large events}, given by a Poisson point process $\Pi^B_L$ on $\R\times \T(L) \times (0,\infty)$, independent of $\Pi^s_L$ and with intensity measure $(\rho_L\psi_L^2)^{-1}dt\otimes dx \otimes \mu^B(dr)$. If $(t,x,r)$ is a point of $\Pi^B_L$, then the centre of the reproduction event is $x$, its radius is $\psi_Lr$ and the fraction of individuals replaced during the event is chosen according to $\nu_r^B$.
\end{itemize}
Notice that we allow $\rho_L$ to be infinite, in which case large events do not occur. Since $\Pi^s_L$ and $\Pi^B_L$ are independent, the reproduction events could be formulated in terms of a single Poisson point process to fit into the Definition \ref{def slfv} of the spatial $\Lambda$-Fleming-Viot process. However, our aim here is to disentangle the effects of events of different scales, hence our decomposition into two point processes.
\begin{remark}
Observe that, although the intensity of $\Pi_L^B$ is proportional to $(\rho_L\psi_L^2)^{-1}$, the rate at which a lineage is affected by (that is, jumps because of) a large event is of order $\mathcal{O}(\rho_L^{-1})$. Indeed, the volume of possible centres for such an event is proportional to $\psi_L^2$, so that the jump rate of a lineage due to the large events is given by
$$
\frac{1}{\rho_L\psi_L^2}\int_0^{R^B}\int_0^1 \pi (\psi_Lr)^2 u\ \nu_r^B(du)\mu^B(dr)= \frac{\pi}{\rho_L}\int_0^{R^B}\int_0^1 r^2 u\ \nu_r^B(du)\mu^B(dr).
$$
\end{remark}

In order for the genealogical processes, which we now denote by $\A^L$ to emphasize dependence on $L$, to be well-defined for every $L\in \N$, we assume that Condition~(\ref{condition for convergence}) is fulfilled. In this setting, the condition can be written
$$
\int_0^{R^s}\int_0^1r^2u\ \nu_r^s(du)\mu^s(dr) +\frac{1}{\rho_L}\int_0^{R^B}\int_0^1 r^2 u\ \nu_r^B(du)\mu^B(dr)<\infty.
$$

Let us introduce some more notation. We write
$$
\Gamma(L,1)\equiv \bigg\{x\in \T(L): |x|\geq \frac{L}{\log L}\bigg\},
$$
and for each integer $n\geq 2$,
\begin{eqnarray*}
\Gamma(L,n)& \equiv& \Big\{\{x_1,\ldots,x_n\} \in \T(L)^n: |x_i-x_j|\geq \frac{L}{\log L} \mathrm{\ \ for\ all\ }i\neq j\Big\},\\
\GA(L,n)&\equiv &\Big\{\big\{(a_1,x_{a_1}),\ldots,(a_k,x_{a_k})\big\}\in {\mathcal P}_n^{\ell}:\ \{x_{a_1},\ldots,x_{a_k}\}\in \Gamma(L,k)\Big\},
\end{eqnarray*}
where as before ${\mathcal P}_n^{\ell}$ denotes the labelled partitions of $\{1,\ldots ,n\}$. When we require an element $A$ of $\GA(L,n)$ to have exactly $n$ blocks, we shall write $A\in \GA(L,n)^*$.

In order to obtain a non-trivial limit, we rescale time for the process $\A^L$ by a factor that we denote $\vp_L$. Recall that if $A\in {\mathcal P}_n^{\ell}$ for some $n\in \N$, $\mathrm{bl}(A)$ stands for the unlabelled partition of $\{1,\ldots,n\}$ induced by $A$. For each $L\in\N$, let us define the (non-Markov) process $\A^{L,u}$ by
$$
\A^{L,u}(t)=\mathrm{bl}\big(\A^L(\vp_L t)\big),\qquad t\geq 0.
$$
Note that for each $L\in \N$, if we start $\A^L$ from $A_L$, a labelled partition of $\{1,\ldots,n\}$ with labels from $\IT(L)$, then $\A^{L,u}$ takes its values in the Skorohod space $D_{{\mathcal P}_n}[0,\infty)$ of all c\`adl\`ag paths with values in ${\mathcal P}_n$ (the set of partitions of $\{1,\ldots,n\}$), $\prob_{A_L}$-a.s.

Recall the definition of $\alpha$ given in (\ref{def alpha}). In the absence of large events, our model is similar in many respects to the two-dimensional stepping stone model and so it comes as no surprise that just as for the stepping stone model, the genealogy of a random sample from the torus should converge (on a suitable timescale) to a Kingman coalescent as the size of the torus tends to infinity (see in particular Cox \& Griffeath 1986,1990, Cox \& Durrett 2002 and Z\"ahle et al. 2005 for precise statements of this result in different contexts). Our first result says that if $\alpha<1$, then we still obtain a Kingman coalescent, but the {\em timescale} will be influenced by the large events: the latter reduce the effective population size.

Before stating the result formally, let us try to understand why we should expect something like this to be true. To understand the appropriate timescale we just need to consider two lineages. The time they need to coalesce will be decomposed into two phases. If $\rho_L$ is not too big, the first phase will be the time until they first come within distance $2R^B\psi_L$ and the second will be the additional time required for them to coalesce. During the first phase they evolve according to independent compound  Poisson processes. If $\rho_L$ is small enough, the coalescence event that will eventually occur during the second phase will, with probability close to one, be triggered by a large event. For larger values of $\rho_L$, large events will not be frequent enough to hit the two lineages when they are at a distance that would allow them to coalesce (i.e., less than $2R^B\psi_L$), and coalescence will instead be caused by a small-scale event. The first phase is then taken to be the time until the lineages first come within distance $2R^s$ of one another. The fact that with high probability they will not be hit by the same large-scale event means that once again they evolve (almost) independently of one another during this first phase. The second phase is now the time taken for them to coalesce due to a small event. The transition between these two regimes is when $\rho_L\propto\psi_L^2\log L$. Now suppose that we start from a sample in $\Gamma(L,n)$.  The first phase is then long enough that, when it ends, the spatial location of lineages is no longer correlated with their starting points. Finally, why do large-scale events not lead to multiple mergers? The key point is that, when a pair of lineages ancestral to our sample first comes within $2R^B\psi_L$ of one another, all {\em other} pairs are still well-separated. So if $\rho_L$ is not too big, this pair will coalesce before a third lineage can come close enough to be affected by a common event. If we take larger $\rho_L$, the reason is exactly the same but now lineages have to come within distance $2R^s$ and coalescence is driven by small events.

Here then is the formal result which makes explicit the convergence in distribution of our spatial genealogies to a nonspatial coalescent process. In the following, $\sigma_s^2$ (resp., $\sigma_B^2 \psi_L^2\rho_L^{-1}$) is the variance of the displacement of a lineage during one unit of time due to small (resp., large) events, see~(\ref{variances}) below.
\begin{thm}
\label{result alpha<1} Let $\mathcal{K}$ denote Kingman's coalescent, and recall that for each $n\in \N$, $\wp_n$ denotes the partition of $\{1,\ldots,n\}$ into singletons. In the notation of (\ref{def alpha}), suppose $\alpha<1$ (and (\ref{coal at boundary}) holds). Then, for each integer $n\geq 2$ and any sequence $(A_L)_{L\in \N}$ such that $A_L\in \GA(L,n)^*$ for every $L$,
$$
\mathcal{L}_{\prob_{A_L}}(\A^{L,u})\Rightarrow \mathcal{L}_{\prob_{\wp_n}}(\mathcal{K}) \qquad \mathrm{as\ }L\rightarrow \infty,
$$
where
$$
\vp_L = \left\{ \begin{array}{ll} \frac{(1-\alpha)\rho_LL^2\log L}{2\pi \sigma_B^2\psi_L^2} &\qquad\mathrm{if\ }\rho_L^{-1}\psi_L^2\rightarrow \infty, \vspace{4pt}\\
\frac{(1-\alpha)L^2\log L}{2\pi(\sigma_s^2+b\sigma_B^2)}& \qquad \mathrm{if\ }\rho_L^{-1}\psi_L^2\rightarrow b\in[0,\infty)\ \mathrm{and}\ \frac{\psi_L^2\log L}{\rho_L}\rightarrow \infty, \vspace{4pt}\\
\frac{L^2\log L}{2\pi \sigma_s^2}& \qquad \mathrm{if}\ (\rho_L^{-1}\psi_L^4)_{L\geq 1} \ \mathrm{is\ bounded\ or}\ \frac{L^2\log L}{\rho_L}\rightarrow 0.
\end{array}\right.
$$
Here $\mathcal{L}_{\mathrm{P}}(X)$ denotes the law under the probability measure $\mathrm{P}$ of the random variable $X$ and $\Rightarrow$ refers to weak convergence of probability measures.
\end{thm}
For $\alpha=1$, things are more complicated. When $\psi_L$ is commensurate with $L$, large scale events cover a non-negligible fraction of the torus. If they happen too quickly, then they will be able to capture multiple lineages while the locations of those lineages are still correlated with their starting points. For intermediate ranges of $\rho_L$, lineages will have homogenised their positions on $\IT(L)$ through small events, but not coalesced, before the first large event occurs and we can expect a $\Lambda$-coalescent limit. If the large events are too rare, then coalescence will be through small events and we shall recover the Kingman coalescent again.

To give a precise result we need to define the limiting objects that arise. In the case $\alpha =1$, for each $L\in \N$, we set
$$
\vp_L=\left\{\begin{array}{ll}\rho_L& \qquad \mathrm{if\ }\rho_L/(L^2\log L) \mathrm{\ has\ a\ finite\ limit}, \vspace{3pt}\\
\frac{L^2\log L}{2\pi \sigma_s^2}& \qquad \mathrm{if\ }\rho_L/(L^2\log L) \rightarrow +\infty, \end{array}\right.
$$
and define $\A^{L,u}$ as before. Since we shall need to keep track of the labels (spatial positions) of the ancestral lineages in some cases, it will also be convenient to introduce the following rescaling of $\A^L$, evolving on $\T(1)$ for all $L\in \N$:
$$
\bar{\A}^L(t)= \frac{1}{L}\ \A^L(\vp_L t),\qquad t\geq 0,
$$
where by this notation we mean that the labels are rescaled by a factor $L^{-1}$. Similarly, for $\mathbf{x}\in\IT(1)^n$ we write $L\mathbf{x}$ for $(Lx_1,\ldots ,Lx_n)\in\IT(L)^n$. Finally, let us introduce the processes which will appear as the limits of our rescaled genealogical processes.
\begin{defn}\label{def limit with space}
Let $b\in [0,\infty)$ and $c>0$. We call $\bar{\A}^{\infty,b,c}$ the Markov process with values in $\bigcup_{n\in \N}{\mathcal P}_n^{\ell}$ (with labels in $\T(1)$) such that \begin{enumerate}
\item The labels of the lineages perform independent Brownian motions on $\T(1)$ at speed $b\sigma_s^2$ (if $b=0$, the labels are constant), until the first large event occurs.
\item Large events are generated by a Poisson point process $\overline{\Pi}^B$ on $\R\times \T(1)\times (0,1/\sqrt{2}]$ with intensity measure $c^{-2}dt\otimes dx \otimes \mu^B(dr)$. At a point $(t,x,r)$ of $\ov{\Pi}^B$, a number $u\in [0,1]$ is sampled from the probability measure $\nu_r^B$, and each lineage whose label belongs to $B_{\T(1)}(x,cr)$ is affected (resp., is not affected) by the event with probability $u$ (resp., $1-u$), independently of each other. A label $z$ is chosen uniformly at random in $B_{\T(1)}(x,cr)$, and all the lineages affected merge into one block which adopts the label $z$. The other lineages (blocks and labels) remain unchanged.
\item The evolution of the labels starts again in the same manner.
\end{enumerate}
\end{defn}
\begin{remark}
Notice that this process looks like another spatial $\Lambda$-coalescent, except that now ancestral lineages perform independent spatial motions in between coalescence events. This process is dual (in the obvious way) to a spatial $\Lambda$-Fleming-Viot process in which, during their lifetimes, individuals move around in space according to independent Brownian motions.
\end{remark}
For each $r\in [0,1/\sqrt{2}]$, let $V_r$ denote the volume of the ball $B_{\T(1)}(0,r)$.
\begin{defn}\label{def lambda coal}Let $\beta\in [0,\infty)$ and $c>0$. We use $\Lambda^{(\beta,c)}$ to denote the $\Lambda$-coalescent, defined on $\bigcup_{n\in \N}{\mathcal P}_n$, for which if there are currently $m$ ancestral blocks, then each transition involving $k$ of them merging into one happens at rate
$$
\lambda^{(\beta,c)}_{m,k}= c^{-2}\int_0^{(\sqrt{2})^{-1}}\int_0^1(V_{cr}u)^k(1-V_{cr}u)^{m-k}\nu_r^B(du)\mu^B(dr)+ \beta\ \delta_{\{k=2\}}.
$$
\end{defn}
Recall the notation $\wp_n$ and $\wp_n({\mathbf x})$ introduced in Notation~\ref{notation for partitions}, and $\mathcal{L}_{\mathrm{P}}(X)$ and $\Rightarrow$ introduced in the statement of Theorem~\ref{result alpha<1}. We can now state the result for $\alpha=1$.

\begin{thm}
\label{result alpha=1} Suppose there exists $c>0$ such that for every $L\in \IN$, $\psi_L = cL$. Let $n\in\N$, $\mathbf{x}\in \T(1)^n$ such that $x_i\neq x_j$ whenever $i\neq j$, and let $(A_L)_{L\in \N}$ be such that for every $L$, $A_L\in \GA(L,n)^*$. Then, as $L\rightarrow \infty$,

\noindent $(a)$ If $\rho_LL^{-2}\rightarrow b\in [0,\infty)$,
$$
\mathcal{L}_{\prob_{\wp_n(L\mathbf{x})}}\big(\bar{\A}^L\big)\Rightarrow \mathcal{L}_{\prob_{\wp_n(\mathbf{x})}}\big(\bar{\A}^{\infty,b,c}\big),
$$

\noindent $(b)$ If $\rho_LL^{-2}\rightarrow \infty$, $\frac{2\pi \sigma_s^2 \rho_L}{L^2 \log L}\rightarrow \beta \in [0,\infty)$ and if the total rate of occurrence of large events is finite (i.e., $\mu^B$ has finite total mass),
$$
\mathcal{L}_{\prob_{A_L}}\big(\A^{L,u}\big)\Rightarrow \mathcal{L}_{\prob_{\wp_n}}\big(\Lambda^{(\beta,c)}\big).
$$

\noindent $(c)$ If $\frac{\rho_L}{L^2 \log L}\rightarrow \infty$,
$$
\mathcal{L}_{\prob_{A_L}}\big(\A^{L,u}\big)\Rightarrow \mathcal{L}_{\prob_{\wp_n}}\big(\mathcal{K}\big).
$$
\end{thm}
Notice that the case $(a)$ differs from all other cases in that the influence of space does not disappear as $L\rightarrow \infty$ and the evolution of the limiting genealogy still depends on the precise locations of the lineages.

The intuition behind Theorem \ref{result alpha=1} is as follows. If $\psi_L \propto L$ large events cover a non-negligible fraction of the torus, and so only a few large events are sufficient to gather two lineages at a distance at which they can coalesce. However, a local central limit theorem will give us that on a timescale of order at most $\mathcal{O}(L^2)$, a lineage subject to only small events behaves approximately like Brownian motion, whereas after a time $t_L\gg L^2$, its distribution is nearly uniform on $\T(L)$ (for $L$ large enough, see Lemma \ref{lemma local TCL}). Since the mean time before a large event affects a lineage is of order ${\mathcal O}(\rho_L)$, the limiting genealogical process (when we include both large and small reproduction events) will depend on how $\rho_L$ scales with $L^2$. If $\rho_L$ is of order at most $\mathcal{O}(L^2)$, then space matters and the process $\A^L$ rescaled to evolve on $\T(1)$ on the timescale $\rho_L$ converges to a system of coalescing Brownian motions, whereas if $\rho_L\gg L^2$, the homogenisation of the labels/locations of the lineages before the occurrence of the first large event which affects them leads to a limiting unlabelled genealogical process given by an exchangeable coalescent with multiple mergers.
\begin{remark}\label{rk finite rate}
It is somehow disappointing that we must impose a finite rate of large events to obtain the convergence of Theorem \ref{result alpha=1}(b). Indeed, it seems that case (a) should give us the right picture: in the limit, in between large events lineages perform Brownian motions on the torus of sidelength 1 due to small events, except that now the time required for at least one lineage to be affected by a large event is so long that lineages exhaust space and their locations become uniformly distributed over the torus before they are taken by a coalescence event. However, when $\mu^B$ has infinite mass, lineages are infinitely often in the (geographical) range of a large reproduction event over any interval of time, and we need good control of their complete paths to actually be able to say something about the epoch and outcome of the first potential coalescence event. Now, observe that Equation (\ref{eq homogen}) can only be generalized to the finite-dimensional distributions of these paths, and does not guarantee that a large event cannot capture some of the lineages at a time when they are not uniformly distributed over $\T(1)$.
\end{remark}

Theorem \ref{result alpha=1} deals with the case where $\psi_L$ is proportional to $L$. Let us now comment on the remaining cases, in which $\alpha=1$ but $\psi_L\ll L$. First, it is easy to see that the convergence in $(c)$ still holds, since it is based on the fact that large events are so rare that none of them occurs before small events reduce the genealogical process to a single lineage.

Second, since the total rate of large events on the timescale $\rho_L$ is $\mu^B(\IR_+)L^2/\psi_L^2$, it cannot be bounded unless $\mu^B\equiv 0$ (a situation we excluded in (\ref{coal at boundary})). On the other hand, for the reason expounded in Remark \ref{rk finite rate} we are unable to derive a limiting behaviour for the genealogy when large events can accumulate, and so the result of Theorem \ref{result alpha=1}$(b)$ has no counterpart when $\psi_L\ll L$.

Third, as explained above, when $\rho_L\leq b L^2$ any limiting process will necessarily have a spatial component. Now, because we start with lineages at distance $\mathcal{O}(L)$ of each other, we need to rescale space by $L$ in order to obtain a non trivial initial condition. The last parameter we need is the timescale $\varpi_L$ on which to consider the genealogical process. But a separation of timescales will not occur here, and so the computations done in Section \ref{levy processes} will show that the suitable choice of $\varpi_L$ depends on the precise behaviour of $\rho_L/L^2$ and $\rho_L/\psi_L^2$. Several limiting processes are thus possible, and since all the arguments needed to derive these limits are scattered in Sections \ref{levy processes} and \ref{alpha=1}, we chose not to detail them here.

\section{Existence and uniqueness of the forwards-in-time process}
\label{section existence}

Our spatial $\Lambda$-Fleming-Viot process associates a probability measure on type space to each point in $\IR^2$. In other words, it takes its values among functions from $\IR^2$ to ${\mathcal M}_1([0,1])$. \nocite{evans:1997} Evans~(1997) uses duality with a system of coalescing Borel right processes on a Lusin space $E$ to construct a family of Markov processes with values in the set of functions from $E$ to $\mathcal{M}_1(\{0,1\}^{\N})$ (or equivalently, to $\mathcal{M}_1([0,1])$). He also obtains uniqueness in distribution of the process. In his setting, coalescing particles evolve independently until they meet, at which point they instantly coalesce. In our case, the particles in the candidate dual do not move independently and nor do two particles hit by the same reproduction event necessarily coalesce, but nonetheless the key ideas from his construction remain valid. Note that, although we present the result in two dimensions, the proof carries over to other dimensions.

First we give a formal description of the coalescing dual and then we use the Evans' construction to give existence and uniqueness in law of a process $\rho$ which assigns a probability measure on $[0,1]$ to each point in $\IR^2$. We then identify $\rho$ as the spatial $\Lambda$-Fleming-Viot process in which we are interested.

\subsection{State-space of the process and construction via duality}
We shall only present the main steps of the construction, and refer to Evans~(1997) for more details.

Let us define $\tilde{\Xi}$ as the space of all Lebesgue-measurable maps $\rho:\R^2 \rightarrow \mathcal{M}_1([0,1])$. Two elements $\rho_1$ and $\rho_2$ of $\tilde{\Xi}$ are said to be equivalent if $\mathrm{Leb}(\{x\in \R^2:\ \rho_1(x)\neq \rho_2(x)\})=0$. Let $\Xi$ be the quotient space of $\tilde{\Xi}$ by this equivalence relation. If $E$ is a compact space, let us write $C(E)$ for the Banach space of all continuous functions on $E$, equipped with the supremum norm $\| \cdot \|_{\infty}$. For each $n\in \N$, let $L^1(C([0,1]^n))$ be the Banach space of all Lebesgue-measurable maps $\Phi:(\R^2)^n\rightarrow C([0,1]^n)$ such that $\int_{(\R^2)^n} \|\Phi(x)\|_{\infty}\ dx <\infty$. A remark in Section 3 of Evans~(1997) tells us that the separability of $L^1(C([0,1]))$ and a functional duality argument guarantee that $\Xi$, equipped with the relative weak* topology, is a (compact) metrisable space. Finally, if $\lambda$ is a measure on a space $E'$, let us write $L^1(\lambda)$ for the set of all measurable functions $f:E'\rightarrow \R$ such that $\int_{E'}|f(e)|\lambda(de)<\infty$.

Let $n\in \N$. Given $\Phi\in L^1(C([0,1]^n))$, let us define a function $I_n(\cdot\ ;\Phi)\in C(\Xi)$ by
$$
I_n(\rho;\Phi)\equiv \int_{(\R^2)^n} \Big\langle \bigotimes_{1\leq i\leq n}\rho(x_i),\Phi(x_1,\ldots,x_n)\Big\rangle\ dx_1\ldots dx_n,
$$
where as before the notation $\langle \nu,f\rangle$ stands for the integral of the function $f$ against the measure $\nu$. We have the following lemma, whose proof is essentially that of Lemma~3.1 in Evans~(1997).
\begin{lemma}\label{lemm set of functions}The linear subspace spanned by the constant functions and functions of the form $I_n(\cdot\ ;\Phi)$, with $\Phi= \psi\otimes \big(\prod_{i=1}^n\chi_i\big)$, $\psi\in L^1(dx^{\otimes n})\cap C((\R^2)^n)$ and $\chi_i\in C([0,1])$ for all $1\leq i\leq n$ is dense in $C(\Xi)$.
\end{lemma}

We need a last definition before stating the existence and uniqueness result. Let $n\in \N$. For any $\rho\in \Xi$, $\pi\in {\mathcal P}_n^{\ell}$ such that $\mathrm{bl}(\pi)=\{a_1,\ldots,a_k\}$, and any bounded measurable function $F:[0,1]^n\rightarrow \R$, we set
$$
\Upsilon_n(\rho;\pi;F)\equiv \int_{[0,1]^k}F(v_{a^{-1}(1)},\ldots,v_{a^{-1}(n)}) \rho(x_{a_1})(dv_{a_1})\ldots \rho(x_{a_k})(dv_{a_k}),
$$
where $a^{-1}(i)$ is the (unique) block $a_j$ which contains $i$ and $v_{a_j}$ is the variable used for the measure $\rho(x_{a_j})$. In words, we assign the same variable to all coordinates which belong to the same block in the partition $\pi$.  (Recall that $x_a$ is our notation for the label of block $a$.)  Recall also the notation $\wp_n(\mathbf{x})$ and $\mathcal{A}$ introduced in Notation \ref{notation for partitions} and the following paragraph.
\begin{thm}\label{theo existence}There exists a unique, Feller, Markov semigroup $\{Q_t,t\geq 0\}$ on $\Xi$ such that for all $n\in \N$ and $\Phi \in L^1(C([0,1]^n))$, we have
\begin{equation}\label{def semigroup}
\int Q_t(\rho,d\rho')I_n(\rho';\Phi)=\int_{(\R^2)^n} \E_{\wp_n(\mathbf{x})}\big[\Upsilon_n\big(\rho;\mathcal{A}(t);\Phi(x_1,\ldots,x_n)\big)\big] dx_1\ldots dx_n.
\end{equation}
Consequently, there exists a Hunt process $\{\rho(t),t\geq 0\}$ with state-space $\Xi$ and transition semigroup $\{Q_t,t\geq 0\}$.
\end{thm}

Before proving Theorem \ref{theo existence}, let us make two comments on this result. First, since the $\Xi$-valued process we obtain is a Hunt process it is c\`adl\`ag and quasi-left continuous, that is, it is almost surely left-continuous at any previsible stopping time (see e.g. Rogers \& Williams 1987 \nocite{rogers/williams:1987} for a definition of quasi-left continuous filtrations). However, more precise statements on its space-time regularity seem to be a delicate question, which will require a thorough investigation.

Second, as in Kimura's stepping stone model introduced in (\ref{stepstone model}), the duality relation (\ref{def semigroup}) can be interpreted in terms of genealogies of a sample of individuals. Indeed, recall the stepping stone model is dual to the system $(\{n_i(t);\ i\in I\})_{t\geq 0}$ of particles migrating from deme $i$ to deme $j$ at rate $m_{ji}$ and coalescing in pairs at rate $1/N_e$ when in the same deme: for any $t\geq 0$, we have
$$
\E\bigg[\prod_{i\in I}p_i(t)^{n_i(0)}\bigg]=\E\bigg[\prod_{i\in I}p_i(0)^{n_i(t)}\bigg].
$$
These equations show that a function (here the $n_i(0)$-th moments) of the frequencies at different sites of $\Z^2$ and at (forward) time $t$ can be expressed in terms of the genealogy of a sample made of $n_i(0)$ individuals in deme $i$ for every $i\in I$, and run for a (backward) time $t$: all lineages having coalesced by time $t$ necessarily carry the same type, whose law is given by the type distribution at the site where their ancestor lies at backward time $t$ (or forward time $0$). Equation~(\ref{def semigroup}) can be interpreted in exactly the same manner, but holds for a much wider collection of functions of $\rho$ and $\mathcal{A}$.

\medskip
\noindent\emph{Proof of Theorem \ref{theo existence}: } The observation that the construction of Evans~(1997) can also be justified in our setting follows from Remark~(a) at the end of his Section~4.

Existence and uniqueness of $\mathcal{A}$ are easy from Assumptions (\ref{condition 3}) and (\ref{condition 2}). Next, we must verify consistency of $\mathcal{A}$ in the sense of his Lemma~2.1. In fact, this is the `sampling consistency' described in the introduction and was a primary consideration in writing down our model. It follows since the movement of the labels of a collection of blocks does not depend on the blocks themselves and from the fact that a coalescence event of the form $\{(\{1\},x_1),(\{2\},x_2)\}\rightarrow \{(\{1,2\},x)\}$ for a pair of particles corresponds to a jump $\{(\{1\},x_1)\}\rightarrow \{(\{1\},x)\}$ onto the same site $x\in \R^2$ if we restrict our attention to the first particle.

The next property needed in the construction is that provided it is true at $t=0$, for every $t>0$ the distribution of the labels in $\mathcal{A}(t)$ has a Radon-Nikodym derivative with respect to Lebesgue measure, and furthermore an analogue of Evans' Equation (4.2) holds. In the setting of Evans~(1997), the first requirement stems from
the independence of the spatial motions followed by different labels and the corresponding result for a single label. Here, since the motion of all lineages is driven by the same Poisson process of events, their movements are correlated. However, the desired property is still satisfied. To see this, note that each jump experienced by a lineage in
the interval $[-t,0]$ takes it to a position that is uniformly distributed over the open ball affected by the corresponding reproduction event. Thus, if $\A(t)$ has $k$ blocks and $D\subset (\R^2)^k$ has zero Lebesgue measure, the probability that the labels of the blocks of $\A(t)$ belong to $D$ is equal to $0$. Equation (4.2) of Evans~(1997) then still holds, without Evans' additional assumption of the existence of a dual process for the motion of one lineage (which anyway is satisfied since our lineages perform symmetric L\'evy processes).

The last step is to check the strong continuity of the semigroup $\{Q_t,t\geq 0\}$, but this readily follows from the relation (\ref{def semigroup}) and the Feller property of $\A$ (which is itself evident since jumps do not accumulate in our dual process).

The desired conclusion now follows from Theorem~4.1 in Evans~(1997).$\hfill\square$

\subsection{Identification of the process}

We can use (\ref{def semigroup}) to derive an expression for the infinitesimal generator of $\{\rho(t),t\geq 0\}$ acting on the functions $I_n(\cdot\ ;\Phi)$ considered in Lemma~\ref{lemm set of functions}. This lemma and the uniqueness result stated in Theorem~\ref{theo existence} guarantee that it will be sufficient to characterize the process $\rho$ and to show that it corresponds to the evolution we described in Section~\ref{model} in terms of a Poisson point process of reproduction events.

Let $n\in \N$ and $\Phi\in C(\Xi)$ be such that $\Phi= \psi\otimes \big(\prod_{i=1}^n\chi_i\big)$, where $\psi\in L^1(dx^{\otimes n})\cap C((\R^2)^n)$ and $\chi_i\in C([0,1])$ for all $1\leq i\leq n$. Writing $G$ for the generator of the process $\rho$ and $\G_n$ for the generator of the coalescing L\'evy processes $\mathcal{A}$ acting on functions of ${\mathcal P}_n^{\ell}$, we obtain from (\ref{def semigroup}) that \setlength\arraycolsep{1pt}
\begin{eqnarray}
 GI_n(&\rho&; \Phi)= \lim_{t\rightarrow 0}\frac{\E_{\rho}[I_n(\rho(t),\Phi)]-I_n(\rho,\Phi)}{t} \nonumber\\
&=& \lim_{t\rightarrow 0}\frac{1}{t}\int_{(\R^2)^n} \psi(x_1,\ldots,x_n) \bigg\{ \E_{\wp_n(\mathbf{x})}\Big[\Upsilon_n\Big(\rho;\mathcal{A}(t);\prod_{i=1}^n \chi_i\Big)\Big] - \prod_{i=1}^n \langle \rho(x_i),\chi_i\rangle \bigg\}\ dx^{\otimes n} \nonumber\\
& =& \int_{(\R^2)^n} \psi(x_1,\ldots,x_n)\ \G_n\Big[\Upsilon_n\Big(\rho;\ \cdot\ ;\prod_{i=1}^n \chi_i\Big)\Big](\wp_n(\mathbf{x}))\ dx^{\otimes n}\label{def generator on Xi}.
\end{eqnarray}
Note that the quantity on the right-hand side of (\ref{def generator on Xi}) is well-defined (and the interchange of limit and integral is valid) since $\psi$ belongs to $L^1(dx^{\otimes n})$ and the rate at which at least one of $k\leq n$ blocks is affected by a reproduction event is bounded by $n$ times the integral in (\ref{condition for convergence}), so that $\A$ is a jump-hold process and its generator satisfies
$$
\Big\|\G_n\Big[\Upsilon_n\Big(\rho;\ \cdot\ ;\prod_{i=1}^n \chi_i\Big)\Big]\Big\|_{\infty}\leq 2Cn\ \Big\|\Upsilon_n\Big(\rho;\ \cdot\ ;\prod_{i=1}^n \chi_i\Big)\Big\|_{\infty}\leq 2Cn \prod_{i=1}^n\|\chi_i\|_{\infty} <\infty
$$
for a given constant $C<\infty$.

Using the description of the evolution of $\A$ in terms of events in $\Pi$, the right-hand side of (\ref{def generator on Xi}) is equal to \setlength\arraycolsep{1pt}
\begin{eqnarray}
\int_{(\R^2)^n}& &dx^{\otimes n} \psi(x_1,\ldots,x_n)\int_{\R^2}dy\int_0^{\infty}\mu(dr)\int_0^1\nu_r(du) \int_{B(y,r)}\frac{dz}{\pi r^2}\nonumber\\
&  \times& \sum_{I\subset\{1,\ldots,n\}}\bigg[\prod_{i\in I}\mathbf{1}_{B(y,r)}(x_i)\prod_{i'\notin I}\mathbf{1}_{B(y,r)^c}(x_{i'})\bigg] \nonumber \\
&  \times& \sum_{J\subset I}u^{|J|}(1-u)^{|I|-|J|}\bigg[\prod_{i\notin J}\big\langle \rho(x_i),\chi_i\big\rangle \bigg]\bigg[\Big\langle \rho(z),\prod_{j\in J}\chi_j\Big\rangle- \prod_{j\in J}\big\langle \rho(x_j),\chi_j\big\rangle \bigg],\phantom{AAA} \label{alt expression}
\end{eqnarray}
where $|\cdot|$ stands for cardinality. Indeed, given $x_1,\ldots,x_n$ in (\ref{alt expression}), only one term in the sum over $I \subset \{1,\ldots,n\}$ is non-zero. For this particular term, each of the $|I|$ blocks whose labels lie in $B(y,r)$ belong to the set $J$ of the blocks affected by the event with probability $u$ (independently of one another), and the affected blocks adopt the label $z$. After some algebra and several uses of Fubini's theorem, we obtain that (\ref{alt expression}) is equal to \setlength\arraycolsep{1pt}
\begin{eqnarray}
\int_{\R^2}dy \int_0^{\infty} &\mu&(dr)\int_0^1\nu_r(du)\int_{B(y,r)}\frac{dz}{\pi r^2}\int_0^1\rho_z(dk)\int dx_1\ldots dx_n\ \psi(x_1,\ldots,x_n)\nonumber\\
 &\times& \sum_{I\subset \{1,\ldots,n\}}\prod_{j\notin I}\big\{\mathbf{1}_{B(y,r)^c}(x_j)\langle \rho_{x_j},\chi_j\rangle\big\}\prod_{i\in I}\mathbf{1}_{B(y,r)}(x_i)\nonumber
\\
&&\qquad \qquad\times \bigg(\prod_{i\in I}\big\langle(1-u)\rho_{x_i}+u\delta_k,\chi_i\big\rangle- \prod_{i\in I}\big\langle \rho_{x_i},\chi_i\big\rangle\bigg),\label{mp for lfv}
\end{eqnarray}
which is precisely the generator of the forwards in time process of Section \ref{model}. Using Theorem~\ref{theo existence}, we arrive at the following result. \begin{propn}\label{prop mp for lfv} The martingale problem associated to the operator $G$ defined by (\ref{mp for lfv}) on functions of the form given in Lemma~\ref{lemm set of functions} is well-posed. Furthermore, the spatial $\Lambda$-Fleming-Viot process $\rho$ of Theorem~\ref{theo existence} is the solution to it.
\end{propn}

\section{Some estimates for symmetric L\'evy processes}
\label{levy processes}

In this section, we gather some results on symmetric L\'evy processes that we shall need to call upon in our proofs of Theorem~\ref{result alpha<1} and Theorem~\ref{result alpha=1}. For the sake of clarity, the proofs of the three lemmas are given in Appendix \ref{appendix 1}.

First, we introduce some notation that we shall use repeatedly.
\begin{notn}
\label{entrance time}
\begin{enumerate}
\item{In the following, we shall suppose that all the random objects considered are constructed on the same probability space $(\Omega, \mathcal{F}, \prob)$, and if $X$ is a process defined on $\Omega$ with state-space $E$ and $x\in E$, we shall write $\prob_x$ for the probability measure on $\Omega$ under which $X(0)=x$ a.s. }
\item{For a stochastic process $\{X_t\}_{t\geq 0}$ evolving in $\IT(L)$, we shall write $T(R,X)$ for the {\em first entrance time} of $X$ into $B_{\T(L)}(0,R)$. When there is no ambiguity, we write simply $T(R)$.}
\end{enumerate}
\end{notn}

Let $(\ell^L)_{L\geq 1}$ be a sequence of L\'evy processes such that for each $L\in \N$, $\ell^L$ evolves on the torus $\T(L)$ and $\ell^L(1)-\ell^L(0)$ has a covariance matrix of the form $\sigma_L^2\mathrm{Id}$. Assume that the following conditions hold.
\begin{assumption}
\label{assumptions for levy processes}
\begin{description}
\item [(i)] There exists $\sigma^2>0$ such that $\sigma_L^2\rightarrow \sigma^2$ as $L\rightarrow \infty$.
\item [(ii)] $\E_0\big[|\ell^L(1)|^4\big]$ is bounded uniformly in $L$.
\end{description}
\end{assumption}
Our first lemma describes the time $\ell^L$ needs to reach a ball of radius $d_L\ll L$ around $0$, when it starts at distance $\mathcal{O}(L)$ of the origin (recall the definition of $\Gamma(L,1)$ given in Section \ref{section result}).
\begin{lemma}\label{lemma entrance_levy}Let $(d_L)_{L\geq 1}$ be such that $\liminf_{L\rightarrow \infty}d_L >0$ and $\frac{\log^+(d_L)}{\log L}\rightarrow \gamma \in [0,1)$ as $L\rightarrow \infty$. Then,
\begin{equation}
\label{eq entrance_levy}\lim_{L\rightarrow \infty}\sup_{t\geq 0}\sup_{x_L\in \Gamma(L,1)} \left|\prob_{x_L}\left[T(d_L,\ell^L)> \frac{(1-\gamma)L^2\log L}{\pi \sigma^2}\ t\right]- e^{-t}\right|=0.
\end{equation}
\end{lemma}
The proof of Lemma \ref{lemma entrance_levy} follows that of Theorem~2 in Cox \& Durrett~(2002). \nocite{cox/durrett:2002} In particular, we shall use the following local central limit theorem (which is the counterpart in our setting of Lemma 3.1 in Cox \& Durrett~2002). Let $\lfloor z \rfloor$ denote the integer part of $z\in \R$, and write $p^L(x,t)$ for $\prob_x[\ell^L(t)\in B(0,d_L)]$.
\begin{lemma}\label{lemma local TCL}

\noindent$(a)$ Let $\e_L=(\log L)^{-1/2}$. There exists a constant $C_1<\infty$ such that for every $L\geq 2$,
\begin{equation}
\sup_{t\geq \integ{\e_LL^2}}\ \sup_{x\in \T(L)}\ \frac{\integ{\e_LL^2}}{d_L^2}\ p^L(x,t)\leq C_1.
\end{equation}
$(b)$ If $v_L\rightarrow \infty$ as $L\rightarrow \infty$, then
\begin{equation}
\lim_{L\rightarrow \infty}\ \sup_{t\geq \integ{v_LL^2}}\ \sup_{x\in \T(L)}\ \frac{L^2}{d_L^2}\left|\ p^L(x,t)-\frac{\pi d_L^2}{L^2}\right|=0.
\end{equation}
$(c)$ If $u_L\rightarrow \infty$ as $L\rightarrow \infty$ and $I(d_L,x)\equiv 1+(|x|^2\vee d_L^2)$, then
\begin{equation}
\lim_{L\rightarrow \infty}\ \sup_{x\in \T(L)}\ \sup_{u_LI(d_L,x)\leq t\leq \e_LL^2}\ \left|\frac{2\sigma_L^2t}{d_L^2}\ p^L(x,t)-1 \right|=0.
\end{equation}
$(d)$ There exists a constant $C_2<\infty$ such that for every $L\geq 1$,
\begin{equation}
\sup_{t\geq 0}\sup_{x\in \T(L)}\left(1+\frac{|x|^2}{d_L^2}\right)p^L(x,t)\leq C_2.
\end{equation}
\end{lemma}
In essence, Lemma \ref{lemma local TCL} says that on the timescale $d_L^2\ll t\ll L^2$, the L\'evy process $\ell^L$ behaves like two-dimensional Brownian motion, whereas at any given time $t\gg L^2$, its location is roughly uniformly distributed over $\IT(L)$.

Another consequence of Lemma \ref{lemma local TCL} is the following result, which bounds the probability that $\ell^L$ hits a ball of bounded radius during a `short' interval of time in the regime $t\gg L^2$.
\begin{lemma}\label{lemm no entrance}Fix $R>0$. Let $(U_L)_{L\geq 1}$ and $(u_L)_{L\geq 1}$ be two sequences increasing to infinity such that $U_LL^{-2}\rightarrow \infty$ as $L\rightarrow \infty$ and $2u_L\leq L^2(\log L)^{-1/2}$ for every $L\geq 1$. Then, there exist $C>0$ and $L_0\in \N$ such that for every sequence $(U'_L)_{L\geq 1}$ satisfying $U_L'\geq U_L$ for each $L$, every $L\geq L_0$ and all $x\in \T(L)$,
$$
\prob_x\Big[T(R,\ell^L) \in [U_L'-u_L, U_L'\big]\Big]\leq \frac{C u_L}{L^2}.
$$
\end{lemma}

\section{Proof of Theorem~\ref{result alpha<1}}
\label{alpha<1}

Armed with the estimates of Section \ref{levy processes}, we can now turn to the proofs of our main results.
\begin{notn}
\label{motion of ancestral lineages} For each $L\geq 1$, let $\{\xi^L(t),t\geq 0\}$ be the L\'evy process on $\T(L)$ whose distribution is the same as that of the motion of a single lineage subject to the large and small reproduction events generated by $\Pi_L^s$ and $\Pi_L^B$.
\end{notn}
In the rest of this section, we assume that the assumptions of Theorem \ref{result alpha<1} are satisfied.

\subsection{Coalescence time for two lineages}\label{section coal}
We begin by studying the genealogical process of a pair of lineages starting at distance $\mathcal{O}(L)$ from each other. Since the motions $\xi_1^L$ and $\xi_2^L$ of the lineages are distributed like two independent copies of the process $\xi^L$ until the random time $T_L$ at which they come at distance less than $2R^B\psi_L$, the difference $$
X^L(t)\equiv \xi_1^L(t)-\xi_2^L(t),\qquad 0\leq t\leq T_L
$$
has the same distribution as $\big\{\xi^L(2t),\ 0\leq t\leq \frac{1}{2}\ T(2R^B\psi_L,\xi^L)\big\}$. We shall use Lemma~\ref{lemma entrance_levy} to derive the limiting distribution of $T_L$, but first we need to introduce the relevant variances. Consider a single lineage. Because it jumps at a finite rate owing to small and large events, the following two quantities are well-defined and finite~:
\begin{equation}
\label{variances} \sigma_s^2\equiv \int y^2\ \chi^s(dy,dz) \qquad \mathrm{and}\qquad \sigma_B^2 \equiv \int y^2\ \chi^B(dy,dz),
\end{equation}
where $\chi^s$ stands for the intensity measure of the small jumps experienced by the lineage and $\chi^B$ for that of the large jumps renormalised by $\psi_L^{-1}$ (the form of these two measures is given in (\ref{jump intensity})). We now have all the ingredients we need to describe the asymptotic `gathering time' of two lineages.
\begin{propn}\label{prop gathering} $(a)$ If $\rho_L^{-1}\psi_L^2\rightarrow \infty$ as $L\rightarrow \infty$, then
$$
\lim_{L\rightarrow \infty}\ \sup_{t\geq 0}\ \sup_{A_L\in \GA(L,2)^*}\ \bigg|\ \prob_{A_L}\left[T_L> \frac{(1-\alpha)\rho_LL^2\log L}{2\pi \sigma_B^2\psi_L^2}\ t\right] -e^{-t}\bigg|=0.
$$

\noindent$(b)$ If $\rho_L^{-1}\psi_L^2\rightarrow b \in [0,\infty)$ as $L\rightarrow \infty$, then
$$
\lim_{L\rightarrow \infty}\ \sup_{t\geq 0}\sup_{A_L\in \GA(L,2)^*}\ \bigg|\ \prob_{A_L}\left[T_L> \frac{(1-\alpha)L^2\log L}{2\pi (\sigma_s^2 + b \sigma_B^2)}\ t\right] - e^{-t}\bigg|=0.
$$
\end{propn}

\noindent \emph{Proof of Proposition \ref{prop gathering}: } Let us first recall two results on Poisson point processes, which are consequences of the exponential formula given, for instance, in Section 0.5 of \nocite{bertoin:1996} Bertoin~(1996). Following Bertoin's notation, let $\{e(t),t\geq 0\}$ be a Poisson point process on $\R\times \R_+$ with intensity measure $\kappa(dy)\otimes dt$, where the Borel measure $\kappa$ satisfies
\begin{equation}
\int_{\R}|1-e^y|\kappa(dy)<\infty \qquad \mathrm{and}\qquad \int_{\R}y^m \kappa(dy)=0, \quad m\in \{1,3\}.
\end{equation}
Under these conditions, we have for each fixed $t>0$
\begin{eqnarray}
\E\bigg[\Big(\sum_{s\leq t}e(s)\Big)^2\bigg]&=& t \int_{\R} y^2 \kappa(dy), \label{PPP moment 2} \\
\E\bigg[\Big(\sum_{s\leq t}e(s)\Big)^4\bigg]&=& 3t^2\bigg(\int_{\R}y^2 \kappa(dy)\bigg)^2+t\int_{\R}y^4 \kappa(dy). \label{PPP moment 4}
\end{eqnarray}
These properties will be useful in computing the variances and fourth moments of the random variables considered below.

Let us start with the proof of $(a)$. Consider the process $\ell^L$ defined by: for every $t\geq 0$,
$$
\ell^L(t)=\frac{1}{\psi_L}\ \xi^L\big(2\rho_Lt\big).
$$
This process evolves on the torus of sidelength $\psi_L^{-1}L$, and makes jumps of size $\mathcal{O}(\psi_L^{-1})$ at a rate of order $\mathcal{O}(\rho_L)$, as well as jumps of size $\mathcal{O}(1)$ at a rate of order $\mathcal{O}(1)$.

Let us check that $\ell^L$ satisfies the assumptions of Lemma \ref{lemma entrance_levy}. To this end, we view $\ell^L(1)$ starting at $0$ as the sum of its jumps and adapt the problem to use the results on Poisson point processes given above. First, let us define $\hat{\ell}^L$ as the L\'evy process on $\R^2$ evolving like $\ell^L$ (but without periodic conditions). For $i\in \{0,1\}$ and each $L\geq 1,\ t\geq 0$, let $\hat{\ell}^{L,i}(t)$ denote the $i$-th coordinate of $\hat{\ell}^L(t)$. Note that the distance reached by $\ell^L$ up to a given time $t$ is less than or equal to the distance at which $\hat{\ell}^L$ traveled up to $t$, and so we can write
\begin{eqnarray*}
\E_0\big[|\ell^L(1)|^4\big]\leq \E_0\big[|\hat{\ell}^L(1)|^4\big]&=& \E_0\Big[\Big\{ \hat{\ell}^{L,1}(1)^2 + \hat{\ell}^{L,2}(1)^2\Big\}^2\Big] \\
&\leq& 2 \Big\{\E_0\big[\hat{\ell}^{L,1}(1)^4\big]+\E_0\big[\hat{\ell}^{L,2}(1)^4\big]\Big\}.
\end{eqnarray*}
By symmetry, we need only bound $\E_0\big[\hat{\ell}^{L,1}(1)^4\big]$. Let us denote by $a_1,a_2,\ldots\in [-2R^s/\psi_L,$ $2R^s/\psi_L]^2$ (resp., $b_1,b_2,\ldots\in [-2R^B,2R^B]^2$) the sequence of the jumps of $\hat{\ell}^{L,1}$ before time $1$ due to small (resp., large) events. Using the convexity of $y\mapsto y^4$, we have \begin{equation}\label{eq different jumps}
\E_0\big[\hat{\ell}^{L,1}(1)^4\big]=\E_0\bigg[\Big(\sum_i a_i + \sum_j b_j\Big)^4\bigg]\leq 8\ \E_0\bigg[\Big(\sum_i a_i\Big)^4 + \Big(\sum_j b_j\Big)^4\bigg].
\end{equation}
Applying (\ref{PPP moment 4}) to each term on the right-hand side of (\ref{eq different jumps}) yields
\begin{equation}\label{bound coord}
\E_0\big[(\hat{\ell}^{L,1}(1))^4\big]\leq 96\frac{\rho_L^2}{\psi_L^4}\ \sigma_s^4 + 16\frac{\rho_L}{\psi_L^4}\ \int y^4 \chi^s(dy,dz)+96 \sigma_B^4 +16\int y^4 \chi^B(dy,dz),
\end{equation}
which is bounded uniformly in $L$ since $\rho_L\psi_L^{-2}$ vanishes as $L$ grows to infinity, and each integral is finite. Coming back to the original problem, we obtain that Assumption~\ref{assumptions for levy processes} (ii) holds for the sequence of processes $(\ell^L)_{L\geq 1}$.

Concerning Assumption~\ref{assumptions for levy processes} (i), observe that $\sigma_L^2$ is simply the variance of $\ell^{L,1}(1)$. To obtain the asymptotic behaviour of $\sigma_L^2$, we show that up to time $1$, $\ell^L$ does not see that it is on a torus. Hence, with high probability $\ell^{L,1}(1)^2= \hat{\ell}^{L,1}(1)^2$ and so
$$
\E_0\big[\ell^{L,1}(1)^2\big] \approx \E_0\big[\hat{\ell}^{L,1}(1)^2\big]= 2\frac{\rho_L}{\psi_L^2}\int y^2 \chi^s(dy,dz)+2\int y^2\chi^B(dy,dz) = 2\sigma^2_B +o(1)
$$
as $L\rightarrow \infty$, where the second equality uses (\ref{PPP moment 2}). To make the first equality rigorous, we apply Doob's maximal inequality to the submartingale $|\hat{\ell}^L|^4$. This yields, with a constant $C>0$ which may change from line to line,
$$
\prob_0\bigg[\sup_{0\leq s\leq 1}|\hat{\ell}^L(s)|>\frac{L}{3\psi_L}\bigg] \leq \frac{C\psi_L^4}{L^4}\ \E_0\big[|\hat{\ell}^L(1)|^4\big].
$$
But the calculation in (\ref{bound coord}) shows that the latter expectation is finite, and so
\begin{equation}\label{proba far}
\prob_0\bigg[\sup_{0\leq s\leq 1}|\hat{\ell}^L(s)|>\frac{L}{3\psi_L}\bigg] \leq C\frac{\psi_L^4}{L^4}.
\end{equation}
On the event $\mathcal{E}_L\equiv\big\{\sup_{0\leq s\leq 1}|\hat{\ell}^L(s)|\leq \frac{L}{3\psi_L}\big\}$, the paths of $\ell^L$ and $\hat{\ell}^L$ can be coupled so that $\ell^L(s)=\hat{\ell}^L(s)$ for every $s\in [0,1]$, and since these quantities are bounded for each $L$ we can write
\begin{eqnarray}
\E_0\big[(\ell^{L,1}(1))^2\big]&=&\E_0\big[(\hat{\ell}^{L,1}(1))^2\ \mathbf{1}_{\mathcal{E}_L}\big]+\E_0\big[(\ell^{L,1}(1))^2 \ \mathbf{1}_{\mathcal{E}_L^c}\big] \nonumber\\
&=& \E_0\big[(\hat{\ell}^{L,1}(1))^2\big]-\E_0\big[(\hat{\ell}^{L,1}(1))^2 \ \mathbf{1}_{\mathcal{E}_L^c}\big] +\E_0\big[(\ell^{L,1}(1))^2 \ \mathbf{1}_{\mathcal{E}_L^c}\big].
\label{var ell1}
\end{eqnarray}
By (\ref{proba far}) and the fact that $\ell^L$ evolves on the torus of size $L\psi_L^{-1}$, the last term on the right-hand side of (\ref{var ell1}) is bounded by
$$
C\ \frac{L^2}{\psi_L^2}\times \frac{\psi_L^4}{L^4}=C\ \frac{\psi_L^2}{L^2}\rightarrow 0\qquad \mathrm{as\ }L\rightarrow \infty.
$$
For the second term on the right-hand side of (\ref{var ell1}), let $\hat{s}_L(1)\equiv \sup_{0\leq s\leq 1}|\hat{\ell}^L(s)|$. Using Fubini's theorem on the second line, we have
\begin{eqnarray}
\E_0\big[(\hat{\ell}^{L,1}(1))^2 \ \mathbf{1}_{\mathcal{E}_L^c}\big]&\leq & \E_0\big[\hat{s}_L(1)^2 \ \mathbf{1}_{\mathcal{E}_L^c}\big] \nonumber\\
&=& \int_0^{\infty}\prob_0\Big[\hat{s}_L(1)> \frac{L}{3\psi_L}\vee \sqrt{y}\Big]\ dy \nonumber\\
&=& \frac{L^2}{9\psi_L^2}\ \prob_0\Big[\hat{s}_L(1)> \frac{L}{3\psi_L}\Big] + \int_{\frac{L^2}{9\psi_L^2}}^{\infty} \prob_0\big[\hat{s}_L(1)> \sqrt{y}\big]\ dy.\phantom{AAA}
\label{calcul variance}
\end{eqnarray}
Now, by the argument leading to (\ref{proba far}), $\prob_0[\hat{s}_L(1)>\sqrt{y}]$ is bounded by $Cy^{-2}$ for each $y>0$, where $C$ is a constant independent of $y$. Consequently, the right-hand side of (\ref{calcul variance}) is bounded by
$$
C'\ \frac{\psi_L^2}{L^2}+C\int_{L^2/(9\psi_L^2)}^{\infty}\frac{dy}{y^2}\ \rightarrow \ 0 \qquad \mathrm{as\ }L\rightarrow \infty.
$$
Coming back to (\ref{var ell1}), we can conclude that
$$
\sigma_L^2=2\sigma_B^2 +o(1)\qquad \mathrm{as\ }L\rightarrow \infty.
$$
If we now recall the equality in distribution described at the beginning of the section, we can use Lemma \ref{lemma entrance_levy} applied to $\ell^L$ on the torus of size $L\psi_L^{-1}$ and the entrance time into $B(0,2R^B)$ to write that
\begin{equation}\label{almost hitting}
\lim_{L\rightarrow \infty}\ \sup_{t\geq 0}\ \sup_{A_L\in \GA(L,2)^*}\ \left|\prob_{A_L}\left[T_L> \frac{\rho_L(L/\psi_L)^2\log (L/\psi_L)}{2\pi \sigma_B^2}\ t\right] - e^{-t}\right|=0.
\end{equation}
By the assumption on $|\alpha\log L-\log(\psi_L)|$ introduced just after (\ref{def alpha}) and Lemma \ref{lemm no entrance} applied to $\ell^L$ to bound the probability that $T_L$ lies between $\frac{\rho_LL^2\log (L/\psi_L)}{2\pi \sigma_B^2\psi_L^2}$ and $\frac{(1-\alpha)\rho_LL^2\log L}{2\pi \sigma_B^2\psi_L^2}$, $(a)$ of Proposition \ref{prop gathering} follows from (\ref{almost hitting}).

Let us now turn to the proof of $(b)$. This time, we define $\ell^L$ for every $t\geq 0$ by
$$
\ell^L(t)=\frac{1}{\psi_L}\ \xi^L(2\psi_L^2t).
$$
Similar calculations give, as $L\rightarrow \infty$,
$$
E_0\left[|\ell^L(1)|^2\right]=2\sigma^2_s + 2 b \sigma_B^2 +o(1) \qquad \mathrm{if\ }\rho_L^{-1}\psi_L^2\rightarrow b\in [0,\infty).
$$
and $E_0\left[|\ell^L(1)|^4\right]$ is bounded uniformly in $L$. We can therefore apply Lemma \ref{lemma entrance_levy} to $\ell^L$ as above.$\hfill\square$

\medskip
Having established the time that it takes for two lineages starting from distance $L$ apart to come close enough together that they have a chance to coalesce, we now calculate the additional time required for them to actually do so. We shall have to distinguish between several regimes, depending on whether large or small events prevail in the evolution of the pair of lineages. Our goal in the rest of this section is to prove the following result.

\begin{thm}\label{theo time coal}
For each $L\in \N$, let $t_L$ denote the coalescence time of the pair of lineages under consideration. Then,

\noindent $(a)$ If $\frac{\psi_L^2}{\rho_L}\rightarrow \infty$ as $L\rightarrow \infty$,
$$
\lim_{L\rightarrow \infty}\ \sup_{t\geq 0}\ \sup_{A_L\in \GA(L,2)^*}\left|\ \prob_{A_L}\left[t_L> \frac{(1-\alpha)\rho_LL^2\log L}{2\pi \sigma_B^2\psi_L^2}\ t\right]-e^{-t}\right|=0.
$$

\noindent $(b)$ If $\frac{\psi_L^2}{\rho_L}\rightarrow b\in [0,\infty)$ and $\frac{\psi_L^2\log L}{\rho_L}\rightarrow \infty$ as $L\rightarrow \infty$,
$$
\lim_{L\rightarrow \infty}\ \sup_{t\geq 0}\ \sup_{A_L\in \GA(L,2)^*}\ \left|\ \prob_{A_L}\left[t_L> \frac{(1-\alpha)L^2\log L}{2\pi (\sigma_s^2+b\sigma_B^2)}\ t\right]-e^{-t}\right|=0.
$$

\noindent $(c)$ If $(\frac{\psi_L^4}{\rho_L})_{L\geq 1}$ is bounded or $\frac{L^2\log L}{\rho_L}\rightarrow 0$ as $L\rightarrow \infty$ (and so $\frac{\psi_L^2\log L}{\rho_L}\rightarrow 0$), then
$$
\lim_{L\rightarrow \infty}\ \sup_{t\geq 0}\ \sup_{A_L\in \GA(L,2)^*}\ \left|\ \prob_{A_L}\left[t_L> \frac{L^2\log L}{2\pi \sigma_s^2}\ t\right]-e^{-t}\right|=0.
$$
\end{thm}
The cases $(a)$ and $(b)$ are separated only because the timescales of interest are not of the same order, but the reasons why they hold are identical: in both cases, large jumps are frequent enough that, once the lineages have been gathered at distance $2R^B\psi_L$, they coalesce in a time negligible compared to $T_L$. In contrast, in $(c)$ we assume that the rate at which the lineages are affected by large events is so slow that we have to wait for the lineages to be gathered at distance less than $2R^s$ before they have a chance to coalesce (and they do so in a negligible time compared to $L^2\log L$). If none of the above conditions hold, then the proof of $(c)$ will show that, also in this case, the probability that a large event affects the lineages when they are at distance less than $2R^B\psi_L$ and before a time of order $\mathcal{O}(L^2\log L)$ vanishes as $L$ tends to infinity. However, we are no longer able to describe precisely the limiting behaviour of $t_L$, see Remark \ref{rk last cases}.

Let us first make more precise the sense in which the additional time to coalescence is negligible once the lineages have been gathered at the right distance.
\begin{propn}\label{prop coal time} Let $(\Phi_L)_{L\geq 1}$ be a sequence tending to infinity as $L\rightarrow \infty$.

\noindent $(a)$ If $(\Phi_L)_{L\geq 1}$ is such that $\frac{\rho_L}{\psi_L^2\log \Phi_L}\rightarrow 0$ as $L\rightarrow \infty$, we have
\begin{equation}
\label{short coal 1}\lim_{L\rightarrow \infty}\sup_{A_L}\ \prob_{A_L}\big[t_L> \Phi_L\rho_L\big]=0,
\end{equation}
where the supremum is taken over all samples $A_L= \big\{(\{1\},x_1^L),(\{2\},x_2^L)\big\}$ such that $|x_1^L-x^L_2|\leq 2R^B\psi_L$.

\noindent $(b)$ Under no additional condition, we have
\begin{equation}
\label{short coal 2}\lim_{L\rightarrow \infty}\sup_{A'_L}\ \prob_{A'_L}\big[t_L> \Phi_L\big]=0,
\end{equation}
where the supremum is now taken over all samples $A'_L= \big\{(\{1\},x_1^L),(\{2\},x_2^L)\big\}$ such that $|x_1^L-x^L_2|\leq 2R^s$.
\end{propn}
Taking $\Phi_L=\frac{L^2}{\rho_L\log L}(1\wedge \rho_L \psi_L^{-2})$, the result in $(a)$ shows that when $\frac{\psi_L^2\log L}{\rho_L}\rightarrow \infty$, the coalescence time of two lineages at distance at most $2R^B\psi_L$ is indeed much smaller than $T_L$ (which is of order $L^2\log L\times \big(1\wedge \rho_L\psi_L^{-2}\big)$ by Proposition \ref{prop gathering}).

\medskip
\noindent \emph{Proof of Proposition \ref{prop coal time}: }Recall that for each $L\in \N$, we defined $X^L$ as the difference between the locations of the lineages $\xi^L_1$ and $\xi^L_2$ on the torus $\T(L)$. In the following, if both lineages are affected by the same event, we shall consider that $X^L$ hits $0$ but the number of lineages remains equal to $2$, which means that they can separate again later (if the measures $\nu^s_r$ and $\nu^B_r$ are not all the point mass at $1$). However, it is the first time at which such an event occurs which will be of interest, and we keep the notation $t_L$ to denote this time. As we already noticed, $X^L$ behaves like $\{\xi^L(2t),\ t\geq 0\}$ outside $B\big(0,2R^B\psi_L\big)$, whereas inside the ball it can hit $0$ owing to reproduction events affecting both lineages $\xi^L_1$ and $\xi^L_2$.

\noindent \textbf{Case $(a)$.} For each $L\in \N$, set $q^L_0=Q^L_0\equiv 0$ and for every $i\geq 1$,
$$
Q^L_i\equiv \inf \Big\{t> q^L_{i-1}:\ X^L(t)\notin B\Big(0,\frac{7}{4}R^B\psi_L\Big) \Big\}
$$
and
$$
q^L_i\equiv \inf \Big\{t>Q^L_i:\ X^L(t)\in B\Big(0,\frac{3}{2}R^B \psi_L\Big) \Big\},
$$
with the convention that $\inf \ \emptyset =+\infty$. We shall use the following lemmas, which will enable us to describe how $X^L$ wanders around in $\T(L)$, independently of whether it ever hits $0$ or not.

\begin{lemma}\label{lemm excursion}There exist a function $g:\R_+\rightarrow \R_+$ vanishing at infinity, $C_q>0$, $u_q>1$ and $L_q\in \N$ such that for every $L\geq L_q$ and
$u\geq u_q$,
$$
\sup_{x\in B(0,4R^B)\setminus B(0,(7/4)R^B)} \prob_{\psi_Lx}\big[q^L_1>\rho_Lu\big] \leq g(u) \qquad  \mathrm{if\ }\rho_L=\mathcal{O}(\psi_L^2),
$$
$$
\sup_{x\in B(0,4R^B)\setminus B(0,(7/4)R^B)} \prob_{\psi_Lx}\big[q^L_1>\psi_L^2u\big] \leq \frac{C_q}{\log u} \qquad \mathrm{if\ }\rho_L^{-1}\psi_L^2\rightarrow 0.
$$
\end{lemma}
Lemma \ref{lemm excursion} will give us good control of the probability of a long excursion outside $B(0,(3/2)R^B\psi_L)$.

\begin{lemma}\label{lemm time in B} Suppose that
\begin{equation}\label{condition nu}
\mathrm{Leb}\big(\big\{r\in [0,R^B]:\ \nu^B_r \notin \{\delta_0,\delta_1\}\big\}\big)>0.
\end{equation}
Then, there exists a constant $C_Q<\infty$ such that for each $L\geq 1$,
$$
\sup_{x\in B(0,(3/2)R^B)}\frac{1}{\rho_L}\ \E_{\psi_Lx}\big[Q_1^L\big]< C_Q.
$$
\end{lemma}
Condition (\ref{condition nu}) guarantees that, whenever $X^L$ hits $0$, it has a chance not to remain stuck at this value for all times. Lemma \ref{lemm time in B} then tells us that $X^L$ starting within $B((3/2)R^B\psi_L)$ needs an average time of order $\mathcal{O}(\rho_L)$ to reach distance $(7/4)R^B\psi_L$ from the origin.

\begin{lemma}\label{lemm proba coal}Suppose that $\rho_L\psi_L^{-2}$ remains bounded as $L\rightarrow \infty$. Then, there exists $\theta_1 \in (0,1)$ such that for every $L\geq 1$,
\begin{equation}\label{hitting prob 1}
\inf_{x\in B(0,(3/2)R^B)}\prob_{\psi_Lx}\big[X^L \mathrm{\ hits\ }0 \mathrm{\ before\ leaving\ }B\big(0,(7/4)R^B\psi_L\big)\big]\geq \theta_1.
\end{equation}
If $\liminf_{L\rightarrow \infty} \rho_L^{-1}\psi_L^2= 0$, there exist $\theta_2 \in (0,1)$ and $\theta_3>0$ such that \setlength\arraycolsep{1pt}
\begin{eqnarray}
\inf_{x\in B(0,(3/2)R^B)}\prob_{\psi_Lx}\big[ X^L \mathrm{\ hits\ }0 \mathrm{\ before\ leaving}& B&\big(0,(7/4)R^B\psi_L\big)\big]\nonumber\\
& \geq & \theta_2\bigg(1-\exp\Big\{-\theta_3\frac{\psi_L^2}{\rho_L}\Big\}\bigg).\label{hitting prob 2}
\end{eqnarray}
\end{lemma}
The proofs of these lemmas are given in Appendix \ref{appendix 2}.

The following technique is inspired by that used in Cox \& Durrett~(2002) and Z\"ahle et al.~(2005), \nocite{cox/durrett:2002} \nocite{zahle/cox/durrett:2005} although the motions of the lineages and the mechanism of coalescence here are more complex and require slightly more work. Our plan is first to find a good lower bound on the number of times the lineages meet at distance less than $(3/2)R^B\psi_L$ (and then separate again) before time $\Phi_L\rho_L$. In a second step, we use the estimates on the probability that during such a gathering the lineages merge before separating again derived in Lemma \ref{lemm proba coal}, and obtain that coalescence does occur before $\Phi_L\rho_L$ with probability tending to $1$. For the sake of clarity, we show (\ref{short coal 1}) in the case where $\rho_L\psi_L^{-2}$ remains bounded, and then comment on how to adapt the arguments in the general case.

Assume first that Condition~(\ref{condition nu}) holds. Recall the definition of $Q_i^L$ and $q_i^L$ given above, and define $k_L$ by
$$
k_L\equiv \max\big\{n:\ Q_n^L \leq \Phi_L\rho_L\big\}.
$$
By Lemma \ref{lemm proba coal}, there exists a positive constant $\theta_1$ such that for every $L\geq 1$ and $x\in B(0,(3/2)R^B\psi_L)$,
$$
\prob_x\big[X^L \mathrm{\ hits\ }0 \mathrm{\ before\ leaving\ }B\big(0,(7/4)R^B\psi_L\big)\big]\geq \theta_1.
$$
Hence, for every $x\in B\big(0,2R^B\psi_L\big)$, we have
\begin{equation}\label{eq k_L}
\prob_x\big[t_L>\Phi_L\rho_L\big]\leq \prob_x\big[t_L>Q_{k_L}^L\big]\leq \E_x\big[\big(1-\theta_1\big)^{k_L}\big].
\end{equation}

Let us fix $x\in B\big(0,2R^B\psi_L\big)$ and show that $k_L \rightarrow \infty$ as $L\rightarrow \infty$, in $\prob_x$-probability. The fact that the bounds obtained below do not depend on $x\in B\big(0,2R^B\psi_L\big)$ will then give us the desired uniformity. Let $M\in \N$. We have
\setlength\arraycolsep{1pt}
\begin{eqnarray}\prob_x\big[k_L<M\big]&=&\prob_x\big[Q^L_M> \Phi_L\rho_L \big] \nonumber\\
& = &\prob_x\bigg[\sum_{i=1}^M(Q_i^L-q^L_{i-1})+\sum_{i=1}^{M-1}(q_i^L-Q_i^L)> \Phi_L\rho_L\bigg]\nonumber \\
&\leq &\sum_{i=1}^M\prob_x\bigg[Q_i^L-q_{i-1}^L>\frac{\Phi_L\rho_L}{2M}\bigg] + \sum_{i=1}^{M-1}\prob_x\bigg[q_i^L-Q_i^L> \frac{\Phi_L\rho_L}{2(M-1)}\bigg],\phantom{AAAA} \label{ineg k_L}
\end{eqnarray}
where the last inequality uses the fact that at least one of the $2M-1$ terms of the sums on the second line must be larger than a fraction $(2M-1)^{-1}$ of the total time. Now, using the Markov inequality, the strong Markov property at time $q_{i-1}^L$ and then Lemma~\ref{lemm time in B}, we can write for each $i$
\begin{eqnarray*}\prob_x\bigg[Q_i^L-q_{i-1}^L>\frac{\Phi_L\rho_L}{2M}\bigg] &\leq & \frac{2M}{\Phi_L\rho_L}\ \E_x \big[Q_i^L-q_{i-1}^L\big] \\
& \leq & \frac{2M}{\Phi_L\rho_L}\sup_{y\in B(0,(3/2)R^B)}\E_{\psi_Ly}\big[Q_1^L\big]\\ &\leq & \frac{2MC_Q}{\Phi_L}.
\end{eqnarray*}
If we now apply the strong Markov property to $X^L$ at time $Q^L_i$ and use Lemma \ref{lemm excursion} together with the fact that $X^L(Q_i^L)\in B(0,4R^B\psi_L)$ with probability one, we obtain for each $i$, and $L$ large enough
$$
\prob_x\bigg[q_i^L-Q^L_i>\frac{\Phi_L\rho_L}{2(M-1)}\bigg]\leq g\bigg(\frac{\Phi_L}{2(M-1)}\bigg).
$$
Coming back to (\ref{ineg k_L}), we arrive at
$$
\prob_x\big[k_L< M\big] \leq \frac{2M^2C_Q}{\Phi_L} + (M-1) g\bigg(\frac{\Phi_L}{2(M-1)}\bigg) \rightarrow  0, \qquad \mathrm{as \ }L\rightarrow \infty.
$$
To complete the proof of $(a)$ when Condition~(\ref{condition nu}) holds and $\rho_L\psi_L^{-2}$ remains bounded, let $\e>0$ and fix $M = M(\e)\in \N$ such that
$$
(1-\theta_1)^M <\e.
$$
Splitting the expectation in (\ref{eq k_L}) into the integral over $\{k_L\geq M\}$ and $\{k_L< M\}$ yields
$$
\limsup_{L\rightarrow \infty}\sup_{x\in B(0,2R^B\psi_L)}\prob_x\big[t_L>\Phi_L\rho_L\big]\leq \e + \limsup_{L\rightarrow \infty}\sup_{x\in B(0,2R^B\psi_L)}\prob_x\big[k_L<M\big] = \e,
$$
and since $\e$ was arbitrary, the desired result follows.

When Condition~(\ref{condition nu}) is fulfilled but $\rho_L\psi_L^{-2}$ is unbounded as $L\rightarrow \infty$, we can apply the same technique to obtain (\ref{short coal 1}). This time, using the second result of Lemma \ref{lemm proba coal} we can write as in (\ref{eq k_L}) that, for every $x\in B\big(0,2R^B\psi_L\big)$,
$$
\prob_x\big[t_L>\Phi_L\rho_L\big]\leq \E_x\bigg[\bigg(1-\theta_2\Big(1-\exp\Big\{-\theta_3
\frac{\psi_L^2}{\rho_L}\Big\}\Big)\bigg)^{k_L}\bigg].
$$
The same arguments as above (using the second part of Lemma \ref{lemm excursion}) yield, for $L$ large enough,
\begin{eqnarray*}
\sup_{x\in B(0,2R^B\psi_L)}\prob_x\bigg[k_L<M \frac{\rho_L}{\psi_L^2}\bigg]& \leq &\frac{2C_QM^2\rho_L^2}{\psi_L^4\Phi_L}  + \frac{C_qM\rho_L}{\psi_L^2\log(\Phi_L/2M)},
\end{eqnarray*}
which tends to $0$ as $L$ tends to infinity by our assumption of $(\Phi_L)_{L\geq 1}$. We conclude in the same manner, using the fact that when $\psi_L^2/\rho_L\rightarrow \infty$,
$$
\bigg(1-\theta_2\Big(1-\exp\Big\{-\theta_3 \frac{\psi_L^2}{\rho_L}\Big\}\Big)\bigg)^{M\rho_L/\psi_L^2}\sim e^{-\theta_2\theta_3 M}.
$$

Let us finish the proof of $(a)$ by removing the assumption (\ref{condition nu}). In the preceding proof, the main idea is that each time $X^L$ passes through $B\big(0,(3/2)R^B\psi_L\big)$, the two lineages have an opportunity to try to coalesce and their success probability is bounded from below by the quantity obtained in Lemma \ref{lemm proba coal}. However, if we do not assume that (\ref{condition nu}) holds, $X^L$ may become stuck at $0$ once it has hit it, and so the number $k_L$ of such sojourns in $B\big(0,(3/2)R^B\psi_L\big)$ may be finite. This makes our arguments break down. Nevertheless, $X^L$ can only hit $0$ through a coalescence event, and so this issue is merely an artefact of the technique of the proof. To overcome it, let us increase the rate of reproduction events by a factor $2$, but divide each probability to be affected by 2. Overall, coalescence will take a longer time in this new setting, but the motions of the lineages before their coalescence time will remain identical in distribution.

More precisely, assume that (\ref{condition nu}) does not hold. Define $\hat{\Pi}^B_L$ as a Poisson point process on $\R\times \T(L)\times (0,\infty)$, independent of $\Pi^s_L$ and $\Pi^B_L$ and with intensity measure $2(\rho_L\psi_L^2)^{-1} dt\otimes dx\otimes \mu^B(dr)$, and for each $r>0$ such that $\nu_r^B=\delta_1$, set $\hat{\nu}^B_r \equiv \delta_{1/2}$. Let also $\hat{\Pi}^s_L$ be a Poisson point process with the same distribution as $\Pi^s_L$ and independent of all the other point processes. Call $\hat{X}^L$  the process defined in the same manner as $X^L$ but with $\Pi^B_L$ (resp., $\Pi^s_L$, $\nu^B_r$) replaced by $\hat{\Pi}^B_L$ (resp., $\hat{\Pi}^s_L$, $\hat{\nu}^B_r$). By computing the intensity of the jumps of a single lineage, one can observe that it is equal to
$$
dt\otimes \bigg(\frac{2}{\rho_L}\int_{|x|/2}^{R^B}\frac{L_r(x)}{2\pi r^2}\ \mathbf{1}_{\{\nu_r^B=\delta_1\}} \mu^B(dr)d(\psi_Lx) +\int_{|x|/2}^{R^s}\int_0^1\frac{L_r(x)}{\pi r^2}\ u\ \nu^s_r(du)\mu^s(dr)dx\bigg),
$$
which is precisely that of $\xi^L$. Here, $L_r(x)$ stands for the volume of $B(0,r)\cap B(x,r)$. If we now compute the coalescence rate of two lineages at distance $z\in [0,2R^B\psi_L]$, we obtain the same term due to small events for $X^L$ and $\hat{X}^L$, to which is added the respective contributions of large events
$$
\frac{1}{\rho_L}\int_{z/2}^{R^B}L_r(z)\mathbf{1}_{\{\nu_r^B = \delta_1\}}\mu^B(dr) \qquad \mathrm{and}\qquad \frac{1}{2\rho_L}\int_{z/2}^{R^B}L_r(z)\mathbf{1}_{\{\nu_r^B= \delta_1\}}\mu^B(dr).
$$
Hence, the evolutions of both processes follow the same law outside $B(0,2R^B\psi_L)$, the contribution of large events whose area encompasses only one of the two lineages is identical even within $B(0,2R^B\psi_L)$, and coalescence occurs at a higher rate for $X^L$ than for $\hat{X}^L$. This gives us for every $L\geq 1$ and $x\in \T(L)$,
$$
\prob_x\big[t_L>\Phi_L\rho_L \big] \leq \prob_x\big[\hat{t}_L>\Phi_L\rho_L \big],
$$
where $\hat{t}_L$ is defined in an obvious manner. But Condition~(\ref{condition nu}) holds for $\hat{X}^L$, and so we can use the result obtained in the previous paragraph to complete the proof of $(a)$ when (\ref{condition nu}) does not hold.

\medskip
\noindent \textbf{Case $(b)$.} The arguments are essentially the same. First of all, since we assumed that $\rho_L$ grows to infinity as $L\rightarrow \infty$, and because
$$
\prob_x\big[t_L>\Phi_L\big]\leq \prob_x\big[t_L>\Phi'_L\big]
$$
whenever $\Phi_L\geq \Phi'_L$, we can restrict our attention to sequences $(\Phi_L)_{L\geq 1}$ such that $\rho_L^{-1}\Phi_L\rightarrow 0$ as $L\rightarrow \infty$. Let $\mathcal{E}_L$ denote the event that no large events affected any of the lineages before time $\Phi_L$. Let $\theta_{\mathrm{max}} \in (0,\infty)$ be such that the maximal rate at which at least one of the two lineages of the sample is affected by a large event is less than $\theta_{\mathrm{max}} \rho_L^{-1}$ (recall that the total rate at which at least one of two lineages is affected is smaller than twice the corresponding rate for a single lineage, which is finite and independent of the location of the lineage). For each $L\in \N$, define $e_L$ as an exponential random variable, with parameter $\theta_{\mathrm{max}} \rho_L^{-1}$. By our assumption on $\Phi_L$, we can write $$
\prob_x [\mathcal{E}_L^c]\leq \prob[e_L\leq \Phi_L] = 1-\exp\bigg\{-\frac{\theta_{\mathrm{max}} \Phi_L}{\rho_L}\bigg\}\rightarrow 0, \quad \mathrm{as\ }L\rightarrow \infty.
$$
The distribution of the process $X^L$ up to the first time at which it is affected by a large event is equal to that of $\tilde{X}^L$ (defined as the process experiencing only small events) up to the random time $e(\tilde{X}^L)$, so that if $\rho_L^{-1}\theta_{B,L}(x)$ is the rate at which at least one of two lineages at separation $x\in \T(L)$ is affected by a large event, then for each $t\geq 0$ and $y\in \T(L)$
$$
\prob_y\big[e(\tilde{X}^L)>t\big]=\E_y\bigg[\exp\bigg\{-\int_0^t \frac{\theta_{B,L}\big(\tilde{X}^L(s)\big)}{\rho_L}\ ds\bigg\}\bigg].
$$
By the definition of $\theta_{\mathrm{max}}$, for each $L\in \N$ the variable $e_L$ is stochastically bounded by $e(\tilde{X}^L)$. Consequently, if $\tilde{t}_L$ denotes the coalescence time associated to $\tilde{X}^L$ (or, more precisely, to the model where lineages are affected only by small events), we have for each $x\in B(0,2R^s)$
\begin{eqnarray*}
\prob_x\big[t_L\geq \Phi_L\big]&\leq & \prob_x\big[t_L\geq \Phi_L;\ \mathcal{E}_L\big] + \prob_x\big[\mathcal{E}_L^c\big]\\
&\leq & \prob_x\big[\tilde{t}_L\geq \Phi_L\big]+o(1) \quad \mathrm{as\ }L\rightarrow \infty,
\end{eqnarray*}
where the remaining terms converge to $0$ uniformly in $x\in \T(L)$. Then, an easy modification of the proof of $(a)$ with ``$\psi_L=\rho_L=1$'' yields the desired result and completes the proof of Proposition \ref{prop coal time}. $\hfill\square$

\medskip
We can now turn to the proof of Theorem \ref{theo time coal}.

\noindent \emph{Proof of Theorem \ref{theo time coal}: }

\noindent \textbf{Cases $(a)$ and $(b)$. }For $(a)$, let us define $\Phi_L$ for each $L\in \N$ by
$$
\Phi_L= \frac{\rho_LL^2}{\psi_L^2\log L}.
$$
Let $t>0$ and $(A_L)_{L\geq 1}$ be such that $A_L\in \GA(L,2)^*$ for each $L\in \N$. Introducing the time $T_L$ needed for the two lineages of the sample to come at distance less than $2R^B\psi_L$, we can write \setlength\arraycolsep{1pt}
\begin{eqnarray}
\prob_{A_L}&\bigg[&t_L > \frac{(1-\alpha)\rho_LL^2\log L}{2\pi \sigma^2_B\psi_L^2}\ t\bigg]\nonumber \\
& =& \prob_{A_L}\bigg[t_L>\frac{(1-\alpha)\rho_LL^2\log L}{2\pi \sigma^2_B\psi_L^2}\ t;\ T_L>\frac{(1-\alpha)\rho_LL^2\log L}{2\pi \sigma^2_B\psi_L^2}\ t- \Phi_L\bigg] \label{proof coal 1}
\\ & & + \prob_{A_L}\bigg[t_L>\frac{(1-\alpha)\rho_LL^2\log L}{2\pi \sigma^2_B\psi_L^2}\ t;\ T_L\leq \frac{(1-\alpha)\rho_LL^2\log L}{2\pi \sigma^2_B\psi_L^2}\ t- \Phi_L\bigg].\phantom{AAA}\label{proof coal 2}
\end{eqnarray}
Using the strong Markov property at time $T_L$ and the uniform convergence derived in Proposition \ref{prop coal time}$(a)$, we obtain that the expression in (\ref{proof coal 2}) tends to $0$ as $L\rightarrow \infty$ independently of the choice of $t>0$ and $(A_L)_{L\in\N}$. For (\ref{proof coal 1}), note that \setlength\arraycolsep{1pt}
\begin{eqnarray}
\bigg|\ \prob_{A_L}& \bigg[& t_L>\frac{(1-\alpha)\rho_LL^2\log L}{2\pi \sigma^2_B\psi_L^2}\ t;\ T_L>\frac{(1-\alpha)\rho_LL^2\log L}{2\pi \sigma^2_B\psi_L^2}\ t- \Phi_L\bigg] \nonumber \\
& & - \prob_{A_L}\bigg[T_L>\frac{(1-\alpha)\rho_LL^2\log L}{2\pi \sigma^2_B\psi_L^2}\ t\bigg]\bigg|\nonumber \\
&\leq & \prob_{A_L}\bigg[\frac{(1-\alpha)\rho_LL^2\log L}{2\pi \sigma^2_B\psi_L^2}\ t- \Phi_L\leq T_L \leq \frac{(1-\alpha)\rho_LL^2\log L}{2\pi \sigma^2_B\psi_L^2}\ t\bigg].\phantom{AA}\label{bound T_L}
\end{eqnarray}
Since $X^L$ (defined at the beginning of Section \ref{section coal}) has the same law as $\{\xi^L(2t),t\geq 0\}$ until the random time $T_L$, we can bound the quantity in (\ref{bound T_L}) by working directly with the latter process. In order to apply Lemma~\ref{lemm no entrance} to $\big\{\psi_L^{-1}\xi^L(2\rho_Lt),t\geq 0\big\}$, with
$$
U_L=\frac{(1-\alpha)L^2\log L}{4\pi \sigma^2_B\psi_L^2},\quad u_L=\frac{\Phi_L}{2\rho_L}=\frac{L^2}{2\psi_L^2\log L}\quad \mathrm{and} \ R=2R^B,
$$
we need to check that $U_L\psi_L^2L^{-2}\rightarrow \infty$ and $u_L\leq \frac{L^2}{\psi_L^2\sqrt{\log (L/\psi_L)}}$ (recall that this process evolves on the torus of size $\psi_L^{-1}L$.) Both conditions are fulfilled here, and so by Lemma~\ref{lemm no entrance}, the right-hand side of (\ref{bound T_L}) is bounded by
$$
\frac{C\Phi_L\psi_L^2}{\rho_LL^2}=\frac{C}{\log L}\rightarrow 0\qquad \mathrm{as\ }L\rightarrow \infty.
$$
Hence, coming back to (\ref{proof coal 1}), we can use the result of Proposition \ref{prop gathering} and the uniformity in $t>0$ and  $(A_L)_{L\geq 1}$ of our estimates to obtain
$$
\lim_{L\rightarrow \infty}\sup_{t\geq 0}\sup_{A_L\in \GA(L,2)^*}\left|\ \prob_{A_L}\bigg[t_L>\frac{(1-\alpha)\rho_LL^2\log L}{2\pi \sigma^2_B\psi_L^2}\ t\bigg] - e^{-t}\right|=0.
$$

The proof of $(b)$ follows exactly the same lines, with $\Phi_L\equiv L^2(\log L)^{-1}$ and Lemma~\ref{lemm no entrance} applied to $\psi_L^{-1}\xi^L(2\psi_L^2\cdot)$.

\medskip
\noindent \textbf{Case $(c)$. }In contrast with the two previous cases, where coalescence in the limit is due to large events only, here the pair of lineages can coalesce only through a small event. To see this, let us define $T^*_L$ as the first time at which the two lineages (indexed by $L$) come at distance less than $2R^s$ from each other, and $\tau_L$ as the first time at which at least one of them is affected by a large event while they are at distance less than $2R^B\psi_L$ (i.e., while $X^L\in B(0,2R^B\psi_L)$). Note that for each $L$, $T^*_L$ and $\tau_L$ are stopping times with respect to the filtration $\{\mathcal{F}_t,\ t\geq 0\}$ associated to $\Pi_L^s\cup \Pi_L^B$ as we trace backwards in time. In addition, define $\tilde{T}_L^*$ as the entrance time of $\xi^L$ into $B(0,2R^s)$ and $\tilde{\tau}_L$ as the first time $\xi^L$ makes a jump of size $\mathcal{O}(\psi_L)$ while it is lying in $B(0,2R^B\psi_L)$. These two random times are stopping times with respect to the filtration $\{\tilde{\mathcal{F}}_t,\ t\geq 0\}$ associated to $\xi^L$. We claim that for each $L\in \N$,
\begin{equation}\label{claim identity}
\big\{X^L(t),\ t<\tau_L \wedge T^*_L\big\} \stackrel{(d)}{=} \big\{\xi^L(2t),\ 2t<\tilde{\tau}_L \wedge \tilde{T}^*_L\big\},
\end{equation}
where the notation $\stackrel{(d)}{=}$ refers to equality in distribution. Indeed, as long as $X^L$ has not entered $B(0,2R^s)$ and no large event has affected it while it lay in $B(0,2R^B\psi_L)$, coalescence events are impossible and the rates and distributions of the jumps of both processes are identical. We cannot include the terminal times in (\ref{claim identity}) since the values of the processes will differ if $\tau_L\wedge T^*_L= \tau_L$ and the corresponding event is a coalescence, but since $X^L$ and $\xi^L$ are jump processes with finite rates, we can easily see that the event $\{\tau_L\wedge T^*_L= \tau_L\}$ (resp., $\tilde{\tau}_L\wedge \tilde{T}^*_L= \tilde{\tau}_L$) is $\mathcal{F}_{(\tau_L\wedge T^*_L)-}$ (resp., $\tilde{\mathcal{F}}_{(\tilde{\tau}_L\wedge \tilde{T}^*_L)-}$) -measurable. Hence, for each $L\in \N$, $A = \wp_2(x_1,x_2)$ and $x\equiv x_1-x_2\in \T(L)$, we have
\begin{equation}\label{equality proba tau}
\prob_A\big[\tau_L<T_L^*\big]=\prob_x\big[\tilde{\tau}_L<\tilde{T}_L^*\big].
\end{equation}
Let us now bound the right-hand side of (\ref{equality proba tau}) under the assumption that $(\rho_L^{-1}\psi_L^4)_{L\in \IN}$ is bounded. Analogous computations to those in the proof of Proposition \ref{prop gathering} show that $\{\xi^L(2t),t\geq 0\}$ itself satisfies Assumption~\ref{assumptions for levy processes} with $\sigma_L^2=2\sigma_s^2+o(1)$ as $L\rightarrow \infty$. Hence, Lemma \ref{lemma entrance_levy} applied with $d_L=2R^s$ gives us
\begin{equation}\label{coal small}
\lim_{L\rightarrow \infty}\ \sup_{t\geq 0}\ \sup_{x_L\in \Gamma(L,1)}\bigg|\ \prob_{x_L}\bigg[\tilde{T}_L^*>\frac{L^2\log L}{2\pi \sigma_s^2}\ t\bigg]-e^{-t}\bigg|=0.
\end{equation}
Let $\theta_{\mathrm{max}}<\infty$ be such that for every $L\in \N$, the rate at which $\xi^L$ makes a jump of size $\mathcal{O}(\psi_L)$ is bounded by $\theta_{\mathrm{max}}/\rho_L$. Fixing $\e>0$ and $K>0$ such that $e^{-2\pi \sigma_s^2 K}< \e$, we have for $L$ large enough and any sequence $(x_L)_{L\geq 1}$ such that $x_L\in \Gamma(L,1)$ for every $L$: \setlength\arraycolsep{1pt}
\begin{eqnarray}
\prob_{x_L}\big[&\tilde{\tau}_L&<\tilde{T}_L^*\big]\nonumber\\
& = & \prob_{x_L}\big[\tilde{\tau}_L<\tilde{T}_L^*\leq KL^2\log L \big]+\prob_{x_L}\big[\tilde{\tau}_L<\tilde{T}_L^*\ ;\ \tilde{T}_L^*>KL^2\log L\big]\nonumber \\
&\leq & \prob_{x_L}\big[\tilde{\tau}_L< KL^2\log L \big]+\prob_{x_L}\big[\tilde{T}_L^*>KL^2\log L\big]\nonumber \\
& \leq & \E_{x_L}\bigg[1-\exp\Big\{-\frac{\theta_{\mathrm{max}}}{\rho_L}\int_0^{KL^2\log L}\mathbf{1}_{B(0,2R^B\psi_L)}\big(\xi^L(2s)\big)ds\Big\}\bigg]+ \e. \label{eq large event}
\end{eqnarray}
Splitting the integral below into the sum $\int_0^{\psi_L^2\sqrt{\log L}}+ \int_{\psi_L^2\sqrt{\log L}}^{L^2/\sqrt{\log L}}+ \int_{L^2/\sqrt{\log L}}^{L^2\sqrt{\log L}}+ \int_{L^2\sqrt{\log L}}^{KL^2\log L}$ and using the four results of Lemma \ref{lemma local TCL}, there exists $L_0\in \N$, and $a_1,a_2>0$ independent of $L$, $(x_L)_{L\geq 1}$ and $K>0$, such that for every $L\geq L_0$,
$$
\E_{x_L}\bigg[\int_0^{KL^2\log L}\mathbf{1}_{B(0,2R^B\psi_L)} \big(\xi^L(2s)\big)ds\bigg]\leq (a_1+a_2K)\psi_L^2 \log L.
$$
Hence, the first term on the right-hand side of (\ref{eq large event}) is bounded by
$$
\E_{x_L}\bigg[\frac{\theta_{\mathrm{max}}}{\rho_L}\int_0^{KL^2\log L}\mathbf{1}_{B(0,2R^B\psi_L)}\big(\xi^L(2s)\big)ds\bigg] \leq \theta_{\mathrm{max}}(a_1+a_2K)\frac{\psi_L^2\log L}{\rho_L},
$$
which tends to $0$ as $L\rightarrow \infty$, independently of the sequence $(x_L)_{L\geq 1}$ considered. As $\e$ in (\ref{eq large event}) is arbitrary, we can conclude that $$
\lim_{L\rightarrow \infty}\sup_{x_L\in \Gamma(L,1)}\prob_{x_L}\big[\tilde{\tau}_L<\tilde{T}_L^*\big]=0,
$$
and by (\ref{equality proba tau}), the same result holds for $X^L$ and any sequence $(A_L)_{L\in \N}$ such that $A_L\in \GA(L,2)^*$ for every $L$. In words, we have obtained that with probability tending to $1$, any pair of lineages starting at distance $\mathcal{O}(L)$ from each other gather at distance $2R^s$ before having a chance to coalesce through a large reproduction event. By using the same method as in $(a)$ but this time with the result of Proposition \ref{prop coal time} $(b)$ and with Proposition \ref{prop gathering} replaced by (\ref{coal small}), we obtain the desired conclusion under the assumption that $(\rho_L^{-1}\psi_L^4)_{L\in \IN}$ is bounded.

When $\rho_L\gg L^2\log L$, with probability increasing to $1$ no large events at all affect any of the lineages by the time they are gathered at distance $2R^s$ by small events. The result then follows from the same arguments, with $\xi^L$ replaced by the motion of a single lineage subject to only small reproduction events. $\hfill\square$ \medskip

\begin{remark}\label{rk last cases} Let us comment on the cases not covered by the theorem, that is $\psi_L^4\gg \rho_L$, $\rho_L$ is of order at most $L^2\log L$ and $\rho_L^{-1}\psi_L^2\log L$ has a finite limit (possibly $0$). When the latter limit is positive, from the results obtained so far coalescence events due to small and to large reproduction events occur on the same timescale and depend on the precise paths of the two lineages. Therefore, we do not expect $t_L$ to be exponentially distributed (with a deterministic parameter). When $\rho_L^{-1}\psi_L^2\log L$ tends to $0$, the same reasoning as in the proof of $(c)$ gives us that the probability that a large reproduction event causes the two lineages to coalesce before a time of order $L^2\log L$ vanishes as $L\rightarrow \infty$. However, $X^L$ does not satisfy the conditions of Section \ref{levy processes} (Assumption~\ref{assumptions for levy processes}) as it does when the assumptions of $(c)$ hold. Using instead $\ell^L\equiv \psi_L^{-1}X^L(\psi_L^2\cdot)$, the time needed for the lineages to come at distance less than $2R^s$ translates into $T(\ell^L,2R^s/\psi_L)$, which is not covered by Lemma \ref{lemma entrance_levy} and requires estimates of the entrance time of the jump process into a ball of shrinking radius, which we have been unable to obtain.
\end{remark}

\subsection{Convergence to Kingman's coalescent}
\label{subsection kingman} To complete the proof of Theorem~\ref{result alpha<1}, we now turn to the genealogy of a finite sample, starting at distance $\mathcal{O}(L)$ from each other on $\T(L)$.

We can already see from our analysis for a single pair of lineages that our spatial $\Lambda$-coalescent is similar in several respects to the coalescing random walks dual to the two-dimensional voter and stepping-stone models with short-range interactions (see e.g. \nocite{cox/griffeath:1986} \nocite{cox/griffeath:1990} Cox \& Griffeath~1986, 1990 for a study on $\Z^2$, and Cox~1989 or Z\"ahle et al.~2005 \nocite{cox:1989} \nocite{zahle/cox/durrett:2005} for examples on the torii $\T(L)\cap \Z^2$). It will therefore be no surprise that the analogy carries over to larger samples. In most of the papers cited above, the authors are interested in the sequence of processes giving the number of blocks in the ancestral partition. They show that, when the initial distance between the lineages grows to infinity, the finite-dimensional distributions of these counting processes converge to those of a pure death process corresponding to a time-change of the number of blocks of Kingman's coalescent. In Cox \& Griffeath~(1990), \nocite{cox/griffeath:1990} more elaborate arguments yield the convergence of the finite-dimensional distributions of the unlabelled genealogical processes to those of Kingman's coalescent. Instead of adding a new instance of such proofs to the literature, we shall simply explain why the same method applies to our case. This will also enable us to prove the tightness of the unlabelled genealogical processes.

\medskip
\noindent \emph{Proof of Theorem~\ref{result alpha<1}: } \textbf{(i) Convergence of the finite-dimensional distributions.}

We follow here the proofs in Cox \& Griffeath~(1986) \nocite{cox/griffeath:1986} (for the number of blocks of the ancestral partition) and Cox \& Griffeath~(1990) \nocite{cox/griffeath:1990} (for the unlabelled genealogical process of a system of coalescing simple random walks on $\Z^2$). Notice that, since we work on the torii $\T(L)$, our rescaling of time differs from Cox and Griffeath's. Another significant difference is the fact that, in their model, lineages move independently of each other until the first time two of them are on the same site, upon which they coalesce instantaneously. In our setting, the movements of lineages are defined from the same Poisson point processes, and two lineages having reached a distance that enables them to coalesce can separate again without coalescing.

Despite these differences, Lemma \ref{lemm lineages far} below shows that a key ingredient of their proof is still valid here: at the time when two lineages coalesce, the others are at distance $\mathcal{O}(L)$ from each other and from the coalescing pair. To state this result, we need some notation. Let $\tau_{ij}$ be the first time lineages  $i$ and $j$ come within distance less than $2R^B\psi_L$ (resp., $2R^s$) if $\rho_L\ll \psi_L^2\log L$ (resp., $\rho_L\gg \psi_L^2\log L$) and $\tau$ be the minimum of the $\tau_{ij}$'s over all pairs considered. Let also $\tau^*_{ij}$ be the coalescence time of the ancestral lines of $i$ and $j$, and $\tau^*$ be the minimum of the $\tau^*_{ij}$ over all lineages considered. Finally, for each $i$ we shall denote the motion in $\T(L)$ of the block containing $i$ by $\xi_i^L$.
\begin{lemma}\label{lemm lineages far}Under the conditions of Theorem \ref{result alpha<1}, we have
\begin{eqnarray}
\lim_{L\rightarrow\infty}\sup_{A_L\in \GA(L,4)^*}\prob_{A_L}\bigg[\tau^*=\tau^*_{12}\ ;\ |\xi_1^L(\tau^*)-\xi_3^L(\tau^*)|\leq \frac{L}{\log L}\bigg] = 0, \label{lineage distribution 1} \\
\lim_{L\rightarrow\infty}\sup_{A_L\in \GA(L,4)^*}\prob_{A_L}\bigg[\tau^*=\tau^*_{12}\ ;\ |\xi_3^L(\tau^*)-\xi_4^L(\tau^*)|\leq \frac{L}{\log L}\bigg] = 0.
\label{lineage distribution 2}
\end{eqnarray}
\end{lemma}
The proof of Lemma \ref{lemm lineages far} is deferred to Appendix \ref{appendix 2}.

The other ingredients required to apply Cox and Griffeath's techniques are a control on the probability of ``collision'' for two lineages during a short interval of time, obtained here in Lemma \ref{lemm no entrance}, and the uniform convergence of the coalescence time of two lineages, which constitutes our Theorem \ref{theo time coal}. With these estimates, one can obtain the limiting rates of decrease of the number of blocks of $\A^{L,u}$ (namely those of the number of blocks in Kingman's coalescent), and the fact that mergers are only binary as in Cox \& Griffeath~(1986). In particular, the counterpart of their Proposition~2 here gives us that for each $n\in \N$,
\begin{equation}\label{conv holding time}
\lim_{L\rightarrow \infty}\sup_{t\geq 0}\sup_{A_L\in \GA(L,n)^*}\Big|\ \prob_{A_L}\big[|\A^{L,u}(t)|=n\big]-\exp\Big\{-\frac{n(n-1)}{2}\ t\Big\}\Big|=0,
\end{equation}
which we state here because we shall need it for the case $\alpha=1$ (observe that our $L$ corresponds to their $t$). Note that in Proposition 2 of Cox \& Griffeath~(1986), the right-hand side of their equation gives the probability that the number of blocks is less than $n$, instead of equal to $n$ as it is stated. Furthermore, in (\ref{conv holding time}) the supremum is over $t\geq 0$ instead of $t\in [0,T]$ for some $T>0$ (as in Cox \& Griffeath~1986). Our argument for this modification is the fact that the two quantities we are comparing are monotone decreasing in $t$ and both tend to $0$.

Then, the same arguments lead to the proof that any pair of lineages is equally likely to be the first one to coalesce, as in Lemma 1 of Cox \& Griffeath~(1990). The uniformity of the estimates obtained enables us to proceed by induction to show the uniform convergence (on a compact time-interval) of the one-dimensional distributions of $\A^{L,u}$ to those of $\mathcal{K}$, which translate into the uniform convergence of the finite-dimensional distributions, still on intervals of the form $[0,T]$. We refer to \nocite{cox/griffeath:1990} Cox \& Griffeath~(1990) for the complete proof of these results.

\medskip
\noindent \textbf{(ii) Tightness.}

This follows easily from the fact that the labelled partition $\A^L$ with initial value in $\GA(L,n)^*$ for some $n\in \N$ lies in $\GA(L,n)$ immediately after each coalescence event, with probability tending to $1$. Indeed, for each $L\in \N$, let $\gamma_1^L< \ldots < \gamma_{n-1}^L$ be the ranked epochs of jumps of $\A^{L,u}$ (if less than $n-1$ jumps occur, then the last times are equal to $+\infty$ by convention). Let also $n\in \N$, $A_L\in \GA(L,n)^*$ for every $L\geq 1$, and following Ethier \& Kurtz~(1986), for \nocite{ethier/kurtz:1986} every $\delta,T>0$ let $w'(\A^{L,u},T,\delta)$ denote the modulus of continuity of the process $\A^{L,u}$ on the time interval $[0,T]$ and with time-step $\delta$. Let $\e>0$. With the convention that $(+\infty)-(+\infty)=+\infty$, we have
\begin{equation}\label{eq tightness}
\prob_{A_L}\big[w'(\mathcal{A}^L,T,\delta)>\e\big]\leq \sum_{k=2}^n \prob_{A_L}[\gamma_k^L-\gamma_{k-1}^L<\delta].
\end{equation}
An easy recursion using the fact that we consider only finitely lineages and the uniform bounds obtained in Lemma \ref{lemm lineages far} enables us to write that for all $k\in \{1,\ldots,n-1\}$,
$$
\sup_{A_L'\in \GA(L,n)^*}\prob_{A_L'}[\gamma_k^L<\infty\ ;\ \A^L(\vp_L\gamma_k^L)\notin \GA(L,n)]\rightarrow 0,\qquad \mathrm{as\ }L\rightarrow \infty.
$$
This result and an application of the strong Markov property at time $\gamma_{k-1}^L$ yield
\begin{eqnarray}
\prob_{A_L}[\gamma_k^L-\gamma_{k-1}^L<\delta]&=&\E_{A_L}[\ind{\A^L(\vp_L\gamma_{k-1}^L)\in \GA(L,n)}\prob_{\A^L(\vp_L\gamma_{k-1}^L)}[\gamma_1^L<\delta]]+o(1)\nonumber\\
&\leq& \frac{(n-k)(n-k-1)}{2} \sup_{A'_L\in \GA(L,2)^*}\prob_{A'_L}[\gamma_1^L<\delta]+o(1)\label{bound w}
\end{eqnarray}
as $L\rightarrow \infty$, where the last line uses the consistency of the genealogy to bound the probability that a first coalescence event occurs to the sample of lineages before $\delta$ by the sum over all pairs of lineages of this sample of the probability that they have coalesced by time $\delta$ (note that there are at most $(n-k)(n-k-1)/2$ possible pairs just after $\gamma_{k-1}^L$). But these probabilities converge uniformly to $1-e^{-\delta}$ by Theorem \ref{theo time coal}, and so for $\delta$ small enough, we can make the right-hand side of (\ref{bound w}) less than $\e/(n^3)$ for $L$ large enough ($n$ is fixed here). Coming back to (\ref{eq tightness}), this gives us
$$
\limsup_{L\rightarrow \infty} \prob_{A_L}[w'(\mathcal{A}^L,T,\delta)>\e]\leq \e.
$$
Since ${\mathcal P}_n$ is a compact metrisable space, we can apply Corollary 3.7.4 in Ethier \& Kurtz~(1986) to complete the proof.$\hfill\square$

\section{Proof of Theorem~\ref{result alpha=1}}
\label{alpha=1}

We now turn to the case $\psi_L\propto L$. We still have small reproduction events of size $\mathcal{O}(1)$, but now large events have sizes $\mathcal{O}(L)$ (and rate $\mathcal{O}(\rho_L^{-1})$), so that they cover a non-negligible fraction of the torus $\T(L)$. By Lemma \ref{lemma local TCL}, if the lineages were only subject to small reproduction events, the location of a single lineage would be nearly uniformly distributed on $\T(L)$ after a time $t_L\gg L^2$. This suggests several limiting behaviours for the genealogical process $\A^L$, according to how $\rho_L$ scales with $L^2$:
\begin{itemize}
\item If $\rho_L$ is order at most $\mathcal{O}(L^2)$, then large reproduction events occur at times when the locations of the lineages are still correlated with their starting points, and so we expect space (i.e., labels in the representation we adopted) to matter in the evolution of $\A^L$.
\item If $L^2 \ll \rho_L \ll L^2\log L$, then the lineages have the time to homogenise their locations over $\T(L)$ before the first large event occurs, but not to come at distance $2R^s$ from each other. Hence, large events should affect lineages independently of each other, and bring the genealogy down to the common ancestor of the sample before any pair of lineages experiences a coalescence due to small events.
\item If $\rho_L\approx L^2 \log L$, the fact that pairs of lineages have now the time to gather at distance $2R^s$ should add a Kingman part (i.e., almost surely binary mergers) to the genealogical process obtained in the previous point.
\item If $\rho_L\gg L^2\log L$, Kingman's coalescent due to small reproduction events should bring the ancestry of a sample of lineages down to a single lineage before any large event occurs, so that the limiting genealogy will not see these large events.
\end{itemize}
\emph{Proof of Theorem~\ref{result alpha=1}: }For $(a)$, let us write down the generator $\ov{\G}_L$ of $\bar{A}^L$ applied to functions of the $\T(1)$-labelled partitions of $\{1,\ldots,n\}$. Recall the notation $x_a$ for the label of the block $a$ of a labelled partition $A\in {\mathcal P}_n^{\ell}$ (introduced in Notation \ref{notation for partitions}), and write $|A|$ for the number of blocks of $A$. For each $L\geq 1$, $f$ of class $C^3$ with respect to the labels and $A\in {\mathcal P}_n^{\ell}$ such that any pair $(a_1,a_2)$ of blocks of $A$ satisfies $|x_{a_1}-x_{a_2}|\geq (2R^s)/L$, we have
\begin{eqnarray}
\ov{\G}_Lf(A)& = & \rho_L\sum_{i=1}^{|A|} \int_{\T(L)} dy \int_0^{R^s}\mu^s(dr)\frac{L_r(y)}{\pi r^2}\int_0^1\nu^s_r(du)u \nonumber \\
& & \qquad \quad \times \Big[f\Big(A\setminus \big\{(a_i,x_{a_i})\big\}\cup \Big\{\Big(a_i,x_{a_i}+\frac{y}{L}\Big)\Big\}\Big)-f(A)\Big] + \G^{(B)}(A),\phantom{AAA} \label{generator a}
\end{eqnarray}
where we wrote $A=\big\{(a_1,x_{a_1}),\ldots,(a_{|A|},x_{a_{|A|}})\big\}$ and \setlength\arraycolsep{1pt}
\begin{eqnarray*}
&\G^{(B)}&(A)\\
& =& \frac{1}{c^2}\int_{\T(1)} dz\int_0^{(\sqrt{2})^{-1}}\mu^B(dr)\int_{B(z,cr)}\frac{dy}{V_{cr}}\sum_{I\subset \{1,\dots,|A|\}}\prod_{i\in I}\ind{x_i\in B(z,cr)}\prod_{j\notin I}\ind{x_j\notin B(z,cr)}\\
& \ \  \times& \sum_{J\subset I}\int_0^1u^{|J|}(1-u)^{|I|-|J|}\nu^B_r(du)\bigg[f\bigg(A\setminus\Big(\bigcup_{i\in J}\{(a_i,x_{a_i})\}\Big)\cup\Big\{\Big(\bigcup_{i\in J}a_i,y\Big)\Big\}\bigg)-f(A)\bigg]
\end{eqnarray*}
is the generator of the coalescence events due to large reproduction events (recall $V_r$ is the volume of the ball $B_{\IT(1)}(0,r)$). Note that $\G^{(B)}$ does not depend on $L$. Let us look at a particular term in the sum on the right-hand side of (\ref{generator a}). Since $f$ is of class $C^3$ with respect to the labels of the blocks, a Taylor expansion and the symmetry of the jumps due to small events give us
\begin{eqnarray*}
\rho_L\int & dy & \int_0^{R^s}\mu^s(dr)\frac{L_r(y)}{\pi r^2}\int_0^1\nu_r(du)u \Big[f\Big(A\setminus \big\{(a_i,x_{a_i})\big\}\cup \Big\{\Big(a_i,x_{a_i}+\frac{y}{L}\Big)\Big\}\Big)-f(A)\Big] \\
& = & \frac{\rho_L}{L^2}\ \frac{\sigma^2_s}{2}\ \Delta_i f(A) +\mathcal{O}\Big(\frac{\rho_L}{L^3}\Big),
\end{eqnarray*}
where $\Delta_i$ is the Laplacian operator on $\T(1)$ applied to the label of the block $a_i$ only. Since $\rho_LL^{-2}\rightarrow b\in [0,\infty)$ by assumption and because $f$ is continuous on a compact space, we obtain that $\ov{\G}_Lf$ defined on the compact set $E_L\equiv \big\{A\in {\mathcal P}_n^{\ell}: L|x_{a_i}-x_{a_j}|\geq 2R^s\ \forall i\neq j\}$ converges uniformly towards
$$
\ov{\G}f(A)\equiv \frac{b\sigma_s^2}{2} \sum_{i=1}^{|A|}\Delta_i f(A) + \G^{(B)}f(A).
$$
Now, by the same technique as in Section \ref{levy processes}, one can prove that the gathering time at distance $2R^s$ of two lineages starting at distance $\mathcal{O}(L)$ on $\T(L)$ and subject only to small events converges uniformly on the time scale $\frac{L^2\log L}{\pi \sigma_s^2}$ to an $\mathrm{Exp}(1)$ random variable (in the sense of Lemma \ref{lemma entrance_levy}). In addition, since the new location of a lineage affected by a large event is chosen uniformly over a ball of $\T(L)$ whose radius is of order $\mathcal{O}(L)$, if a large event affects a pair of lineages but does not lead to their coalescence, then the probability that the lineages are at distance less than $L(\log L)^{-1}$ just after the event vanishes as $L\rightarrow 0$. If we call $\check{T}^*_L$ the first time at which two lineages on $\T(L)$ are gathered at distance $2R^s$ and $t_L^*$ their coalescence time in the original timescale, we readily obtain that for any $u>0$, and $x'_1\neq x'_2\in \T(1)^2$,
$$
\lim_{L\rightarrow \infty}\prob_{\wp_2(Lx'_1,Lx'_2)}\big[t_L^*>\check{T}^*_L\ ;\ \check{T}^*_L\leq \rho_Lu\big]= 0.
$$
Indeed, as we already mentioned, if a large event does not make the lineages coalesce then with probability tending to one, the latter start at separation $\mathcal{O}(L)$ and do not have the time to meet at distance $2R^s$ before the next large event. Now, the number of large reproduction events that the pair of lineages experiences before time $\rho_Lu$ can be stochastically bounded by a Poisson random variable whose parameter is finite and independent of $L$. Hence, if none of them leads to a coalescence then with probability tending to $1$, $\check{T}^*_L> \rho_Lu$. It follows that, if $u>0$ is fixed, we can use the consistency of the genealogy and write
$$
\prob_{\wp_n(L\mathbf{x})}[\exists t\in [0,u]: \bar{\A}^L(t)\notin E_L]\leq \sum_{i< j=1}^{n}\prob_{\{(\{i\},Lx_i),(\{j\},Lx_j)\}}\big[t_L^*>\check{T}^*_L ; \check{T}^*_L\leq \rho_Lu\big] \rightarrow 0.
$$
Consequently, one can use Corollary 4.8.7 in Ethier \& Kurtz~(1986) (with $E_L$ as the subspace of interest in condition $(f)$) to conclude that the law under $\prob_{\wp_n(L\mathbf{x})}$ of $\bar{A}^L$ converges to that of $\bar{\A}^{\infty,b,c}$ as processes in the Skorohod space of all c\`adl\`ag paths with values in the $\T(1)$-labelled partitions of $\{1,\ldots,n\}$.

Let us now prove $(b)$. Recall the assumption that the total rate at which large events occur is finite, that is $M\equiv c^{-2}\mu^B([0,(\sqrt{2})^{-1}])<\infty$. Let us first analyse what happens during the first event which may affect the unlabelled ancestral partition.

Define for each $L\geq 1$ the stopping time $e^L_1$ by the following property: $\rho_Le_1^L$ is the first time on the original timescale at which either a large event occurs, or $\A^L$ undergoes a coalescence event due to small reproduction events. Since large and small reproduction events are independent, $\rho_Le_1^L$ has the same distribution as the minimum of two following independent random times:
\begin{itemize}
\item the first time of occurrence of a large event, that is an $\mathrm{Exp}\big(M/\rho_L\big)$-random variable.
\item the time $t^*_L$ at which a first coalescence event occurs between lineages of the genealogical process $\tilde{\A}^L$ evolving only owing to small reproduction events.
\end{itemize}
By (\ref{conv holding time}) applied to the case $\rho_L\equiv +\infty$ (i.e., no large events occur), $\frac{2\pi \sigma^2_s}{L^2\log L}\ t^*_L$ converges to an $\mathrm{Exp}\big(n(n-1)/2\big)$-random variable under $\prob_{A_L}$, uniformly in $(A_L)_{L\in \N}$ such that $A_L\in \GA(L,n)^*$ for every $L$. It is then straightforward to obtain
\begin{equation}\label{conv coal times}
\lim_{L\rightarrow \infty}\sup_{t\geq 0}\sup_{A_L\in \GA(L,n)^*}\bigg|\ \prob_{A_L}[e_1^L>t] - \exp\Big(-\Big\{M+\beta \frac{n(n-1)}{2}\Big\}t\Big)\bigg|=0,
\end{equation}
where the formulation is also valid for $\beta =0$. Also, by the independence of $\Pi_L^s$ and $\Pi_L^B$, for every $(A_L)_{L\in \N}$ as above we have (with an abuse of notation)
$$
\prob_{A_L}\big[\rho_Le^L_1 = t_L^*\big]= \E_{A_L}\Big[\exp\Big\{-\frac{M}{\rho_L}\ t_L^*\Big\}\Big].
$$
Using Fubini's theorem and a change of variable, we can write
\begin{eqnarray*}
\E_{A_L}\Big[\exp\Big\{-\frac{M}{\rho_L}\ t_L^*\Big\}\Big] &=& \int_0^1 \prob_{A_L}\Big[\exp\Big\{-\frac{M}{\rho_L}\ t_L^*\Big\}>s\Big]ds\\
 &=& \int_0^1 \prob_{A_L}\bigg[\frac{2\pi \sigma_s^2}{L^2\log L}\ t_L^*< -\frac{2\pi\sigma^2_s\rho_L}{ML^2\log L}\ \log s\bigg]ds\\
&=& \frac{ML^2\log L}{2\pi \sigma_s^2\rho_L}\int_0^{\infty}e^{-\frac{ML^2\log L}{2\pi \sigma_s^2\rho_L} u}\ \prob_{A_L}\bigg[\frac{2\pi \sigma_s^2}{L^2\log L}\ t_L^*<u\bigg]du\\
&=& 1- \frac{ML^2\log L}{2\pi \sigma_s^2\rho_L}\int_0^{\infty}e^{-\frac{ML^2\log L}{2\pi \sigma_s^2\rho_L} u}\ \prob_{A_L}\bigg[\frac{2\pi \sigma_s^2}{L^2\log L}\ t_L^*\geq u\bigg]du.
\end{eqnarray*}
When $\beta>0$, we have $\frac{ML^2\log L}{2\pi \sigma_s^2\rho_L}\rightarrow \frac{M}{\beta}$ and so we can use the uniform convergence derived in (\ref{conv holding time}) and the fact that the distribution of $t_L^*$ does not charge points to conclude that
$$
\lim_{L\rightarrow \infty}\sup_{A_L\in \GA(L,n)^*}\bigg|\prob_{A_L}\big[\rho_Le_1^L = t_L^*\big]-\frac{\beta n(n-1)}{\beta n(n-1)+2M}\bigg|=0.
$$
The limit holds also for $\beta=0$ by a trivial argument. A byproduct of this result is the existence of a constant $C_0>0$ and $L_0\in \N$ such that, for all $L\geq L_0$ and $(A_L)_{L\in \N}$ as above, $\prob_{A_L}[\rho_Le_1^L < t_L^*]\geq C_0$. We shall need this fact in the next paragraph.

By Theorem \ref{result alpha<1} in the case $\rho_L\equiv\infty$, up to an error term tending uniformly to $0$, on the event $\big\{\rho_Le_1^L=t_L^*\big\}$ the transition occurring to $\A^{L,u}$ at time $\rho_Le_1^L$ is the coalescence of a pair of blocks, each pair having the same probability to be the one which coalesces. Let us show that, conditioned on $\big\{\rho_Le_1^L < t^*_L\big\}$, the locations of the lineages at time $(\rho_Le_1^L)-$ are approximately distributed as $n$ independent uniform random variables on $\T(L)$. We use again the notation $\tau_{ij},\tau_{ij}^*$ and $\tau,\tau^*$($=t_L^*$ here) introduced in the proof of Theorem \ref{result alpha<1} for the gathering time at distance $2R^s$ and the coalescence time of lineages $i$ and $j$, and their minima (once again on the original timescale). These quantities depend on $L$ but, for the sake of clarity, we do not reflect that in our notation. In order to use our results on L\'evy processes, we need to make sure that no pairs of lineages have come at  distance less than $2R^s$ before time $\rho_Le_1^L$. We have for each $L\in \N$
\begin{equation}\label{eq big event first}
\prob_{A_L}\big[\tau <\rho_Le_1^L\big|\ \rho_Le_1^L< t_L^*\big]\leq \sum_{i<j=1}^n \prob_{A_L}\big[\tau_{ij}<\rho_Le_1^L\big|\ \rho_Le_1^L< t_L^*\big],
\end{equation}
Each term $(i,j)$ on the right-hand side of (\ref{eq big event first}) is bounded by \setlength\arraycolsep{1pt}
\begin{eqnarray}
\prob_{A_L}\big[&\tau_{ij}&<\rho_Le_1^L-\log L\big|\ \rho_Le_1^L< t_L^*\big] \nonumber\\ & &\qquad +\prob_{A_L}\big[\rho_Le_1^L- \log L \leq \tau_{ij}<\rho_Le_1^L\big|\ \rho_Le_1^L< t_L^*\big]\nonumber \\
&\leq & C_0^{-1}\Big\{\prob_{A_L}\big[\tilde{\tau}_{ij}^* >\tilde{\tau}_{ij}+ \log L\big] + \prob_{A_L}\big[\tilde{\tau}_{ij}\in [\varsigma_L-\log L,\varsigma_L)\big]\Big\}, \label{independence}
\end{eqnarray}
where for each $L\in \N$, $\varsigma_L$ is an $\mathrm{Exp}(M/\rho_L)$-random variable independent of all other variables, and $\tilde{\tau}_{ij}$ and $\tilde{\tau}_{ij}^*$ are defined as above, but for the process $\tilde{\A}^L$. By the strong Markov property applied at time $\tilde{\tau}_{ij}$ and the result of Proposition \ref{prop coal time} $(b)$, the first term on the right-hand side of (\ref{independence}) converges to $0$ uniformly in $A_L\in \GA(L,n)^*$. By a simple change of variable, the second term is equal to
$$
M\int_0^{\infty} e^{-Ms}\ \prob_{A_L} \big[\tilde{\tau}_{ij}\in [\rho_Ls-\log L,\rho_Ls)\big]ds \leq M\int_0^{\infty} e^{-Ms}\ C\ \frac{\log L}{L^2}\ ds\ \rightarrow\ 0,
$$
where the inequality comes from Lemma \ref{lemm no entrance}. Therefore, back to (\ref{eq big event first}) we obtain that
\begin{equation}\label{independence motions tilde A}
\lim_{L\rightarrow \infty}\sup_{A_L\in \GA(L,n)^*}\prob_{A_L}\big[\tau<\rho_Le_1^L\big|\ \rho_Le_1^L< t_L^*\big]= 0.
\end{equation}
Now, let $D_1,\ldots,D_n$ be $n$ measurable subsets of $\T(1)$, and for each $i\in \{1,\ldots,n\}$ and $L\geq 1$, let $LD_i\subset \T(L)$ be the dilation of $D_i$ by a factor $L$. Let us show that
\begin{eqnarray}
\lim_{L\rightarrow \infty}\sup_{A_L\in \GA(L,n)^*}\Big|\ \prob_{A_L}&\Big[&(\xi_1^L,\ldots,\xi_n^L)(\rho_Le_1^L-)\in (LD_1)\times \ldots \times (LD_n)\big|\rho_Le_1^L< t^*_L\Big]\nonumber \\& & \qquad \qquad \qquad \qquad\qquad \qquad \qquad  - \prod_{i=1}^n \mathrm{Leb}(D_i)\Big|=0,\label{eq homogen}
\end{eqnarray}
where $\xi_i^L(t)$ denotes the location of the $i$-th lineage of $\A^L$ at time $t$. To do so, let us use the fact that on the event $\big\{\rho_Le_1^L< t^*_L\big\}$, the genealogical process $\A^L$ up to time $\rho_Le_1^L$ has the same distribution as $\tilde{\A}^L$ up to time $\varsigma_L$ and on the event $\{\tilde{\tau}^*>\varsigma_L\}$. We have \setlength\arraycolsep{1pt}
\begin{eqnarray}
&\prob_{A_L}&\Big[(\xi_1^L,\ldots,\xi_n^L)(\rho_Le_1^L-)\in \prod_{i=1}^n(LD_i)\Big|\rho_Le_1^L< t^*_L\Big]\nonumber\\
&=&\frac{1}{\prob_{A_L}[\rho_Le_1^L< t_L^*]}\ \prob_{A_L}\Big[(\xi_1^L,\ldots,\xi_n^L)(\rho_Le_1^L-)\in \prod_{i=1}^n(LD_i);\ \rho_Le_1^L< t^*_L\Big]\nonumber\\
&=& \frac{1}{\prob_{A_L}[\rho_Le_1^L< t_L^*]}\ \prob_{A_L}\Big[(\tilde{\xi}_1^L,\ldots,\tilde{\xi}_n^L)(\varsigma_L-)\in \prod_{i=1}^n(LD_i);\ \varsigma_L< \tilde{\tau}^*\Big]\nonumber \\
&=& \frac{1}{\prob_{A_L}[\rho_Le_1^L< t_L^*]}\ \prob_{A_L}\Big[(\tilde{\xi}_1^L,\ldots,\tilde{\xi}_n^L)(\varsigma_L-)\in \prod_{i=1}^n(LD_i);\ \varsigma_L< \tilde{\tau}\Big]+\eta_L(A_L)\nonumber \\
&=&\frac{M}{\prob_{A_L}[\rho_Le_1^L< t_L^*]}\ \int_0^\infty ds\ e^{-Ms} \prob_{A_L}\Big[(\tilde{\xi}_1^L,\ldots,\tilde{\xi}_n^L)(\rho_Ls-)\in \prod_{i=1}^n(LD_i);\nonumber \\
& & \qquad \qquad \qquad \qquad\qquad \qquad \qquad \qquad\qquad  \tilde{\tau}>\rho_Ls\Big]+\eta_L(A_L), \label{au1}
\end{eqnarray}
where $\eta_L(A_L)$ tends to $0$ uniformly in $(A_L)_{L\in \N}$ by (\ref{independence motions tilde A}) and the fact that $\prob_{A_L}\big[\rho_Le_1^L< t_L^*\big]$ does not vanish.

Let us fix $s>0$ for a moment, and consider the corresponding probability within the integral. Up to time $\tilde{\tau}$, the movements of the lineages are distributed as $n$ independent copies $\hat{\xi}^L_1,\ldots,\hat{\xi}^L_n$ of the motion of a single lineage, for which an easy modification of Lemma \ref{lemma local TCL} $(b)$ tells us that,  if $(\e_L)_{L\in \N}$ is such that $\e_L\rightarrow 0$ but $\e_L\rho_L\gg L^2$ as $L\rightarrow \infty$,
\begin{equation}\label{asymptotic uniformity}
\lim_{L\rightarrow \infty}\ \sup_{v\geq \e_L}\ \sup_{x\in \T(L)}\ \left|\prob_{x}\big[\hat{\xi}^L(v\rho_L)\in (LD)\big]-\mathrm{Leb}(D)\right|=0.
\end{equation}
However, it is not entirely clear that this convergence will still hold for $n$ independent lineages on the event $\{\hat{\tau}>\rho_Ls\}$ (where $\hat{\tau}$ is the first time at which at least two of them come at distance less than $2R^s$). Keeping the notation $A_L$ for the initial value of the set of lineages and denoting the set of $n$ (non-coalescing) motions by $\hat{\A}^L$, we have
\begin{eqnarray*}
&\prob_{A_L}&\big[(\hat{\xi}_1^L,\ldots,\hat{\xi}_n^L)(\rho_L s-)\in (LD_1)\times \ldots \times (LD_n);\ \hat{\tau}\leq \rho_Ls\big]\\
&=& \E_{A_L}\bigg[\ind{\hat{\tau}\leq \rho_Ls}\ \prob_{\hat{\A}^L(\hat{\tau})}\Big[(\hat{\xi}_1^L,\ldots,\hat{\xi}_n^L)\big((\rho_L s-\hat{\tau})-\big)\in (LD_1)\times \ldots \times (LD_n)\Big]\bigg].
\end{eqnarray*}
Splitting the preceding integral into $\big\{\rho_L(s-\e_L)\leq \hat{\tau}\leq \rho_Ls\big\}$ and $\big\{\hat{\tau}<\rho_L(s- \e_L)\big\}$, we can use (\ref{asymptotic uniformity}) in the latter case to write \setlength\arraycolsep{1pt}
\begin{eqnarray}
\E_{A_L}&\bigg[&\ind{\hat{\tau}\leq \rho_Ls}\ \prob_{\hat{\A}^L(\hat{\tau})}\Big[(\hat{\xi}_1^L,\ldots,\hat{\xi}_n^L)\big((\rho_L s-\hat{\tau})-\big)\in \prod_{i=1}^n(LD_i)\Big]\bigg]\nonumber\\
& & =\E_{A_L}\bigg[\ind{\rho_L(s-\e_L)\leq \hat{\tau}\leq \rho_Ls}\ \prob_{\hat{\A}^L(\hat{\tau})}\Big[(\hat{\xi}_1^L,\ldots,\hat{\xi}_n^L)\big((\rho_L s-\hat{\tau})-\big)\in \prod_{i=1}^n(LD_i)\Big]\bigg]\nonumber\\
& &  \ \ + \Big(\prod_{i=1}^n \mathrm{Leb}(D_i)\Big)\prob_{A_L}[\hat{\tau}<\rho_L(s- \e_L)]+ \delta_L(A_L), \label{au2}
\end{eqnarray}
where $(\delta_L(A_L))_{L\in \N}$ tends to zero uniformly in $(A_L)_{L\in \N}$ as $L$ tends to infinity (we still impose that $A_L\in \GA(L,n)^*$ for every $L$). By the convergence of the distribution function of $\frac{\tilde{\tau}}{L^2\log L}$ to that of an exponential random variable, uniformly in the time variable and in $(A_L)_{L\in \N}$, we obtain that $\prob_{A_L}[\rho_L(s-\e_L)\leq \hat{\tau}\leq \rho_Ls]$ converges to $0$ uniformly in $(A_L)_{L\in \N}$ (which is also true if $\beta=0$, i.e., $\rho_L\ll L^2\log L$). Hence, we can find a sequence $(\delta'_L(A_L))_{L\in\N}$ decreasing to $0$ uniformly in $(A_L)_{L\in \N}$, such that the whole sum on the right-hand side of (\ref{au2}) is equal to
$$
\Big(\prod_{i=1}^n \mathrm{Leb}(D_i)\Big)\prob_{A_L}[\hat{\tau}\leq \rho_Ls ]+ \delta'_L(A_L).
$$
Likewise, we can find another sequence $(\delta_L'')_{L\in \N}$ decreasing to zero uniformly in $(A_L)_{L\in \N}$ such that
$$
\prob_{A_L}\big[(\hat{\xi}_1^L,\ldots,\hat{\xi}_n^L)(\rho_Ls-)\in (LD_1)\times \ldots \times (LD_n)\big]= \prod_{i=1}^n \mathrm{Leb}(D_i) + \delta''_L(A_L).
$$
Subtracting the two last equalities, we obtain
$$
\prob_{A_L}\Big[(\hat{\xi}_1^L,\ldots,\hat{\xi}_n^L)(\rho_L s-)\in \prod_{i=1}^n(LD_i); \hat{\tau}> \rho_Ls\Big]= \bigg\{\prod_{i=1}^n \mathrm{Leb}(D_i)\bigg\}\prob_{A_L}[\hat{\tau}> \rho_Ls ]+ o(1),
$$
where the remainder decreases to $0$ uniformly in $s>0$ and $(A_L)_{L\geq 1}$ such that $A_L\in \GA(L,n)^*$ for each $L$. Coming back to (\ref{au1}), we obtain that it is equal to \setlength\arraycolsep{1pt}
\begin{eqnarray*}
\frac{M}{\prob_{A_L}[\rho_Le_1^L< t_L^*]}\int_0^\infty& ds& e^{-Ms} \bigg\{\Big(\prod_{i=1}^n \mathrm{Leb}(D_i)\Big)\prob_{A_L}[\tilde{\tau}> \rho_Ls ] +o(1)\bigg\}\\
&= & \frac{\prob_{A_L} [\tilde{\tau}> \varsigma_L]}{\prob_{A_L}[\tilde{\tau}^*>\varsigma_L]}\ \prod_{i=1}^n \mathrm{Leb}(D_i) + o(1) \\
&=& \frac{\prob_{A_L} [\tilde{\tau}^* > \varsigma_L]+o(1)}{\prob_{A_L}[\tilde{\tau}^*>\varsigma_L]}\ \prod_{i=1}^n \mathrm{Leb}(D_i)+o(1),
\end{eqnarray*}
where the last line uses (\ref{independence motions tilde A}). We can thus conclude that (\ref{eq homogen}) holds.

Condition on the first event being a large reproduction event. By the description of such an event, the result for the genealogical process is the merger of at most one group of blocks into a bigger block. Furthermore, the transitions depend only on the number of blocks and their labels, so for convenience we derive the transition probabilities for $A_L$ of the form $\wp_n(\mathbf{x})$ only, although we shall use the result later for more general labelled partitions. Let $\pi$ be a partition of $\{1,\ldots,n\}$ such that $\pi$ has exactly one block of size greater than $1$, which we call $J$. Then if the large event has centre $x$ and radius $cr$ in $\T(1)$, the probability that the transition undergone by $\A^{L,u}$ is $\wp_n\rightarrow \pi$ is the probability that at this time, at least all the lineages in $J$ have labels in $B(x,cr)$ and are really affected by the event, and all the other lineages present in $B(x,cr)$ are not affected by the event. Summing over all possible choices $I\subset \{1,\ldots,n\}\setminus J$ for these ``other lineages'' ($I$ can be empty) and using (\ref{eq homogen}), the probability of the transition $\wp_n\rightarrow \pi$ up to a vanishing error is given by \begin{eqnarray}
\sum_{I}&V_{cr}^{|J|+|I|}&(1-V_{cr})^{n-|J|-|I|}\int_0^1u^{|J|}(1-u)^{|I|}\nu_r^B(du)\nonumber\\
&=& \int_0^1 (uV_{cr})^{|J|}\sum_{i=0}^{n-|J|}\binom{n-|J|}{i}V_{cr}^i(1-V_{cr})^{n-|J|-i}(1-u)^i \nu_r^B(du) \nonumber \\
&=& \int_0^1 (uV_{cr})^{|J|}((1-u)V_{cr}+1-V_{cr})^{n-|J|}\nu_r^B(du)\nonumber \\
&=& \int_0^1 (uV_{cr})^{|J|}(1-uV_{cr})^{n-|J|}\nu_r^B(du).\label{transition lambda}
\end{eqnarray}

We now have the results we need to show $(b)$. For every $L\in \N$, let us consider again the time $e_1^L$ introduced earlier, and define for each integer $i\geq 2$, \begin{eqnarray*}
e_i^L= \inf\big\{t> e_{i-1}^L\ &:&\ \rho_Lt\in \Pi_L^B\mathrm{\ or\ } \rho_Lt \mathrm{\ is\ the\ epoch\ of\ a\ coalescence} \\
& & \mathrm{due\ to\ small\ events}\big\}.
\end{eqnarray*}
Let us also define similar times corresponding to $\Lambda^{(\beta,c)}$. From the expression of its rates given in Definition \ref{def lambda coal}, $\Lambda^{(\beta,c)}$ is composed of a Kingman part (i.e., only binary mergers) run at rate $\beta$, and of a set of multiple mergers due to the part $\Lambda^{(0)}$ of its $\Lambda$-measure with the atom at $0$ removed. Furthermore, the finite measure $\Lambda^{(0)}$ on $[0,1]$ is given by
\begin{eqnarray*}
\Lambda^{(0)}(dv)&= &c^{-2}v^2\int_0^{(\sqrt{2})^{-1}}\nu^B_r\big(\big\{u:\ uV_{cr}\in dv\big\}\big)\mu^B(dr) \\
& = & c^{-2}v^2\int_0^{(\sqrt{2})^{-1}}\ind{V_{cr}\geq v}\nu_r^B\Big(d\frac{v}{V_{cr}}\Big)\mu^B(dr) .
\end{eqnarray*}
Following Pitman's Poissonian construction of a coalescent with multiple mergers (whose $\Lambda$-measure has no atom at $0$, see Pitman~1999), \nocite{pitman:1999} let us define $\Pi$ as a Poisson point process on $\R_+\times [0,1]$ with intensity $dt\otimes v^{-2}\Lambda^{(0)}(dv)$. Note that because of our assumption on $M$, $v^{-2}\Lambda^{(0)}(dv)$ is also a finite measure, with total mass $M$. The atoms of $\Pi$ constitute the times at which $\Lambda^{(\beta,c)}$ acting on the partitions of $\N$ experiences a multiple collision, and the probabilities that any given lineage is affected by the event. The Kingman part of $\Lambda^{(\beta,c)}$ is superimposed on this construction by assigning to all pairs of blocks of the current partition independent exponential clocks with parameter $\beta$, giving the time at which the corresponding pair merges into one block.

From now on, we consider only the restriction of $\Lambda^{(\beta,c)}$ to ${\mathcal P}_n$, although we do not make it appear in the notation. Let $e_1$ be the minimum of the first time a pair of blocks of $\Lambda^{(\beta,c)}$ merges due to the Kingman part and of the time corresponding to the first point of $\Pi$. Define $e_i$ in a similar manner for all $i\geq 2$, so that $(e_i)_{i\in \N}$ is an increasing sequence of random times at which $\Lambda^{(\beta,c)}$ may undergo a transition. Our goal is to show that the finite-dimensional distributions of $\big\{(e_i^L,\A^{L,u}(e_i^L)),\ i\in \N\big\}$ under $\prob_{A_L}$ converge to those of $\big\{(e_i,\Lambda^{(\beta,c)}(e_i)),\ i\in \N\big\}$ under $\prob_{\wp_n}$, as $L\rightarrow \infty$. Since $\A^{L,u}$ (resp., $\Lambda^{(\beta,c)}$) can jump only at the times $e_i^L$ (resp., $e_i$), the fact that only finitely many jumps occur to $\Lambda^{(\beta,c)}$ in any compact time interval, together with Proposition 3.6.5 in Ethier \& Kurtz~(1986) enable us to conclude that this convergence yields $(b)$. We proceed by induction, by showing that for each $i\in \N$:

\medskip
\noindent \textit{$H(i)$ : if $a_L\in \GA(L,n)$ for each $L$ and there exists $\pi_0\in {\mathcal P}_n$ such that for all $L\in \N$, $\mathrm{bl}(a_L)=\pi_0$, then as $L\rightarrow \infty$}
$$
\mathcal{L}_{\prob_{a_L}}\big(\big\{(e_1^L,\A^{L,u}(e_1^L)),\ldots,(e_i^L,\A^{L,u}(e_i^L))\big\}\big) \Rightarrow \mathcal{L}_{\prob_{\pi_0}}\big(\big\{(e_1,\Lambda^{(\beta,c)}(e_1)),\ldots,(e_i,\Lambda^{(\beta,c)}(e_i))\big\}\big).
$$
(Note that $a_L$ can have less than $n$ blocks).

Let us start by $H(1)$. Let $t\geq 0$, $\pi\in {\mathcal P}_n$ and write $n_0$ for the number of blocks of $\pi_0$. We have, in the notation used in the previous paragraph (and with $\tilde{\A}^{L,u}$ defined as the unlabelled partition induced by $\tilde{\A}^L$ on the timescale $\rho_L$), \setlength\arraycolsep{1pt}
\begin{eqnarray}
&\prob_{a_L}&\big[ e_1^L\leq  t;\ \A^{L,u}(e_1^L)=\pi \big]\nonumber\\
&= & \prob_{a_L}\big[e_1^L\leq t;\ \A^{L,u}(e_1^L)=\pi;\ \rho_Le_1^L=t_L^* \big]+ \prob_{a_L}\big[e_1^L\leq t;\ \A^{L,u}(e_1^L)=\pi;\ \rho_Le_1^L< t_L^* \big]\nonumber\\
&=& \prob_{a_L}\big[t_L^*\leq \rho_Lt;\ \tilde{\A}^{L,u}(t_L^*/\rho_L)=\pi;\ t_L^* <\zeta_L \big]\label{H(1) term one}\\
& & + \prob_{a_L}\big[e_1^L\leq t;\ \A^{L,u}(e_1^L)=\pi\big|\ \rho_Le_1^L<t_L^* \big]\prob_{a_L}\big[\rho_Le_1^L <t_L^*\big].\label{H(1) term two}
\end{eqnarray}
By Theorem \ref{result alpha<1} applied with $\rho_L\equiv +\infty$, $\tilde{\A}^{L,u}$ with initial value $a_L$ converges as $L\rightarrow \infty$ to Kingman's coalescent $\mathcal{K}^{(\beta)}$ started at $\pi_0$ and run at rate $\beta$, as a process in $D_{{\mathcal P}_n}[0,\infty)$ (if $\beta=0$, then $\tilde{\A}^{L,u}$ converges to the constant process equal to $\pi_0$). Hence, by the independence of $\tilde{\A}^L$ and $\zeta_L$ for every $L$ and a simple time-change, the quantity in (\ref{H(1) term one}) tends to that corresponding to $\mathcal{K}^{(\beta)}$, that is
\begin{equation}\label{first event kingman}
 \prob_{\pi_0}\big[\mathcal{K}^{(\beta)}(e_1^{\mathcal{K}})=\pi\big]  \prob_{\pi_0}\big[e_1^{\mathcal{K}}< t\wedge \zeta \big],
\end{equation}
where $e_1^{\mathcal{K}}$ is distributed like an $\mathrm{Exp}\big(\beta\frac{n_0(n_0-1)}{2}\big)$-random variable and stands for the epoch of the first event occurring to $\mathcal{K}^{(\beta)}$, and $\zeta$ is an $\mathrm{Exp}(M)$-random variable. By the construction of $\Lambda^{(\beta,c)}$ given in the last paragraph, (\ref{first event kingman}) is the probability that the first event occurring to $\Lambda^{(\beta,c)}$ happens before time $t$, is due to the Kingman part of the coalescent and leads to the transition $\pi_0\rightarrow \pi$. For (\ref{H(1) term two}), note first that because $\Pi_L^B$ and $\Pi_L^s$ are independent, if we condition on $\rho_Le_1^L$ being the time of the first point $(t_1^L,x_1^L,r_1^L)$ of $\Pi_L^B$, then $e_1^L$ and the pair $(x_1^L,r_1^L)$ are independent. Hence, we have for each $L\geq 1$
\begin{eqnarray*}
\prob_{a_L}\big[e_1^L\leq t&;& \A^{L,u}(e_1^L)=\pi\big|\ \rho_Le_1^L<t_L^* \big]\\
& =& \prob_{a_L}\big[e_1^L\leq t\big|\ \rho_Le_1^L<t_L^* \big]\prob_{a_L}\big[\A^{L,u}(e_1^L)=\pi\big|\ \rho_Le_1^L<t_L^* \big].
\end{eqnarray*}
Using (\ref{conv coal times}) and the same reasoning as for (\ref{H(1) term one}), we can write \setlength\arraycolsep{1pt}
\begin{eqnarray*}
\prob_{a_L} \big[e_1^L\leq t\big|\ \rho_Le_1^L<t_L^* \big]&\prob_{a_L}&\big[\rho_Le_1^L<t_L^*\big] \\
&=& \prob_{a_L}\big[e_1^L\leq t;\ \rho_Le_1^L<t_L^*\big]\\
&=& \prob_{a_L}\big[e_1^L\leq t\big]-\prob_{a_L}\big[e_1^L\leq t;\ \rho_L e_1^L=t_L^*\big]\\
&\rightarrow & \exp\Big\{-\Big(M+\beta \frac{n_0(n_0-1)}{2}\Big)t \Big\}-\prob_{\pi_0}\big[e_1^{\mathcal{K}} \leq t \wedge \zeta \big]\\
&=& \prob_{\pi_0}\big[\zeta < t \wedge e_1^{\mathcal{K}}\big],
\end{eqnarray*}
where the last equality comes from the fact that an $\mathrm{Exp}\big(\beta\frac{n_0(n_0-1)}{2}\ +M\big)$-random variable has the same distribution as the minimum of an $\mathrm{Exp}\big(\beta\frac{n_0(n_0-1)}{2}\big)$- and an $\mathrm{Exp}(M)$-random variables, independent of each other. In addition, by the calculation done in (\ref{transition lambda}),
$$
\prob_{a_L}\big[\A^{L,u}(e_1^L)=\pi\big|\ \rho_L e_1^L<t_L^*\big]\rightarrow \prob_{\pi_0}\big[\Lambda^{(0)}(e_1^{\Lambda})=\pi\big],\qquad \mathrm{as}\ L\rightarrow \infty,
$$
where $e_1^{\Lambda}$ is the time of the first event of $\Pi$. Combining the above, and recognizing the transition probability of $\Lambda^{(\beta,c)}$ through the decomposition obtained, we can write
$$
\lim_{L\rightarrow \infty}\prob_{a_L}\big[e_1^L\leq t;\ \A^{L,u}(e_1^L)=\pi \big]=\prob_{\pi_0}\big[e_1\leq t;\ \Lambda^{(\beta,c)}(e_1)=\pi \big].
$$
Since this result holds for each $t\geq 0$ and $\pi_0\in {\mathcal P}_n$, using a monotone class argument we can conclude that the distribution of $\big(e_1^L,\A^{L,u}(e_1^L)\big)$ under $\prob_{a_L}$ converges to the distribution of $(e_1,\Lambda^{(\beta,c)}(e_1))$ under $\prob_{\pi_0}$ as $L\rightarrow \infty$. This proves $H(1)$.

Suppose that $H(i-1)$ holds for some $i\geq 2$. Let $D\subset (\R_+)^{i-1}$, $t\geq 0$ and $\pi_1,\ldots,\pi_i\in {\mathcal P}_n$. Let also $L\in \N$. By the strong Markov property applied to $\A^L$ at time $\rho_Le_{i-1}^L$, we have \setlength\arraycolsep{1pt}
\begin{eqnarray*}
 \prob_{a_L}\big[ \big(e_1^L,\ldots,e_{i-1}^L\big)&\in& D;\ e_i^L-e_{i-1}^L\leq t;\ \A^{L,u}(e_1^L)=\pi_1,\ldots, \A^{L,u}(e_i^L)=\pi_i\big]\\
&=& \E_{a_L}\Big[\ \mathbf{1}_{\{(e_1^L,\ldots,e_{i-1}^L)\in D\}}\ \mathbf{1}_{\{\A^{L,u}(e_1^L)=\pi_1,\ldots, \A^{L,u}(e_{i-1}^L)=\pi_{i-1}\}} \\
& &\qquad\qquad\qquad \times \prob_{\A^L(\rho_Le_{i-1}^L)}\big[e_1^L\leq t;\ \A^{L,u}(e_1^L)=\pi_i\big]\Big].
\end{eqnarray*}
First, using arguments analogous to those leading to Lemma \ref{lemm lineages far}, up to an error term vanishing uniformly in $(a_L)_{L\in \N}$ such that $a_L\in \Gamma(L,n)$ for every $L\in \N$, we can consider that $\A^L(\rho_Le_{i-1}^L)\in \GA(L,n)$. As $\mathrm{bl}\big(\A^L(\rho_Le_{i-1}^L)\big)=\pi_{i-1}$ for each $L$, we can therefore use $H(1)$ to write that
$$
\lim_{L\rightarrow \infty}\prob_{\A^L(\rho_Le_{i-1}^L)}\big[e_1^L\leq t;\ \A^{L,u}(e_1^L)=\pi_i\big]= \prob_{\pi_{i-1}}\big[e_1\leq t;\ \Lambda^{(\beta,c)}(e_1)=\pi_i\big],
$$
and so dominated convergence and $H(i-1)$ give us \setlength\arraycolsep{1pt}
\begin{eqnarray*}
& &\lim_{L\rightarrow \infty}\prob_{a_L}\big[\big(e_1^L,\ldots,e_{i-1}^L\big)\in D;\ e_i^L-e_{i-1}^L\leq t;\ \A^{L,u}(e_1^L)=\pi_1,\ldots, \A^{L,u}(e_i^L)=\pi_i\big]\\
& &= \E_{\pi_0}\Big[\ \mathbf{1}_{\{(e_1,\ldots,e_{i-1})\in D\}}\ \mathbf{1}_{\{\Lambda^{(\beta,c)}(e_1)=\pi_1,\ldots, \Lambda^{(\beta,c)}(e_{i-1})=\pi_{i-1}\}} \prob_{\pi_{i-1}}\big[e_1\leq t;\ \Lambda^{(\beta,c)}(e_1)=\pi_i\big]\Big]\\
& & = \prob_{\pi_0}\big[\big(e_1,\ldots,e_{i-1}\big)\in D;\ e_i-e_{i-1}\leq t;\ \Lambda^{(\beta,c)}(e_1)=\pi_1,\ldots, \Lambda^{(\beta,c)}(e_i)=\pi_i\big],
\end{eqnarray*}
which again yields $H(i)$ by standard arguments. The induction is now complete, and so we can conclude that the finite-dimensional distributions of the embedded Markov chain and the holding times of $\A^{L,u}$ under $\prob_{a_L}$ converge as $L\rightarrow \infty$ towards those of $\Lambda^{(\beta,c)}$ under $\prob_{\pi_0}$. The proof of $(b)$ is then complete.

To finish, suppose that $\rho_L\gg L^2\log L$. Then, we can find a sequence $\Phi_L$ increasing to $+\infty$ such that
$$
\sup_{A\in \GA(L,n)}\prob_A[\ \mathrm{a\ large\ event\ affects\ at\ least\ one\ lineage\ before\ time\ }\Phi_LL^2\log L]\rightarrow 0
$$
as $L\rightarrow \infty$. Hence, we can couple $\A^L$ with the process $\tilde{\A}^L$ which experiences only small events, so that the time by which they differ at step $L$ is larger than $\Phi_L$ with probability tending to one, uniformly in the sequence $(A_L)_{L\geq 1}$ chosen as above. By the results obtained in Section \ref{alpha<1} with $\rho_L\equiv +\infty$, we know that $\tilde{\A}^{L,u}$ converges in distribution towards $\mathcal{K}$, as a process in $D_{{\mathcal P}_n}[0,\infty)$. Since the sample size $n$ is finite and under Kingman's coalescent, a sample of $n$ lineages reaches a common ancestor in finite time almost surely, $(c)$ follows.$\hfill\square$

\appendix
\section{Proofs of the results of Section \ref{levy processes}}\label{appendix 1}
Since the proofs of Lemmas \ref{lemma entrance_levy} and \ref{lemma local TCL} are highly reminiscent of those of Theorem~2 and Lemma~3.1 in Cox \& Durrett~(2002), we shall only give the arguments we need to modify and refer to their paper for more extensive proofs.

\medskip \noindent \emph{Proof of Lemma \ref{lemma local TCL}: }Since $\ell^L$ is a L\'evy process, for any integers $n$ and $L$ one can decompose $\ell^L(n)$ into
$$
\ell^L(n)=\ell^L(0)+\sum_{k=1}^n\{\ell^L(k)-\ell^L(k-1)\},
$$
where the $n$ terms in the sum are i.i.d. random variables whose common distribution is that of $\ell^L(1)$ under $\prob_0$. Using Bhattacharya's local central limit theorem (see Theorem~1.5 in Bhattacharya~1977) \nocite{bhattacharya:1977} and the boundedness assumption on $\E_0[|\ell^L(1)|^4]$, we can control the deviation of $p^L(x,n)$ from the corresponding probabilities for Brownian motion up to an error of order $o(n^{-1})$ independent of $L$. Following Cox and Durrett's arguments, we obtain the desired results for integer times. For arbitrary times $t$, the Markov property applied to $\ell^L$ at time $\integ{t}$ (plus, for (d), the fact that the variations of $\ell^L$ are bounded on a time interval $[n,n+1]$) completes the proof.$\hfill\square$
\medskip

\noindent \emph{Proof of Lemma \ref{lemma entrance_levy}: }To simplify notation, we shall write $T(d_L)$ instead of $T(d_L,\ell^L)$ in the rest of the proof. For every $L\geq 1$, $x\in \T(L)$ and $\lambda>0$, let us define the following quantities :
\begin{eqnarray*}
F_L(x,\lambda)&=& \E_x\big[\exp(-\lambda T(d_L))\big],\\
G_L(x,\lambda)&=&\int_0^{\infty}e^{-\lambda t}p^L(x,t) dt=\E_x\bigg[\int_0^{\infty}e^{-\lambda t}\ \ind{\ell^L(t)\in B(0,d_L)}dt\bigg].
\end{eqnarray*}
Applying the strong Markov property to $\ell^L$ at time $T(d_L)$ and using a change of variables, we obtain (for any $x_L$)
\begin{equation}
G_L(x_L,\lambda) = \E_{x_L}\Big[e^{-\lambda T(d_L)}G_L\big(\ell^L(T(d_L)),\lambda\big)\Big].\label{eq lemma entrance levy}
\end{equation}
From the results of Lemma \ref{lemma local TCL}, we can derive the asymptotic behaviour of $G_L(x_L,\lambda)$. To this end, let $(v_L)_{L\geq 1}$ and $(u_L)_{L\geq 1}$ be two sequences growing to infinity such that $v_L(\log L)^{-1/2}\rightarrow 0$ and $u_L(\log L)^{-1}\rightarrow 0$ as $L\rightarrow \infty$. Splitting the integral in the definition of $G_L\big(x_L,\frac{\lambda}{L^2\log L}\big)$ into four pieces, we obtain first by $(b)$ of Lemma \ref{lemma local TCL} \setlength\arraycolsep{1pt}
\begin{eqnarray*}
\frac{1}{d_L^2\log L}\int_{v_LL^2}^{\infty} &\exp&\left(-\frac{\lambda t}{L^2\log L}\right)p^L(x_L,t)dt \\
&=&\frac{1}{d_L^2\log L}\int_{v_LL^2}^{\infty} \exp\left(-\frac{\lambda t}{L^2\log L}\right)\frac{\pi d_L^2}{L^2}\ (1+\delta_{L,1})\ dt \\
&=&\frac{\pi}{\lambda}\ \exp\left(-\frac{\lambda v_L}{\log L}\right)(1+\delta_{L,1}) = \frac{\pi}{\lambda}\ (1+\delta'_{L,1})
\end{eqnarray*}
as $L\rightarrow \infty$, where $\delta_{L,1}, \delta'_{L,1}\rightarrow 0$ uniformly in $x\in \T(L)$. By $(a)$ of Lemma \ref{lemma local TCL}, we have
\begin{eqnarray*}
\frac{1}{d_L^2\log L}\int_{\e_LL^2}^{v_LL^2} \exp\left(-\frac{\lambda t}{L^2\log L}\right)p^L(x_L,t)\ dt\ &\leq&\ \frac{1}{d_L^2\log L}\frac{C_1d_L^2}{\integ{L^2\e_L}}\ v_LL^2\\
&\sim&\ \frac{C_1v_L}{\sqrt{\log L}}\rightarrow 0,\quad \mathrm{as\ } L\rightarrow \infty
\end{eqnarray*}
by our assumption on $v_L$. By $(c)$ of Lemma \ref{lemma local TCL},\setlength\arraycolsep{1pt}
\begin{eqnarray*}
\frac{1}{d_L^2\log L}& &\int_{u_L(1+|x_L|^2\vee d_L^2)}^{\e_LL^2} \exp \left(-\frac{\lambda t}{L^2\log L}\right)p^L(x_L,t)\ dt \\
&=& \frac{1}{d_L^2\log L}\int_{u_L(1+|x_L|^2\vee d_L^2)}^{\e_LL^2}\frac{d_L^2}{2\sigma_L^2t}\ (1+\delta_{L,2})\ dt \\
&=& \frac{1}{2\sigma_L^2\log L}\ \Big(2\log L-\log(1+|x_L|^2\vee d_L^2)+\log\e_L - \log u_L\Big)(1+\delta_{L,2}) \\
&= & \frac{1-\beta\vee \gamma}{\sigma^2}\ (1+\delta'_{L,2}),
\end{eqnarray*}
whenever $\frac{\log^+|x_L|}{\log L}\rightarrow \beta$ as $L$ grows to infinity. Here again, $\delta_{L,2},\delta'_{L,2}\rightarrow 0$ uniformly in $x\in \T(L)$ as $L\rightarrow \infty$. Finally, by $(d)$ of Lemma \ref{lemma local TCL}, we can write
\begin{multline*}
\frac{1}{d_L^2\log L}\int_0^{u_L(1+|x_L|^2\vee d_L^2)} \exp\left(-\frac{\lambda t}{L^2\log L}\right)p^L(x_L,t)dt\\
\leq \frac{C_2}{d_L^2\log L}\frac{u_L(1+|x_L|^2\vee d_L^2)}{1+ d_L^{-2}|x_L|^2}\rightarrow 0,
\end{multline*}
independently of $(x_L)_{L\geq 1}$ since $d_L$ does not vanish and $u_L (\log L)^{-1}\rightarrow 0$.

Combining the above, we obtain that if $\frac{\log^+|x_L|}{\log L}\rightarrow \beta$, then
$$
\frac{1}{d_L^2\log L}\ G_L\Big(x_L,\frac{\lambda}{L^2\log L}\Big)= \frac{\pi}{\lambda}+\frac{1-(\beta \vee \gamma)}{\sigma^2}+ o(1),\qquad \mathrm{as}\ L\rightarrow \infty,
$$
where the remainder does not depend on $(x_L)_{L\geq 1}$. Coming back to (\ref{eq lemma entrance levy}) with $x_L \in \Gamma(L,1)$, the uniform convergence obtained above, together with the fact that $\ell^L(d_L)\in B(0,d_L)$ a.s. yield
\begin{equation}\label{conv to exponential}
\lim_{L\rightarrow \infty}\E_{x_L}\left[\exp\left(-\frac{\lambda \pi \sigma^2\ T(d_L)} {(1-\gamma)L^2\log L} \right)\right]=\frac{(1-\gamma)/(\sigma^2\lambda)}{(1-\gamma)/(\sigma^2\lambda)+(1-\gamma)/\sigma^2} =\frac{1}{1+\lambda},
\end{equation}
which we recognize as the Laplace transform of an $\mathrm{Exp}(1)$-random variable. Since the left-hand side of (\ref{conv to exponential}) is monotone in $\lambda$ and the function $\lambda\mapsto (1+\lambda)^{-1}$ is continuous, this convergence is in fact uniform in $\lambda \geq 0$. By standard approximation arguments (see for instance the proof of Theorem 4 in \nocite{cox:1989} Cox~1989), we obtain that for any fixed $t>0$,
$$
\lim_{L\rightarrow \infty}\sup_{x_L\in \Gamma(L,1)}\left|\prob_{x_L}\bigg[\frac{\lambda \pi \sigma^2} {(1-\gamma)L^2\log L}\ T(d_L)>t\bigg]-e^{-t}\right|=0,
$$
and, by monotonicity and the fact that all the quantities involved tend to $0$ as $t\rightarrow \infty$, this convergence is uniform in $t\geq 0$. The interested reader will find all the missing details in the appendix of Cox \& Durrett (2002). $\hfill\square$
\medskip

\noindent \emph{Proof of Lemma \ref{lemm no entrance}: }Let $x\in \T(L)$ and $(U_L')_{L\in \IN}$ be as in the statement of Lemma \ref{lemm no entrance}. Using the strong Markov property at time $T(R,\ell^L)$, we can write
\begin{multline}\label{eq lemm no entrance}
\prob_x\big[\ell^L(U_L'+u_L)\in B(0,R)\big]\\
\geq \int_{U_L'-u_L}^{U_L'}\int_{B(0,R)}\prob_x\big[T(R,\ell^L)\in ds, \ell^L(s)\in dy\big]\prob_y\big[\ell^L(U_L'+u_L-s)\in B(0,R)\big].\phantom{A}
\end{multline}
Note that, on the right-hand side of (\ref{eq lemm no entrance}), the quantity $U_L'+u_L-s$ lies in $[u_L,2u_L]$. We assumed that $2u_L\leq L^2(\log L)^{-1/2}$, and so we can use $(c)$ of Lemma \ref{lemma local TCL} with $d_L\equiv R$ and write
$$
\lim_{L\rightarrow \infty}\ \sup_{y\in B(0,R)}\ \sup_{u_L\leq t\leq 2u_L}\ \bigg|\frac{2\sigma^2_L t}{R^2}\ \prob_y\big[\ell^L(t)\in B(0,R)\big]-1\bigg| = 0,
$$
which gives us the existence of a constant $C_0$ and of an index $L_0$ such that for each $L\geq L_0$, $y\in B(0,R)$ and $t\in [u_L,2u_L]$,
$$
\prob_y\big[\ell^L(t)\in B(0,R)\big] \geq \frac{C_0}{t}\geq \frac{C_0}{2u_L}.
$$
Furthermore, since $U_LL^{-2}\rightarrow \infty$, we can use $(b)$ of Lemma \ref{lemma local TCL} to obtain the existence of $L_1\in \IN$ and a constant $C_1>0$ depending only on $(U_L)_{L\geq 1}$ such that for every $L\geq L_1$,
$$
\sup_{t\geq U_L}\sup_{y\in \IT(L)}\bigg|\prob_y\big[\ell^L(t+u_L)\in B(0,R)\big]- \frac{\pi R^2}{L^2}\bigg|\leq \frac{C_1}{L^2}.
$$
Using these two inequalities in (\ref{eq lemm no entrance}), we have for $L$ large enough and for all $x\in \T(L)$
$$
\frac{C_1+\pi R^2}{L^2}\geq \prob_x\big[T(R,\ell^L)\in [U_L'-u_L,U_L']\big]\times \frac{C_0}{2u_L},
$$
which gives us the desired result.$\hfill\square$

\section{Proof of the technical points of Section \ref{alpha<1}}\label{appendix 2}
\noindent \emph{Proof of Lemma \ref{lemm excursion}: } Let us start
with the case $\rho_L=\mathcal{O}(\psi_L^2)$ as $L\rightarrow
\infty$. The rate of decay of the probability of a long excursion is
known for simple random walks and Brownian motion (see
Ridler-Rowe~1966), and so \nocite{ridler-rowe:1966} the proof of
Proposition \ref{prop gathering} suggests that we should consider
the process $\hat{\ell}^L\equiv \psi_L^{-1}X^L(\rho_L\cdot)$. But
$\hat{\ell}^L$ here is not a L\'evy process, since $X^L$ is the
difference of the locations of two lineages whose motions are not
independent in $B(0,2R^B\psi_L)$. However, it is not difficult to
convince oneself that for each $y\in B(0,(7/4)R^B)^c$, the return
time into $B(0,(3/2)R^B)$ of $\hat{\ell}^L$ starting at $y$ is
smaller than or equal to that of $\ell^L$ defined as the rescaled
process $\psi_L^{-1}\xi^L(\rho_L\cdot)$ also starting at $y$.
Indeed, the rate at which reproduction events affect at least one of
the lineages is bounded from below by the rate at which a single
lineage is affected, the distribution of the jumps of $\hat{\ell}^L$
and $\ell^L$ are identical outside $B(0,2R^B)$ and inside this ball,
coalescence events make it easier for $\hat{\ell}^L$ to enter
$B(0,(3/2)R^B)$. Hence, we shall establish the desired bound for
$\ell^L$. In addition, we shall consider that $\ell^L$ evolves on
$\R^2$ instead of $\T(L)$, since the return time here can only
increase with the available space.

For each $L\in \N$, set $\sigma_0^L=0$ and let
$(\sigma_i^L)_{i\in\N}$ be the sequence of jump times of $\ell^L$.
Let $\rho_L\theta_s$ (resp., $\theta_B$) be the jump rate of
$\ell^L$ due to small events (resp., due to large events). The
quantities $\theta_s$ and $\theta_B$ do not depend on $L$ since
$\mu^B,\ \mu^s$ and the probability measures $\nu_r^{B,s}$ do not.
For each $t\geq 0$, we have
$\ell^L(t)=\ell^L(0)+\sum_{i:\sigma_i^L\leq
t}\big\{\ell^L(\sigma_i^L)-\ell^L(\sigma_{i-1}^L)\big\}$, where
$\big(\ell^L(\sigma_i^L)-\ell^L(\sigma_{i-1}^L)\big)_{i\in \N}$ is a
sequence of i.i.d. random variables with covariance matrix of the form
$\upsilon_L \mathrm{Id}$. Using the distribution of a single small or large jump and the fact that a given jump is a small one with probability $\theta_s\rho_L/(\theta_s\rho_L+\theta_B)$, we easily check that there exists $V>0$, independent of $L$, such that $\upsilon_L\sim V/\rho_L$ as $L\rightarrow \infty$ (recall our assumption $\rho_L=\mathcal{O}(\psi_L^2)$).

Let $x$ in $B(0,4R^B)\setminus B(0,(7/4)R^B)$ and let $W$ be a
two-dimensional Brownian motion starting at $x$. For each $L\in \N$,
by the Skorohod Embedding Theorem (see, e.g., Billingsley~1995) one
can \nocite{billingsley:1995} construct a sequence $(s_i^L)_{i\in
\N}$ of stopping times such that the $W(s_i^L)$ have the same joint
distributions as the $\ell^L(\sigma_i^L)$ : for every $i\in \N$,
conditionally on $W(s_{i-1}^L)$, $s_i^L$ is the first time greater
than $s_{i-1}^L$ at which $W$ leaves
$B\big(W(s_{i-1}^L),r_i^L\big)$, where $r_i^L$ is a random variable
independent of $W$ and of $\{s_j^L,j<i\}$ having the same
distribution as the length of the first jump of $\ell^L$. Now, we
claim that there exists $\gamma>0$ independent of $L$ and $x$, such that
each time $W$ visits $B(0,R^B/2)$ and then leaves $B(0,(3/2)R^B)$,
the probability that one of the $s_i^L$'s falls into the
corresponding period of time that $W$ spends within $B(0,(3/2)R^B)$
is at least $\gamma$. Indeed, set $T_0(W)=\breve{T}_0(W)=0$ and define
the sequences of stopping times $\{T_k(W),k\geq 1\}$ and
$\{\breve{T}_k(W),k\geq 1\}$ by induction in the following manner:
\begin{eqnarray*}
T_k(W)&=&\inf\big\{t> \breve{T}_{k-1}(W):\ W(t)\in B(0,R^B/2)\big\}, \\
\breve{T}_k(W)&=& \inf\big\{t> T_k(W):\ W(t)\notin
B(0,(3/2)R^B)\big\}.
\end{eqnarray*}
(Note that each $T_k$ is a.s. finite due to the recurrence of two-dimensional Brownian motion.)
Then for each $k\in \N$, if $j$ is the index of the last $s_i^L$
before $T_k(W)$ and $s_j^L$ corresponds to a small event, by
construction we have $\big|W(s_j^L)-W(T_k(W))\big|<2R^s\psi_L^{-1}$
and so $W(s_j^L)\in B(0,(3/2)R^B)$ for $L$ large enough. If $s_j^L$
is due to a large event and $W(s_j^L)\notin B(0,(3/2)R^B)$, then
necessarily $W(s_j^L)\in B(0,(5/2)R^B)$. But the exit point from a
ball $B$ of Brownian motion started at the centre of this ball is
uniformly distributed over the boundary of $B$, and so one can
define $\gamma$ as the minimum over $(y,r)$ with $|y|\geq 3R^B/2$ and
$|y|-R^B/2<r\leq 2R^B$ of the probability that $W$ started at $y$
escapes $B(y,r)$, through the part of its boundary which lies within
$B(0,(3/2)R^B)$. Hence, if we define for each $t\geq 0$ the random
variable $N(t)$ as the maximal integer $k$ such that
$\breve{T}_k(W)\leq t$, we can write for each $L$
$$
\prob_{\psi_Lx}\big[q_1^L>\rho_Lu\big]= \prob_x\big[\ell^L(\sigma_j^L)\notin B(0,(3/2)R^B),\ \forall \ j\leq i(u,L)\big]  \leq \E_x\big[(1-\gamma)^{N(s^L_{i(u,L)})}\big],
$$
where $i(u,L)=\max\{j:\ \sigma_j^L\leq u\}$. Since $N$ is a.s. a
non-decreasing function of $t$, we have for any given $m\in \R_+$
\begin{equation}\label{decay q1}
\prob_{\psi_Lx}\big[q_1^L>\rho_Lu\big]\leq
\E_x\big[(1-\gamma)^{N(mu)}\big]+\prob_x\big[s^L_{i(u,L)}<mu\big].
\end{equation}
Now, $i(u,L)$ is the number of points of the Poisson point processes
$\Pi_L^s$ and $\Pi_L^B$ which fall into the time interval
$[0,u\rho_L]$ on the original timescale, it is therefore a Poisson
random variable with parameter $u(\theta_s\rho_L+\theta_B)$. If
$a>0$, then by the Markov inequality
$$
\prob_x\big[i(u,L)\leq au\theta_s\rho_L\big]\leq
e^{au\theta_s\rho_L}\E\big[e^{-i(u,L)}\big]=
\exp\big\{u\theta_s\rho_L(a+e^{-1}-1)+u\theta_B(e^{-1}-1)\big\},
$$
so that this quantity converges exponentially fast to $0$ for $a>0$
small enough. On the event $\{i(u,L)>au\theta_s\rho_L\}$,
$s^L_{i(u,L)}$ is the sum of at least $au\theta_s\rho_L$ i.i.d.
random variables, each of which corresponds to the exit time of Brownian motion from a ball of radius at most $2R^s/\psi_L$ with probability $\theta_s\rho_L/(\theta_s\rho_L+\theta_B)$ and to the exit time of Brownian motion from a ball of radius at most $2R^B$ otherwise. Therefore, one can find $V'>0$ independent of $L$ such that $\E[s_1^L]\sim V'\rho_L^{-1}$ as $L\rightarrow \infty$. Using the same technique as above then gives us that for $m>0$ small
enough, there exists $\kappa(m)>0$ and $L(m)\in \N$ such that for all $L\geq L(m)$ and $u\geq 0$,
$$
\prob_x\big[i(u,L)>au\theta_s\rho_L,\ s^L_{i(u,L)}<mu\big]\leq e^{-\kappa(m)\rho_Lu}.
$$
Let us now prove that
$$
\prob_x[N(mu)\leq K\log\log u]\leq C\ \frac{\log\log u}{\log u}
$$
for a constant $C>0$ independent of $x\in B(0,4R^B)\setminus
B(0,(7/4)R^B)$ and $u$ large enough (again independently of $x$).
The reasoning is identical to that made to arrive at (\ref{ineg
k_L}), with $q_i^L$ (resp., $Q_i^L$) replaced by $T_i(W)$ (resp.,
$\breve{T}_i(W)$). Using the fact that $C_1\equiv \sup_{x\in
B(0,R^B/2)}\E_x[\breve{T}_1(W)]<\infty$ and
\begin{equation}\label{excursion Brownian}
\sup_{y\in B(0,4R^B)}\prob_x\big[T_1(W)>u\big]\leq \frac{C_2}{\log
u}
\end{equation}
for a constant $C_2$ and $u$ large enough (see Theorem 2 in
Ridler-Rowe~1966), we can conclude that for each $x\in
B(0,4R^B)\setminus B(0,(7/4)R^B)$, and $u$ large enough,
$$
\prob_x\big[N(mu)\leq \log\log u\big]\leq \frac{2C_1(\log\log
u)^2}{mu}+\frac{C_2\log\log u}{\log\big(mu/(2\log\log u)\big)}\leq
\frac{C'\log \log u}{\log u},
$$
again for $C'>0$ and $u$ large enough independently of $x$. Coming
back to (\ref{decay q1}), we obtain for a constant $C''>0$ and for
all $x\in B(0,4R^B)\setminus B(0,(7/4)R^B)$,
\begin{eqnarray}
\prob_{\psi_Lx}\big[q_1^L>\rho_Lu\big]& \leq &(1-\gamma)^{\log\log u} +
\prob_x\big[N(mu)\leq \log \log u\big]+e^{-C''u}\nonumber \\
&\leq & (1-\gamma)^{\log\log u}+ \frac{C' \log\log u}{\log
u}+e^{-C''u}.\label{decay excursion}
\end{eqnarray}
Define $g(u)$ as the expression on the right-hand side of
(\ref{decay excursion}) to obtain the result.

When $\psi_L^2\rho_L^{-1}\rightarrow 0$, the probability that a
large event occurs by time $u\psi_L^2$ is given by
$$
1-\exp\Big\{-\theta_Bu\frac{\psi_L^2}{\rho_L}\Big\}\rightarrow 0
\qquad \mathrm{as\ }L\rightarrow \infty.
$$
On the event that no large events occur by time $u\psi_L^2$, the
first visit of $W$ into $B(0,R^B/2)$ will produce a time $s_i^L$
such that $W(s_i^L)\in B(0,(3/2)R^B)$ with probability $1$ for the
reason expounded above, and so the first term on the right-hand side
of (\ref{decay q1}) is now the probability that $T_1(W)$ is greater
than $mu$. The inequality in (\ref{excursion Brownian}) and the
exponential decay of $\prob_x[s_{i(u,L)}^L<mu]$ now imply the
result.$\hfill\square$

\medskip
\noindent \emph{Proof of Lemma \ref{lemm time in B}: } The arguments
are slightly different according to whether $\rho_L\psi_L^{-2}$ is
bounded or tends to infinity as $L\rightarrow \infty$. Let us
consider the first case. Recall the definition of
$\rho_L^{-1}\theta_B$ given in the proof of Lemma \ref{lemm
excursion} as the maximal rate at which a lineage is affected by a
large event. The coalescence rate of two lineages is then bounded by
$2\rho_L^{-1}\theta_B$, regardless of their locations. By our
assumption (\ref{condition nu}), there exist $r\in (0, R^B)$ and
$\delta>0$ such that $\mathrm{Leb}\big(\big\{r'\in [r,r+\delta]:\
\nu_{r'}^B\notin \{\delta_0,\delta_1\}\big\}\big)>0$. We shall use
these events to send the two lineages at distance at least
$(7/4)R^B\psi_L$ from each other, whatever their initial separation
was. The proof is quite natural, so we just give the main arguments.
If only large jumps occurred, then if a sequence of at least
$7R^B/(2r)$ large events increased $|X^L|$ by at least
$(r/2)\psi_L$ each before the first coalescence happened, $X^L$
starting within $B(0,(3/2)R^B\psi_L)$ would certainly leave
$B(0,(7/4)R^B\psi_L)$. Moreover, a large event affecting $X^L$ and
conditioned on not leading to a coalescence biases the jump towards
increasing $|X^L|$ (we do not allow some centres that are too close
to both lineages). This remark and (\ref{condition nu}) guarantee
that the rate at which these separating events occur (that is, events increasing $|X^L|$ by at least $(r/2)\psi_L$) is bounded from
below by $\rho_L^{-1}\theta_{\mathrm{sep}}$, where
$\theta_{\mathrm{sep}}$ is a positive constant. The total rate at
which large events affect $X^L$ is bounded by
$2\rho_L^{-1}\theta_B$, and so there is a positive probability
$p_{\mathrm{sep}}$, independent of the starting point of $X^L$, that
$X^L$ leaves $B(0,(7/4)R^B\psi_L)$ before coming back to $0$
through a (large) coalescence event. As regards the effect of small
events, recall that we assumed that $\rho_L\psi_L^{-2}$ is
bounded. Hence, the probability that $X^L$ starting from
$B(0,r\psi_L/2)^c$ does not enter $B(0,2R^s)$ after a time of order
$\mathcal{O}(\rho_L)$ only through small jumps is bounded from
below and by the symmetry of these small jumps, with probability at
least $1/2$ the radius of $X^L$ increases between two large jumps.
Hence, up to modifying $p_{\mathrm{sep}}$ to take into account the
effect of the small jumps, the probability that $X^L$ leaves
$B(0,(7/4)R^B\psi_L)$ before coming back to $0$ is still bounded
from below by $p_{\mathrm{sep}}>0$.

By the definition of $R^B$ and Assumption (\ref{coal at boundary}), large events of size close to $R^B$ occur at a positive
rate and lead to the coalescence of the lineages with positive
probability, so that the waiting time for the coalescence of two
lineages at distance at most $(7/4)R^B\psi_L$ is bounded by
$\rho_L$ times an exponential with positive parameter $\gamma$.
This gives us that $\rho_L^{-1}Q_1^L$ is stochastically bounded by
$\sum_{i=1}^k N_i$, where $k$ is geometric with success probability
$p_{\mathrm{sep}}>0$ and $\{N_i,i\in \N\}$ is a sequence of i.i.d.
$\mathrm{Exp}(\gamma)$ random variables, all of them independent of
the initial value $x\in B(0,(3/2)R^B\psi_L)$ of $X^L$. We can
therefore choose $C_Q=(\gamma p_{\mathrm{sep}})^{-1}$.

When $\rho_L^{-1}\psi_L^2\rightarrow 0$, if we use the same
reasoning as above there is a positive probability that a large
event separates the two lineages at distance at least $r\psi_L$,
regardless of their separation just before this event. In addition,
the rate of these separating events is at least equal to
$\rho_L^{-1}\theta_{\mathrm{sep}}>0$. Between two large events,
$X^L$ only does small jumps, and as long as $X^L\notin B(0,2R^s)$,
the Skorohod Embedding Theorem (see the proof of Lemma \ref{lemm
excursion}) enables us to assert that $X^L$ will leave
$B(0,(7/4)R^B\psi_L)$ in a time of order $\mathcal{O}(\psi_L^2)$.
Moreover, for $\e>0$ small, the same argument shows that the
probability that $X^L$ leaves $B(0,(7/4)R^B\psi_L)$ before entering
$B(0,\e\psi_L)$ is bounded from below by a constant
$p_{\mathrm{esc}}>0$ independent of $L$ and of the value $y\in
B(0,r\psi_L)^c$ of $X^L$ just after the large jump described above.
A fortiori, $p_{\mathrm{esc}}$ is also a lower bound for the
probability that $X^L$ started at $y$ leaves $B(0,(7/4)R^B\psi_L)$
before entering $B(0,2R^s)$ only through small jumps, and so we
obtain that between two large events such that the first large jump
sends (or keeps) $X^L$ out of $B(0,r\psi_L)$, $X^L$ escapes
$B(0,(7/4)R^B\psi_L)$ with probability at least $p_{\mathrm{esc}}$
(recall that the total rate of large events affecting at least one
of the lineages is bounded by $2\theta_B\rho_L^{-1}$ and
$\rho_L\gg \psi_L^2$). Consequently, $Q_1^L$ is this time
stochastically bounded by $\sum_{i=1}^k N_i(L)$, where $k$ is a
geometric random variable with success probability
$p_{\mathrm{esc}}>0$ and for each $L\in \N$, $\{N_i(L),i\in \N\}$ is
a sequence of i.i.d.
$\mathrm{Exp}(\rho_L^{-1}\theta_{\mathrm{sep}})$ random variables,
all of them independent of the initial value $x\in
B(0,(3/2)R^B\psi_L)$ of $X^L$. The desired result follows, with
$C_Q=(\theta_{\mathrm{sep}}p_{\mathrm{esc}})^{-1}.$ $\hfill\square$

\medskip
\noindent \emph{Proof of Lemma \ref{lemm proba coal}: }The
inequality in (\ref{hitting prob 1}) can be restated as in
(\ref{hitting prob 2}) (the quantity inside the brackets then tends
to $1$), so we prove both inequalities using this form. Let
$\theta_c$ be such that $\rho_L^{-1}\theta_c$ is the minimum rate
at which two lineages at distance at most $(1+\delta)R^B\psi_L$
from each other coalesce (where $\delta>0$ is defined at the
beginning of the proof of Lemma \ref{lemm time in B}). By the
definition of $R^B$ and assumption (\ref{coal at boundary}), the
rate at which a reproduction event of radius $r\in
\big[R^B(1-\frac{\delta}{4})\psi_L,R^B\psi_L\big]$ occurs and
leads to the coalescence of the lineages does not vanish as $L$
tends to infinity (when multiplied by $\rho_L$), and so
$\theta_c>0$. Let us show that if $\eta>0$ is small enough, the
probability that $X^L$ starting within $B(0,R^B\psi_L)$ does not
leave $B(0,(1+\delta)R^B\psi_L)$ through only small jumps by time
$\eta\psi_L^2$ is bounded from below by a positive constant,
independent of $L$ large. The term inside the brackets in
(\ref{hitting prob 2}) will then come from the probability that a
large event occurs before time $\eta\psi_L^2$ and the first such event leads
to the coalescence of the lineages (i.e., a jump onto $0$ for
$X^L$).

Let $\eta>0$ and $x\in B(0,R^B\psi_L)$, and let $\tau_B^L$ denote
the epoch of the first large event affecting $X^L$. By the argument
given above, the probability that $X^L$ starting at $x$ hits $0$
before leaving $B(0,(1+\delta)R^B\psi_L)$ is bounded from below by
the probability that $X^L$ started at $x$ stays within this ball
until $\tau_B^L$, $\tau_B^L$ is less than or equal to
$\eta\psi_L^2$ and the first large event leads to the coalescence
of the lineages. Writing $\mathcal{E}_{L,\eta}$ for the event that $X^L$
stays within $B(0,(1+\delta)R^B\psi_L)$ before $\tau_B^L$ and
$\tau_B^L\leq \eta\psi_L^2$, this probability is equal to
\begin{equation}\label{eq hitting 1}
\prob_x\big[\mathrm{the\ first\ large\ event\ is\ a\ coalescence\
}|\ \mathcal{E}_{L,\eta}\big]\prob_x\big[\mathcal{E}_{L,\eta}\big].
\end{equation}
If, for each $L\in \N$, $\rho_L^{-1}E<\infty$ denotes the rate at
which a single lineage on $\T(L)$ is affected by a large
reproduction event, then the rate at which at least one of two
lineages are affected is bounded by twice this quantity, and so the
first probability in (\ref{eq hitting 1}) is bounded from below by
$\theta_c/(2E)$. Now, $X^L$ experiences no large reproduction event
before time $\tau_B^L$, and so we can again use the equality in
distribution stated in the proof of Proposition \ref{prop coal time}
$(b)$ (we also keep the notation introduced there). Write
$t_{\mathrm{exit}}$ for the first time $X^L$ leaves
$B(0,(1+\delta)R^B\psi_L)$, and $\tilde{t}_{\mathrm{exit}}$ for the
corresponding time for $\tilde{X}^L$ (which sees only small events).
We have
\begin{eqnarray}
\prob_x\big[\mathcal{E}_{L,\eta}\big]&=&\prob_x\big[t_{\mathrm{exit}}\geq
\tau_B^L;\ \tau_B^L\leq \eta\psi_L^2\big] \nonumber\\
&=& \prob_x\big[\tilde{t}_{\mathrm{exit}}\geq e(\tilde{X}^L);\
e(\tilde{X}^L)\leq \eta\psi_L^2\big] \nonumber\\
&\geq & \prob_x\big[\tilde{t}_{\mathrm{exit}}\geq \eta \psi_L^2;\
e(\tilde{X}^L)\leq \eta\psi_L^2\big] \nonumber\\
&=&\prob_x\big[e(\tilde{X}^L)\leq \eta\psi_L^2\ \big|\
\tilde{t}_{\mathrm{exit}}\geq \eta
\psi_L^2\big]\prob_x\big[\tilde{t}_{\mathrm{exit}}\geq \eta
\psi_L^2\big]. \label{lower bound hitting}
\end{eqnarray}
Since a pair of lineages is affected by a large event at rate at
least $\rho_L^{-1}E$, the first probability on the right-hand side of (\ref{lower bound
hitting}) is bounded below for all $x\in B(0,R^B\psi_L)$ by
$$
1-\exp\bigg\{-\eta E\frac{\psi_L^2}{\rho_L}\bigg\}.
$$
Now, if $\tilde{X}^L$ starts within $B(0,R^B\psi_L)$, it needs to
cover a distance of at least $\delta R^B\psi_L$ to exit
$B(0,(1+\delta)R^B\psi_L)$. Furthermore, coalescence events tend to
keep $\tilde{X}^L$ within $B(0,(1+\delta)R^B\psi_L)$, and so for
each $x$, the second probability on the right-hand side of (\ref{lower bound hitting}) is
larger than $\prob_0[\hat{t}_{\mathrm{exit}}\geq \eta\psi_L^2]$,
where $\hat{t}_{\mathrm{exit}}$ is the exit time from $B(0,\delta
R^B\psi_L)$ of the process $\{\hat{\xi}^L(2t),t\geq 0\}$ which
experiences only small jumps. Decomposing this L\'evy process into
the sum of its jumps and applying Doob's maximal inequality to the
submartingale $|\hat{\xi}^L|^2$, we obtain
$$
\prob_0\bigg[\sup_{0\leq t\leq
\eta\psi_L^2/2}|\hat{\xi}^L(2t)|^2>(\delta R^B\psi_L)^2\bigg]\leq
\frac{1}{(\delta R^B\psi_L)^2}\ \E_0\big[|\hat{\xi}^L(\eta
\psi_L^2)|^2\big] = \frac{2\eta \sigma_s^2}{\delta^2(R^B)^2},
$$
where the last equality comes from (\ref{PPP moment 2}). Choosing
$\eta>0$ small enough so that the quantity above is less than $1$,
we obtain that for all $x\in B(0,R^B\psi_L)$
$$
\prob_x\big[\tilde{t}_{\mathrm{exit}}\geq \eta \psi_L^2\big]\geq
\prob_0[\hat{t}_{\mathrm{exit}}\geq \eta\psi_L^2]\geq 1-\frac{2\eta
\sigma_s^2}{\delta^2(R^B)^2}\equiv \theta_4>0.
$$
Combining the above and choosing $\theta_2=\theta_4 \theta_c/(2E)$
and $\theta_3=\eta E$, we obtain (\ref{hitting prob 2}).$\hfill\square$
\medskip

\noindent\emph{Proof of Lemma \ref{lemm lineages far}: }
If we were considering the times $\tau_{ij}$ rather than
$\tau^*_{ij}$, Lemma \ref{lemm lineages far} would follow from the
same arguments as in Cox \& Griffeath~(1986)
\nocite{cox/griffeath:1986} (see Lemma 1). Here, we have to work a
bit harder and decompose the event in (\ref{lineage distribution 1})
into more cases. Recall the definition of $\varpi_L$ given in the statement of Theorem \ref{result alpha<1}. For each $L\in \N$, the probability in
(\ref{lineage distribution 1}) is bounded by
\setlength\arraycolsep{1pt}
\begin{eqnarray}
& &\prob_{A_L}\Big[\tau <\frac{\vp_L}{\sqrt{\log L}}
\Big]+\prob_{A_L}\Big[\tau \geq\frac{\vp_L}{\sqrt{\log L}}\ ;\
\tau^*=\tau^*_{12}\ ;\ \tau\neq \tau_{12}\Big] \label{lineages far
1}\\ &+& \prob_{A_L}\Big[\tau^*=\tau^*_{12}\ ;\
\frac{\vp_L}{\sqrt{\log L}}\leq
\tau=\tau_{12}<\tau^*_{12}-\frac{\vp_L}{(\log
L)^2}\Big]\label{lineages far
2}\\
&+& \prob_{A_L}\Big[\tau^*=\tau^*_{12} ; \frac{\vp_L}{\sqrt{\log
L}}\leq \tau=\tau_{12} ; \tau_{12}\geq\tau^*_{12}-\frac{\vp_L}{(\log
L)^2} ; \exists i\in\{1,2\}, \tau_{i3}\in
(\tau_{12},\tau_{12}^*]\Big]\label{lineages far 3}\\ &+&
\prob_{A_L}\Big[\tau^*=\tau^*_{12} ; \frac{\vp_L}{\sqrt{\log L}}\leq
\tau=\tau_{12} ; \forall i\in \{1,2\}, \tau_{i3}>\tau^*_{12} ;
|\xi_1^L(\tau^*)-\xi^L_3(\tau^*)|\leq \frac{L}{\log
L}\Big].\phantom{AAA}\label{lineages far 4}
\end{eqnarray}
Suppose first that $\rho_L\ll \psi_L^2\log L$. The first term in
(\ref{lineages far 1}) is bounded by the sum over $i\neq j\in
\{1,\ldots,4\}^2$ of $\prob_{A_L}[\tau_{ij}<\vp_L(\log L)^{-1/2}]$,
which tends to $0$ uniformly in $A_L$ by Proposition \ref{prop gathering} and the
consistency of the genealogy. The quantity in (\ref{lineages far 2}), expressing the probability that the first pair to meet is the pair $(1,2)$ but then coalescence of these lineages takes longer than $\vp_L/(\log L)^2$ units of time, is therefore bounded by
$$
\prob_{A_L}\Big[\tau^*_{12}-\tau_{12}>\frac{\vp_L}{(\log L)^2}\Big],
$$
which converges to zero as $L\rightarrow \infty$, uniformly in $A_L$
(apply the strong Markov property at time $\tau_{12}$ and use $(a)$
of Proposition \ref{prop coal time}). The expression in
(\ref{lineages far 3}) corresponds to the event in which $(1,2)$ is the first pair to meet and ``quickly'' merge, but another pair of lineages manages to meet between $\tau_{12}$ and $\tau_{12}^*$. It is thus bounded by
\begin{eqnarray*}
\prob_{A_L}&\bigg[&\frac{\vp_L}{\sqrt{\log L}}\leq \tau=\tau_{12}\
;\ \tau_{13}\in \Big(\tau_{12},\tau_{12}+\frac{\vp_L}{(\log
L)^2}\Big]\bigg]\\ & & + \prob_{A_L}\bigg[\frac{\vp_L}{\sqrt{\log
L}}\leq \tau=\tau_{12}\ ;\ \tau_{23}\in
\Big(\tau_{12},\tau_{12}+\frac{\vp_L}{(\log L)^2}\Big]\bigg].
\end{eqnarray*}
Applying the strong Markov property at time $\tau_{12}$ and using
Lemma \ref{lemm no entrance} with $(\ell^L(t))_{t\geq 0}\equiv
(\psi_L^{-1}\{\xi_i^L-\xi_3^L\}((\psi_L^2\wedge \rho_L)t))_{t\geq 0}$ for
each $i\in \{1,2\}$ (as in the proof of Theorem \ref{theo time
coal}), we can conclude that each of the above terms tends to $0$
uniformly in $A_L$. On the event described by (\ref{lineages far
4}), that is $(1,2)$ is the first pair to meet and merge, no other pair meets in between but the distance between lineages $1$ and $3$ at time $\tau^*$ is smaller than $L/\log L$, the differences $\{\xi_1^L-\xi_2^L\}$ and $\{\xi_1^L-\xi_3^L\}$
have the same distribution as two independent copies $\hat{\xi}^L$
and $\check{\xi}^L$ of the process $\xi^L$ run at speed 2 up until $\tau$, and so if
we write $\hat{T}_L$ (resp., $\check{T}_L$) for the entrance time of
$\hat{\xi}^L$ (resp., $\check{\xi}^L$) into $B(0,2R^B\psi_L)$, with
a slight abuse of notation for the initial value to simplify the notation, (\ref{lineages far 4}) is bounded by
$$
\prob_{A_L}\Big[ \check{T}_L>\hat{T}_L\geq \frac{\vp_L}{\sqrt{\log
L}} ; \big|\check{\xi}^L\big(\hat{T}_L\big)\big|\leq \frac{L}{\log
L}\Big]\leq \prob_{A_L}\Big[\hat{T}_L\geq \frac{\vp_L}{\sqrt{\log
L}} ; \big|\check{\xi}^L\big(\hat{T}_L\big)\big|\leq \frac{L}{\log
L}\Big].
$$
A straightforward application of Lemma \ref{lemma local TCL} $(b)$
with $(\ell^L(t))_{t\geq 0}\equiv (\psi_L^{-1}\check{\xi}^L((\rho_L\wedge
\psi_L^2)t))_{t\geq 0}$ yields the uniform convergence of the last term to
$0$. Finally, the second term in (\ref{lineages far 1}), i.e., the probability that $(1,2)$ is the first pair to meet but not to merge, is bounded
by the sum over all pairs $\{i,j\}\in \{1,\ldots,4\}^2$ such that
$i\neq j$ and $\{i,j\}\neq \{1,2\}$ of
\begin{eqnarray*}
\prob_{A_L}\Big[\tau \geq\frac{\vp_L}{\sqrt{\log L}} ;
\tau^*=\tau^*_{12} ; \tau =& \tau_{ij}&\Big] \leq
\prob_{A_L}\Big[\frac{\vp_L}{\sqrt{\log L}}\leq \tau= \tau_{ij}
<\tau^*_{12}-\frac{\vp_L}{(\log L)^2} ;
\tau^*_{ij}>\tau^*_{12}\Big]\\& &+ \prob_{A_L}\Big[
\tau^*=\tau^*_{12}\geq \frac{\vp_L}{\sqrt{\log L}} ; \tau_{ij}\in
\Big[\tau^*_{12}-\frac{\vp_L}{(\log L)^2},\tau^*_{12}\Big] \Big].
\end{eqnarray*}
We can now conclude as we did for (\ref{lineages far 2}) and
(\ref{lineages far 3}).

When $\rho_L\gg \psi_L^2\log L$, we saw in the proof of Theorem
\ref{theo time coal} that with probability increasing to $1$, a pair
of lineages will not be affected by a large event during the periods
of time when the lineages are at distance less than $2R^B\psi_L$
from each other, until they come at distance less than $2R^s$.
Consequently, we could consider the evolution of the lineages to be
independent until their gathering time at distance $2R^s$. Because
we are still considering a finite number of lineages, the arguments
we used are applicable here again, and the proof of the last
paragraph also yields (\ref{lineage distribution 1}) in this case.
The proof of (\ref{lineage distribution 2}) is analogous, and is
therefore omitted. $\hfill \square$

\bigskip
\noindent \textbf{Acknowledgements.} The authors would like to thank the referees for their very careful reading and useful comments which helped to improve the presentation and to correct some inaccuracies. A. V\'eber would like to thank the Department of Statistics of the University of Oxford for hospitality.

\end{document}